\documentclass[12pt,a4paper]{amsart}
\usepackage{amsmath,amscd,amssymb,amsthm,array}
\usepackage{bbm}

\usepackage{mathtools}

\usepackage{hyperref}

\usepackage{comments}
\newComments\DL{D}{red}
\newComments\Bj{Bj}{red}
\newComments\IS{IS}{magenta}
\newComments\AO{AO}{blue}

\newcommand {\lra} {{\leftrightarrow}}

\usepackage{DLdef1}
\providecommand{\headspaces}{\spaceskip=.33em plus .05em minus .1em \xspaceskip=.5em plus .05em}
\makeatletter
\def\subsection{\@startsection{subsection}{2}{0pt}{-1ex \@plus -.4ex \@minus -.2ex}%
  {-.7em plus -.3em}{\normalfont\normalsize\bfseries\boldmath\headspaces}}
\def\subsubsection{\@startsection{subsubsection}{3}{0pt}{-1ex \@plus -.4ex \@minus -.2ex}%
  {-.7em plus -.3em}{\normalfont\normalsize\bfseries\boldmath\headspaces}}
\def\paragraph{\@startsection{paragraph}{3}{0pt}{-1ex \@plus -.4ex \@minus -.2ex}%
  {-.7em plus -.3em}{\normalfont\normalsize\bfseries\boldmath\headspaces}}
\makeatother

\makeatletter
\theoremstyle{thm}
\@ifundefined{Theorem}{}{}
\@ifundefined{Statement}{}{}
\@ifundefined{Corollary}{}{}
\@ifundefined{Lemma}{}{}
\@ifundefined{Question}{}{}
\theoremstyle{definition}
\@ifundefined{Remark}{}{}
\makeatother

\let\ssec\subsection
\def\tocstar{*}
\makeatletter\def\tocrunin{\if@noskipsec\leavevmode\fi}\makeatother

\renewcommand {\ssbegin}[2][*]
 {\tocrunin\refstepcounter{subsection}%
\def\tocttl{#1}\ifx\tocttl\tocstar
\addcontentsline{toc}{subsection}{\thesubsection.\hskip 1pc #2}%
\else
\addcontentsline{toc}{subsection}{\thesubsection.\hskip 1pc #2. #1}%
\fi
 \def \secno {\gdef \secno {}{\ssecfont
\thesubsection.\hskip 2ex}%
 }%
 \begin{#2}}

\renewcommand {\sssbegin}[2][*]
 {\tocrunin\refstepcounter{subsubsection}
\def\tocttl{#1}\ifx\tocttl\tocstar
\addcontentsline{toc}{subsubsection}{\thesubsubsection.\hskip 1pc #2}%
\else
\addcontentsline{toc}{subsubsection}{\thesubsubsection.\hskip 1pc #2. #1}
\fi
 \def \secno {\gdef \secno {}{\ssecfont \thesubsubsection.\hskip 2ex}%
 }%
 \begin{#2}}

\renewcommand {\parbegin}[2][*]
 {\refstepcounter{paragraph}
\if#1*
\addcontentsline{toc}{paragraph}{\theparagraph.\hskip 1pc #2}%
\else
\addcontentsline{toc}{paragraph}{\theparagraph.\hskip 1pc #2. #1}
\fi
 \def \secno {\gdef \secno {}{\ssecfont \theparagraph.\hskip 2ex}%
 }%
 \begin{#2}}

\newcommand {\grd}{{\mathop{\mathrm{grad}}\nolimits}}
\newcommand {\dvg}{{\mathop{\mathrm{div}}\nolimits}}

\newcommand {\Dx}{D_\xi}
\newcommand {\De}{D_\eta}
\newcommand {\Dq}{\cD_q}
\newcommand {\Dxi}{\cD_\xi}
\newcommand {\tp}{\otimes}
\newcommand {\Bol}{\mathrm{Bol}}

\newcommand {\Nord}{\operatorname{ord}}   
\newcommand {\etb}{\overline{\eta}}
\newcommand {\Jop}{J}

\title[GRC brackets on $1|n$-dimensional superstrings]{Gordan-Rankin-Cohen  brackets on $1|n$-dimensional superstrings with a contact structure for $n\leq 6$}

\author{Sofiane Bouarroudj${}^{a, *}$, Dimitry Leites${}^{b}$,\\ Oleksandr Lozhechnyk${}^c$, Irina Shchepochkina${}^{d}$}
\address{${}^*$Corresponding author; ${}^a$Division of Science and Mathematics,
New York University Abu Dhabi,
Po Box 129188, United Arab Emirates; sofiane.bouarroudj@nyu.edu\newline
${}^b$Department of Mathematics, Albanov\"agen 28,
Stockholm University, Stockholm, SE-106 91 Sweden; dimleites@gmail.com\newline
${}^c$Faculty of Mechanics and Mathematics, Taras Shevchenko National University of
Kyiv, Hlushkova Avenue 4e, 03127 Kyiv, Ukraine; alozhechnik@gmail.com\newline
${}^d$Independent University of Moscow,
B. Vlasievsky per., d. 11, RU-119 002 Moscow, Russia; irina@mccme.ru}

\begin{document}

\begin{abstract}The binary projectively $\mathfrak{sl}(2)$-invariant differential operators between spaces of weighted densities on the 1-dimensional manifold are called Gordan transvectants or Rankin-Cohen (GRC) brackets. For the two contact series (with an even and an odd ``time'' variable) on the supercircle or superstring of dimension $1|2$, we describe analogs of GRC brackets. 
We have discovered many new operators. For $\mathfrak{osp}(2|2)\subset\mathfrak{k}(1|2)$, the space of invariant binary operators of every given order is $2$-dimensional for \textit{all} weights of the densities, not only for the weights considered so far. For $\mathfrak{pe}(2)\subset\mathfrak{m}(1)$, there are three independent binary operators of order $2$, not two. While classifying them, we correct an earlier result on the unary operators: the order-$2$ Bol operator over $\mathfrak{m}(1)$ is only one member of a $1$-parameter family invariant for \textit{all} weights, which corrects Theorem 5.1 of arXiv:2110.10504. 

We correct earlier findings in the $1|3$-dimensional case and consider $\mathfrak{osp}_a(4|2)$-invariant (GRC) brackets in $1|4$-dimensional case, where the order-$2$ bracket is \textit{independent} of the parameter $a$ and is invariant under the whole of $\mathfrak{k}(1|4)$. 
In the $1|5$-dimensional case, 
the analogs of Bol operators are odd, of order $5+2k$, where $k=0,1,2,\dots$, and so is the first GRC bracket of odd order, which is of order $5$.
We also settle the exceptional case $\mathfrak{as}\subset\mathfrak{kas}\subset\mathfrak{k}(1|6)$. \end{abstract}


 \makeatletter
\@namedef{subjclassname@2020}{\textup{2020} Mathematics Subject Classification}
\makeatother


\maketitle

\markboth{{\itshape Sofiane Bouarroudj\textup{,} Dimitry Leites\textup{,} Oleksandr Lozhechnyk\textup{,} Irina Shchepochkina}}
{{\itshape Super Gordan-Rankin-Cohen  operators on $1|n$-dimensional superstrings}} 


\thispagestyle{empty}
\setcounter{tocdepth}{3}
\tableofcontents

\section{Introduction} \label{Sgcr}
Hereafter the ground field is $\Cee$, and $\langle\dots\rangle$ denotes the linear span; a~subalgebra generated by a~set is called so in words. For preliminaries, see \cite{BLS1}; here, we correct several statements there and in \cite{BLO}.

\subsection{Feasible classification problems of invariant differential operators}{}~{}

\underline{The unary invariant operators}. The projectively invariant differential operators between spaces of tensor fields on a (super)manifold are called \textit{Bol operators} \cite{Bol} or just \textit{Bols}. The problem of their description was completely solved on superstrings  --- supermanifolds of dimension $1|n$ --- endowed with a contact structure, and only for the spaces of weighted densities, not for general tensor fields. Let $\fg$ be the Lie superalgebra of all vector fields, or its subalgebra preserving a~contact structure, and let $\fg=\bigoplus\fg_i$ be its standard $\Zee$-grading. The Bols of \cite{BLS1} are invariant with respect to a~simple maximal subalgebra $\fs\subset\fg$ containing the negative part $\fg_-=\bigoplus_{i<0}\fg_i$. 

Observe that on contact supermanifolds in dimensions $1|n$ for $n=0, 1, 2$ there are no other \textit{irreducible} tensor fields except for weighted densities, since $\fo(n)$ is commutative for these $n$, so its irreducible finite-dimensional modules are $1$-dimensional (without a~contact structure this fails already for $n=1$), and the amount of possible analogs of generalized Bols between spaces tensor fields on supermanifolds of other superdimensions  remained unclear. 

In \cite{BLS2}, we showed that only cases of spaces of tensor fields on supermanifolds of ``key superdimensions\rq\rq\  ($2|0$, $0|n$ and $1|1$) are of practical interest. For the supermanifolds of other superdimensions $a|b$, the number of invariant operators is indescribably voluminous and depends on $a+b-1$ parameters even for the 1st order  operators. In contrast, on supermanifolds of these other superdimensions $a|b$, there are no non-scalar invariant differential operators between spaces of weighted densities. 

\underline{The binary invariant operators}. Consider binary operators invariant with respect to the maximal subalgebra containing all negative elements of the ambient superalgebra (of all vector fields, or of its subalgebra preserving a contact structure). For them, the answer is richer than for unary operators. (Without the requirement of invariance relative a~\textit{maximal} subalgebra, the binary operators invariant relative the whole of $\fvect(n)$, resp. $\fh(2n)$, were classified in \cite{Gr1}, resp. \cite{Gr2}, for supermanifolds, see \cite{GLS}; on the line, the anti-symmetric $m$-ary ones, for every $m$, were described in \cite{FF}.  Here we consider the binary operators called the \textit{Gordan transvectants} a.k.a. \textit{Rankin-Cohen  brackets} --- briefly: GRC operators in what follows. In this paper, the ambient superalgebra is always a~\textit{contact} superalgebra: $\fk(1|n)$ for $n\leq6$, and the pericontact $\fm(1)$; we make no claim about $\fvect(m|n)$, cf.~\S\ref{whatwedo}. For the answer on distinguished superstrings, see \cite{LShRC}.

\underline{The ternary invariant operators}.  On manifolds and between the spaces of weighted densities, they were classified in \cite{Bter}.

\ssec{Overview of the literature}
The GRC operators were introduced by Gordan \cite{Gor}, see the review by Gi\`eres~\cite{G}. They were rediscovered by Rankin \cite{R} and Cohen \cite{Co} and gained a great deal of currency as \textit{Rankin-Cohen operators}, see~\cite{Z} and the review \cite{El}. For their representation-theoretic counterparts and relatives, see \cite{Y, N, BCK}.

\sssec{The GRC operators between spaces of weighted densities when there are no other types of tensor fields}  As noted in \cite{BL}, the situation when all types of tensor fields are only weighted densities is only possible on superstrings of superdimension $1|n$ for $n=0$, 1, 2. Moreover, a contact structure must be present, given by an odd contact or an even pericontact form.

The first to describe super GRC operators --- bilinear operators invariant under the image of ${\fosp(1|2)}$ embedded into ${\fk(1|1)}$, i.e., on the $(1|1)$-dimensional superstring with a contact structure --- were Gi\`eres and Theisen, see \cite{GT}, see also \cite{H, GMO}. Next, a~slightly modified expression was given by Gargoubi and Ovsienko \cite{GO} who interpreted the super versions of the GRC operators they found in terms of the $(1|2)$-dimensional  superstring $\cC^{1|2}$  endowed with a~contact form.

In \cite{Om}, Omri  described the GRC  operators  invariant under $\fosp(2|2)$ embedded into $\fk(1|2)$, but he considered weighted densities with only 1 parameter, whereas there are 2 parameters, see \cite{BLS1}, so the result is incomplete.

The paper \cite{BLO} is similar to  \cite{Om}, but deals with $\fosp(3|2)$ embedded into $\fk(1|3)$; in this case spaces of weighted densities depend on~ 1 parameter; this is a complete result. For the \textit{unary} $\fosp(n|2)$-invariant operators between spaces of weighted densities, see \cite{BaL}.

\sssec{The GRC operators between spaces of weighted densities when there are other tensor fields}  For the list of maximal subalgebras $\fm\subset\fg$ such that $\fm_-=\fg_-$ in the $\Zee$-graded Lie (super)algebra of vector fields $\fg=\fvect(m|n)$ or its contact subalgebra, see \cite[Table~(4)]{BLS1} and \cite{LShch}. In the non-super case, there are two series of such subalgebras: $\fsl(n+1)$ and $\fo(n+2)$ in $\fvect(n)$.

For the maximal subalgebra $\fsl(n+1)\subset\fvect(n)$, Lecomte and Ovsienko \cite{LO} considered the projectively equivariant symbol calculus on the symmetric contravariant tensor fields, see also \cite{DO1}; their result was extended to these tensor fields tensored with a~weighted density in \cite{Bou}. The analogous problem for the conformal subalgebra $\fo(a+1,b+1)\subset\fvect_\Ree(n)$ and these tensor fields is open. 

Ovsienko and Redou (\cite{OR}, see also \cite{BonLe, DGO, DO2}) considered analogs of GRC brackets between spaces of weighted densities on the $n$-dimensional \textbf{real} manifold with a non-degenerate metric of signature $(a,b)$. They classified GRC brackets over $\Ree$ conformally invariant with respect to the maximal subalgebra $\fo(a+1, b+1)\subset\fvect_\Ree(n)$.

\sssec{Some applications of the super GRC brackets} 1) For the first appearances, see \cite{G, GT, H, CMZ}.

2) These brackets are used in \cite{BLO2} to compute $H^2(\fk(1|2), \fosp(2|2); {\cal D}^2(\cal F_\lambda, \cal F_\mu))$.  Similar computations were performed by Omri, who computed $H^1(\fk(1|2), \fosp(2|2); {\cal D}^2(\cal F_\lambda, \cal F_\mu))$, see~\cite[Thm.~4.1]{Om}. Here, ${\cal D}^2(\cal F_\lambda, \cal F_\mu)$ is the space of all linear differential operators acting between the spaces of weighted densities on~$\Ree^{1|2}$. The superscript~$2$ is the number of odd indeterminates, not the order. (The goal was to superize the result of \cite{BO} where the cohomology $H^1(\fvect(1), \fsl(2); {\cal D}(\cal F_\lambda, \cal F_\mu))$, where ${\cal D}$ is the space of (all) differential operators acting on the spaces of weighted densities was computed. This cohomology appears in the context of \textit{equivariant quantization}, cf. \cite{Me}; for the cohomology of the Lie algebras of vector fields on $\mathbb{RP}^1$ with coefficients in the bilinear differential operators, see \cite{Bc}.)

Next, Omri  described the deformation of the space of symbols; this deformation becomes trivial once restricted to $\fosp(2|2)$. To describe it, he used a~kind of Massey product of cocycles he had found previously. He also gave a~necessary and sufficient condition for an infinitesimal deformation to be a~global deformation.

Ben Ammar with his student Mtaouaa  computed the deformation that does not necessarily become trivial once restricted on $\fosp(2|2)$, see \cite{BeM}. They also obtained some necessary and sufficient conditions for that similar to what Omri did in \cite{Om}.


\ssec{What is done in this paper}\label{whatwedo}
We consider six cases, cf. \cite[Table~(4)]{BLS1}\textup{;} the first five are contact, the sixth is pericontact:

1) $\fosp(n|2)\subset\fk(1|n)$ for $n=1,2$, see \S\ref{ospINk}. For $n=2$ we classify
\textit{all} binary $i(\fosp(2|2))$-invariant operators, for \textit{all} weights of the
densities --- both parameters of the weight being free, which is what \cite{Om} lacked, see
Theorem~\ref{class}. There are exactly two operators of each order $n\geq1$, so the lists in
\cite{GO, Om} are incomplete: for every order there is one more even operator, and there are
odd operators of every odd order, not only of order~1.

2) $\fosp(3|2)\subset\fk(1|3)$, see \S\ref{Sblo}. This case is treated in \cite{BLO}; we
verify their answer and correct three statements of their Theorem 4.1 --- the signs in the
even supertransvectants, formula (4.2) of \cite{BLO}, the missing Levi-Civita sign in the odd ones, formula (4.1) of \cite{BLO}, and the
description of the exceptional weights.

3) $\fosp_a(4|2)\subset\fk(1|4)$, see \S\ref{Sn4}. Here, the subalgebra is not one but
a~one-parameter family. We describe the GRC brackets of every member in orders $\leq4$ (there are none in orders $1$ and $3$). We describe its Bol operators of orders $\leq8$, thereby confirming \cite[Thm~1.6]{BLS1} in these orders by an
independent computation. We also answer the question asked in \cite[\S3.2.3a]{BLS1} about
isomorphic algebras with different answers.

4) $\fosp(5|2)\subset\fk(1|5)$, see \S\ref{Sn5}. Here, on the contrary, the subalgebra $\fosp(5|2)$ is
rigid (the family occurs for $n=4$ and for no other $n$). The odd $n$ shows itself in
the parity: the Bol operators are \textit{odd}, of order $5+2k$ for every $k\in\Zee_+$, namely $T_5\circ\partial_t^{\,k}$, see Theorem~\ref{n5bol}, and so are the GRC brackets of
odd order. The first of these brackets is of order~$5$, not of order~$3$ as for $n=3$. For generic
weights, there is one bracket in each even order and one in each odd order from $5$ on.

5) $\fas\subset\fk\fas\subset\fk(1|6)$, see \S\ref{GRCAS}. This is a~ case
\cite[\S1.3.1]{BLS1} left ``for the future\rq\rq. It is unlike the others: $\fas$ consists of the
contact fields $K_f$ with $\partial_t(f)=0$, so the weight of the densities is invisible to
it. There are no exceptional weights and no analogs of Bol operators: the only invariant unary
operators are the powers of $\partial_t$. Every invariant binary operator of order $\leq4$ is a~linear combination of the operators $\partial_t^a\tp\partial_t^b$ and of the single non-trivial operator $\Omega$ of order~$2$, precomposed with $\partial_t$ in one of the arguments, see Theorem~\ref{kasbin}.  Since
the densities thus carry nothing, we go on to the tensor fields of the other types, see
\S\ref{kastensor}. There, the weights, and the twists by a~power of the determinant, are
invisible to $\fas$ as well. For the half-spinors, \textit{the two chiralities behave
differently}: exactly one of the two Dirac operators is $\fas$-invariant.

6) $\fpe(2)\subset\fm(1)$, see \S\ref{CoRm}. Here, the operators $\cJ_{k,l}$ and $\cJ_{k,r}$ are
the Bol operators of \cite[Thm~5.1]{BLS1} acting in one of the two arguments. The weights on
which they are invariant are computed in Theorem~\ref{m1thm}. For $k=2$, there is a~third,
independent, invariant operator, see Statement~\ref{m1class}. As a~by-product, we correct
\cite[Thm~5.1]{BLS1} itself at $k=2$: there, the order-$2$ Bol operator is not isolated but
a~member of the $1$-parameter family ${\Dxi\Dq+2a\,\partial_\tau:\cF_{a,b}\tto\cF_{a-1,b}}$,
invariant for \textit{all} weights. This family exhausts the invariant operators of its
degree, charge and parity, whereas only the member $a=\nfrac12$ is printed in \cite{BLS1}. For
$k=1,3,4$, the theorem is confirmed as printed, see Statement~\ref{BLSk2}. The complete classification of the GRC brackets over $\fpe(2)\subset\fm(1)$, for every order and every weight, is given in \S\ref{sec:invariant-binary}. 

We make no claim about the general vectorial series $\fvect(m|n)$: nothing in this text
concerns them. One further contact case, $\fsl(1|m|n|1)\subset\fk(2m+1|2n)$, is postponed.

\textbf{``Statements\rq\rq\ are proved by computations. Where the computation covers a~whole finite list, the Statement is final. Where it
covers several first cases, these cases give us reasons to formulate a~Conjecture for the
general case, stated right after the Statement that supports it.} All the classification statements below were
obtained, and all the formulas printed below were checked, by computation with all the
generators of the relevant subalgebra.

\section{Super GRC   operators: the $\fosp(n|2)\subset \fk(1|n)$ case}\label{ospINk}
In this Section, $\cF$ is the space of functions on the superstring, i.e., the coefficients of the weighted densities, and $f^{(a)}:=\left(\nfrac{\partial}{\partial t}\right)^a(f)$ for every $f\in\cF$.  Theorem \ref{Peetre} considerably simplifies calculations; it needs the following notation.

Recall the definition of $K_{f}$ and introduce the operators $D_{f}$, see \cite{BLS1}, for $\deg K_f<0$:
\be\label{KD}
\begin{array}{l}
\text{$D_{1}=\partial_t$,\ \ $K_1=2\partial_t$};\\
\text{$D_{\eta_i}:=\xi_i \partial_t+\partial_{\eta_i}$,\ $D_{\xi_i}:=\eta_i \partial_t+\partial_{\xi_i}$, \ $D_{\theta}:=\theta\partial_t+\partial_{\theta}$};  
\\[2mm]
\text{$K_{\xi_i}=\xi_i \partial_t - \partial_{\eta_i} $,\ $K_{\eta_i}= \eta_i \partial_t-\partial_{\xi_i}$, \ $K_{\theta}= \theta\partial_t - \partial_{\theta}$.} 
\end{array}
\ee

Recall the definition of $M_{f}$ and introduce the operators $\cD_{f}$ for $\deg M_f<0$:
\be\label{MD}
\begin{array}{l}
\text{$\cD_{1}=\partial_\tau$,\ \ $M_1=2\partial_\tau$};\\
\text{$M_{\xi_i}:=\xi_i \partial_\tau+\partial_{q_i}$,\ $M_{q_i}:=q_i \partial_\tau-\partial_{\xi_i}$;}\\[2mm]
\text{$\cD_{\xi_i}:=\xi_i \partial_\tau-\partial_{q_i}$,\ $\cD_{q_i}:=q_i \partial_\tau+\partial_{\xi_i}$}.
\end{array}
\ee

In the description of GRC operators, we will need, for \textbf{odd} operators $X$ and $Y$, the expression denoted in \cite[eq.~(42)]{BLS1} by $[X,Y]_+$ \textup{(}note that the supercommutator of odd operators is $XY+YX$, not $XY-YX$\textup{)}:
\be\label{comm}
[X, Y]_+:=XY-YX.
\ee

\ssbegin[All invariant differential operators have constant coefficients in an invariant  basis \cite{Sh}]{Theorem}[All invariant differential operators  have constant coefficients in an invariant   basis \cite{Sh}]\label{Peetre}
 Let $\fh$ be a~ subalgebra  in a~ vectorial Lie superalgebra $\fg$, and $\fh_-=\fg_-$.

 Let $\fg:=\fk(2m+1|n)$. Any $\fh_-$-invariant differential operator  is a~polynomial with constant coefficients in the operators~$D_f$, see  eq.~\eqref{KD}.

If $\fg=\fm(n)$, then any $\fh_-$-invariant differential operator is a~polynomial with constant coefficients in the operators $\cD_f$, see  eq.~\eqref{MD}.

If $\fg=\fvect(m|n)$, then any $\fh_-$-invariant differential operator  is a~polynomial with constant coefficients in partial derivatives.
\end{Theorem}

\textbf{Warning}. The expressions below look different from those given in \cite{GO} and \cite{Om} because we use different normal shapes of the contact form and bracket.

\sssec{The same on superstrings, with two exceptions occurring for $n=4$}\label{invOpStr}
Theorem~\ref{Peetre} is stated for polynomial coefficients, whereas on a~superstring --- a~supermanifold of superdimension $(1|n)$ --- the coefficients are Laurent polynomials. There are two superstrings to consider. The first superstring is associated with the trivial bundle of rank~$n$; we denote its contact Lie superalgebra by $\fk^{L}(1|n)$.\footnote{Speaking about~$\fk^L(1|n)$ \textup{(}resp.,~$\fk^M(1|n)$\textup{)}, physicists say ``Neveu--Schwarz\rq\rq\ \textup{(}resp., ``Ramond\rq\rq\textup{)} sector, although these family names designate the nontrivial central extensions of these Lie superalgebras, and a~sector, in physics, is the space of states on which such an extension acts. Even the Lie superalgebra~$\fk(1|n)$ with polynomial coefficients is sometimes attributed to the Neveu--Schwarz sector: the inclusion of polynomials into Laurent polynomials makes it a~subalgebra of~$\fk^L(1|n)$, and its subalgebra~$\fosp(n|2)$ is the ``global\rq\rq\ subalgebra of the Neveu--Schwarz superalgebra. For~$n\geq2$, the names no longer match the superscripts: the~$n=2$ Ramond superalgebra of physicists is the Neveu--Schwarz one with all indices integer, isomorphic to it by the spectral flow, i.e., a~central extension of~$\fk^L(1|2)$, whereas~$\fk^M(1|2)$ is what physicists call the ``twisted\rq\rq\ sector.} The second superstring is associated with the Whitney sum of the trivial bundle of rank $n-1$ and the M\"obius bundle; we denote its contact Lie superalgebra by $\fk^{M}(1|n)$. For the definitions, see \cite{LShRC}. In the standard grading, $\fk^{L}(1|n)_0=\fo(n)\oplus\Cee z$, whereas $\fk^{M}(1|n)_0=\fo(n-1)\oplus\Cee z$ for $n\leq4$. Since Theorem~\ref{Peetre} is the tool on which every ansatz of this paper rests, let us record what it becomes there.

Here, as in \cite{BLS1}, an operator $B\colon\cF_a\tto\cF_b$ between spaces of weighted densities is described by the operator $\overline B$ acting on the functions:
\be\label{barB}
\text{$B(f\alpha^a)=\overline B(f)\alpha^b$, where $f\in\cF_0$,}
\ee
where $\alpha$ designates the generator of a~rank-1 module of weighted densities over the algebra of functions; recall that the weight might be a~pair of numbers, see \cite{BLS1}.

Passing to Laurent coefficients changes nothing. Consider the stringy Lie superalgebra preserving a~structure --- a~tensor or a~distribution --- on the superstring considered. All the operators $\overline B$ invariant with respect to it are polynomials with constant coefficients in an invariant basis, that of $Y$-type vectors of \cite{Sh}, precisely as in Theorem~\ref{Peetre}. There are exactly two exceptions, and both occur for $n=4$: for $\fk^{M}(1|4)$ and for the derived algebra $\fk^{L}{}'(1|4)$, one extra generator is needed. It is written as a~product of a~negative power of a~$Y$-vector and a~product of $Y$-vectors, see formulas \eqref{DeltaIS} and \eqref{DeltaM}. This product is \textbf{not} the composition of its factors, and it is \textbf{not} a~differential operator. 
This is essential: it has to be simplified before it is evaluated, i.e., treated as one whole.

For $\fk^{M}(1|n)$, let $\zeta$ be the odd coordinate of the M\"obius summand and $\theta_1,\dots,\theta_{n-1}$ those of the trivial summand. Introduce the operators
\be\label{DzetaM}
\begin{array}{llll}
\cD:=\pder{t} -\nfrac{\zeta}{2t}\pder{\zeta}=\nfrac
12 K^M_1;&&&\\
K_{\zeta}^M :=\zeta\partial_{t} -\nfrac{1}{t}\partial_{\zeta},&D_{\zeta}^M :=
\zeta\partial_{t}
+\nfrac{1}{t}\partial_{\zeta}&\text{for $n\geq 1$};\\
K_{\theta_i}^M :=
\theta_i\partial_t-\partial_{\theta_i}-\frac{\theta_i\zeta}{2t}\partial_{\zeta},&D_{\theta_i}^M :=
\theta_i\partial_t+\partial_{\theta_i}-\frac{\theta_i \zeta}{2t}\partial_{\zeta}&\text{for $n\geq 2$}.\\
\end{array}
\ee
The two extra generators, introduced in \cite{LShRC} in order to write the cocycles of these two exceptional stringy superalgebras, are the following ones. For $\fk^{L}{}'(1|4)$, it is the operator $\Delta_{IS}$ defined by the table in formula~\eqref{DeltaIS}, where $D_1^{-1}=\partial_t^{-1}$ is the formal inverse of $D_1=\partial_t$, i.e., $t^a\mapsto\frac{1}{a+1}t^{a+1}$ for $a\neq-1$ 
\be\label{DeltaIS}
\Delta_{IS}:=D_1^{-1}D_{\theta_1} D_{\theta_2}D_{\theta_3}D_{\theta_4}:\quad
\begin{array}{rcl}
t^a&\longmapsto& a(a-1)(a-2)\,t^{a-3}\theta_1\theta_2\theta_3\theta_4,\\
t^a\theta_1&\longmapsto&-a(a-1)\,t^{a-2}\theta_2\theta_3\theta_4,\\
t^a\theta_1\theta_2&\longmapsto&-a\,t^{a-1}\theta_3\theta_4,\\
t^a\theta_1\theta_2\theta_3&\longmapsto&t^{a}\theta_4,\\
t^a\theta_1\theta_2\theta_3\theta_4&\longmapsto&\nfrac{1}{a+1}\,t^{a+1}
\end{array}
\ee
(and similarly for the other $\theta$-monomials), whose only singularity --- the case $a=-1$ in the last line --- is precisely the element deleted when passing from $\fk^{L}(1|4)$ to $\fk^{L}{}'(1|4)$. For $\fk^{M}(1|4)$, it is the operator $\Delta$ defined, for every $\vf\in\Lambda(\theta_1,\theta_2,\theta_3)$ and every $n$, by the two formulas
\be\label{DeltaM}
\Delta(t^n\vf(\theta))=t^n\zeta \vf(\theta), \qquad \Delta(t^n\zeta\vf(\theta))=\nfrac{t^n}{n+\nfrac 12}\vf(\theta);
\ee
we write $\Delta=D^M_\zeta\cD^{-1}$ for short,
an operator which is defined everywhere and coincides with the composition $D^M_\zeta\circ\cD^{-1}$ on every function on which the composition $D^M_\zeta\cD^{-1}$ is defined --- the pole of $\cD^{-1}$ at $n=-1$ being canceled by the zero of $D^M_\zeta$. In both cases, the literal reading of the composition --- perform one factor, then the other --- is \textit{not} the operator meant: it loses the cancellation. For $\Delta_{IS}$ and $g=t^3$, for instance, the four $D_{\theta_i}$ produce $a(a-1)(a-2)(a-3)t^{a-4}\theta_1\theta_2\theta_3\theta_4$ and $D_1^{-1}$ divides by $a-3$, so the literal reading returns $0$ at $a=3$, whereas eq.~\eqref{DeltaIS} returns $6\,\theta_1\theta_2\theta_3\theta_4$, which is what the cocycle identity requires. For $\Delta$ and $g=t^{-1}\theta_1\theta_2\theta_3$, an element of the algebra, the literal reading is not defined at all.

\ssbegin[Extending Theorem \ref{Peetre} to superstrings \cite{Sh, BLS1, LShRC}]{Theorem}[Extending Theorem~$\ref{Peetre}$ to superstrings; \cite{Sh, BLS1, LShRC}]\label{Ira}
The description of the $\fh_-$-invariant operators $\overline B$ given in Theorem~$\ref{Peetre}$ holds verbatim for the stringy versions of the Lie superalgebras listed there, i.e., for Laurent coefficients.

If $\fg=\fk^M(1|n)$, where $n\leq 3$, then any $\fg$-invariant differential operator $\overline B$ is a~polynomial with constant coefficients in the operators $D_f$ of eq.~\eqref{DzetaM}.

If $\fg=\fk^M(1|4)$, then any $\fg$-invariant operator $\overline B$ is a~polynomial with constant coefficients in the operators $D_f$ of eq.~\eqref{DzetaM} and the composition $\Delta D^{M}_{\theta_3}D^{M}_{\theta_2}D^{M}_{\theta_1}$, see eq.~\eqref{DeltaM}.

If $\fg=\fk^L{}'(1|4)$, then any $\fg$-invariant operator $\overline B$ is a~polynomial with constant coefficients in the operators $D_f$ of eq.~\eqref{KD} and the operator $\Delta_{IS}$, see eq.~\eqref{DeltaIS}.
\end{Theorem}

In the last three claims, $\fg$-invariance, and not $\fh_-$-invariance, is meant: for $\fk^{M}(1|n)$, the operators $D^M_f$ do not supercommute with the whole of the analog of $\fh_-$; e.g., $[\cD,D^M_\zeta]=-\frac{1}{2t}D^M_\zeta$. So the claims describe the shape of an operator invariant with respect to the whole of $\fg$; they say nothing about the operators invariant with respect to a~subalgebra.

\begin{proof} All the claims, except for the cases $\fk^M(1|n)$ and $\fk^L{}'(1|4)$, are proved in \cite{Sh, BLS1}. The exceptional cases are proved by the same arguments applied to the differential part of $\overline B$. It remains to observe that each of the two extra operators, $\Delta$ and $\Delta_{IS}$, is the only --- up to a~factor --- invariant extension of that differential part to the whole algebra. Note that the arguments of \cite{Sh, BLS1} concern \textit{differential} operators, whereas $\Delta$ and $\Delta_{IS}$ are not differential ones. \end{proof}

\sssec{What this means for the ansatz}\label{ansatzstr} For the cases treated in this paper, $n\leq3$, Theorem~\ref{Ira} says that nothing has to be added. The ansatz ``a~polynomial with constant coefficients in the $D_f$\rq\rq , used from Subsection~\ref{orderdeg} on, is complete on the superstring as well as on the supercircle. In the case $n=4$, the algebra $\fk^{L}(1|4)$ is not simple, and the Bol and GRC operators are not described. For $n=4$, the ansatz has to be enlarged by the corresponding operator~\eqref{DeltaIS} or~\eqref{DeltaM}. Since these are not differential operators, the very notion of the order of the operator has to be revised there. %

\ssec{The order and the degree of an operator}\label{orderdeg}
Both words are used in the literature we quote, and they denote two different integers. Here
is the convention of this paper, applied from now on without further comment.

By Theorem~\ref{Peetre}, every operator considered below is a~polynomial with constant
coefficients in the operators $D_f$ of eq.~\eqref{KD} \textup{(}resp. $\cD_f$ of
eq.~\eqref{MD}\textup{)}. For such a~polynomial:

1) Its \textit{order} $\Nord$ is its degree in the $D_f$, where each $D_{\xi_i}$,
$D_{\eta_i}$, $D_{\theta_i}$ \textup{(}resp. $\cD_{\xi_i}$, $\cD_{q_i}$\textup{)} counts
for~$1$, whereas $D_1=\partial_t$ \textup{(}resp. $\cD_1=\partial_\tau$\textup{)} counts
for~$2$; thus $\Nord(D_{\xi_1}D_{\xi_2}\partial_t)=4$.

2) Its \textit{degree} $s$ is its weight relative the grading $\deg t=2$,
$\deg\xi_i=\deg\eta_i=\deg\theta_i=1$, the grading for which $\fk(1|n)_0=\fc\fo(n)$.

Since $\deg\partial_t=-2$ whereas $\deg D_f=-1$ for odd $f$, these are one and the same
integer read with opposite signs:
\be\label{ords}
s=-\Nord .
\ee
We speak of the order everywhere, this being the parameter by which the ansatz is graded. We
use the degree only where the convention of \cite{BLS1} is needed literally. There, an operator
of degree $s$ shifts the weight $\lambda$ of the density by $-s$, and it is the degree that
enters the weight of the target together with the charge, see eq.~\eqref{target}.

\sssec{A~warning concerning the index $k$}\label{indexk} The index of a~bracket is a~subscript, not
a~superscript\textup{:} in $\cJ_2$, $\cJ_4$, $\cJ_{1/2,\xi}$, $\cJ_{k,l}$, the $2$ and the $4$ are labels,
not exponents\textup{;} the second subscript, $l$, $r$, $\xi$ or $\eta$, tells apart brackets of
one order. The weights $\lambda_1,\lambda_2\in\Cee$ \textup{(}and $\nu_1,\nu_2$, if any\textup{)} are those
of the arguments\textup{:} $\cJ_k\colon\cF_{\lambda_1}\tp\cF_{\lambda_2}\tto\cF_{\lambda_1+\lambda_2+\Nord}$.
In the brackets $\cJ_k$ and in the Bol operators $\Bol_k$, the letter $k$ is \textit{not always} the order, and it does not have the same
meaning in all Sections, since in each case we keep the notation of the source:

1) In \S\ref{ospINk}, after \cite{GO, Om}, the bracket $\cJ_k$ is of order $2k$, where
$k\in\nfrac12\Nee$, and it shifts the weight by $2k$.

2) In \S\ref{Sblo}, after \cite{BLO}, the supertransvectant of index $k$ is likewise of
order $2k$, but in the normalization \eqref{eq:conv} of \cite{BLO} it shifts the weight
by~$k$.

3) In \S\ref{CoRm}, after \cite{BLS1}, the Bol operator $\Bol_k$ and the brackets $\cJ_{k,l}$
and $\cJ_{k,r}$ built from it are of order~$k$.

Wherever the order is meant, it is written $\Nord$, or spelled out in words.

\ssec{Conventions for $\fk(1|2)$, and the two realizations of the odd
coordinates}\label{conv2}
We follow \cite{BLS1} throughout; we recall the formulas because everything below depends on
the normalization. The contact form is $\alpha_1=dt+\xi d\eta+\eta d\xi$, see
\cite[(18)]{BLS1}, and, with $E=\xi\partial_\xi+\eta\partial_\eta$,
\be\label{Kf}
K_f=(2-E)(f)\partial_t-H_f+\partial_t(f)E,\quad
\Dx=\eta\partial_t+\partial_\xi,\quad \De=\xi\partial_t+\partial_\eta,
\ee
where $H_f:=-(-1)^{p(f)}\big(\partial_\xi(f)\,\partial_\eta+\partial_\eta(f)\,\partial_\xi\big)$ is the Hamiltonian vector field of~$f$, see \cite[(24),(39)]{BLS1}, in accordance with eq.~\eqref{KD}. The spaces of weighted densities
are $\cF_{\lambda,\nu}=\cF\,\alpha_1^{\lambda/2}\beta^\nu$, where
$\beta=\sqrt{[d\xi(d\eta)^{-1}]}$, with
\be\label{Lie}
L_{K_f}(\varphi)=K_f(\varphi)+\big(\lambda\,\partial_t(f)
+(-1)^{p(f)}\nu\,\Div_2(f)\big)\varphi, \qquad
\Div_2=\widetilde K_\eta\widetilde K_\xi-\partial_t,
\ee
where
\[
\widetilde K_\xi(f)=(-1)^{p(f)}(\partial_\eta-\xi\partial_t)(f),\qquad
\widetilde K_\eta(f)=(-1)^{p(f)}(\partial_\xi-\eta\partial_t)(f),
\]
see \cite[(33),(35),\;2.2.1a]{BLS1}. With this normalization, an operator of degree $s$ shifts
$\lambda$ by $-s$, see \S\ref{orderdeg}, so the bracket $\cJ_k$ in $\sdim=1|1$, which is of
order $\Nord=2k$, maps
$\cF_{\lambda_1}\tp\cF_{\lambda_2}\tto\cF_{\lambda_1+\lambda_2+2k}$, as in
Theorem~\ref{GOthm}.

Set $\deg t=2$ and $\deg\xi=\deg\eta=1$, so that $\deg\partial_t=-2$ and
$\deg\Dx=\deg\De=-1$. Since $\fk(1|2)_0=\Cee K_t\oplus\Cee K_{\xi\eta}$ is
$2$-dimensional, the space of invariant operators is bigraded: besides the degree $s=-\Nord$
of \S\ref{orderdeg}, which is its weight relative $K_t$, every operator has a~\textit{charge},
its weight relative $K_{\xi\eta}=-\xi\partial_\xi+\eta\partial_\eta$, so that
\[
\ch\Dx=+1,\qquad \ch\De=-1,\qquad \ch\partial_t=0 .
\]
For an operator of degree $s$ and charge $w$ the weight of the target is forced:
\be\label{target}
\lambda=\lambda_1+\lambda_2-s,\qquad \nu=\nu_1+\nu_2-w .
\ee
We write $P\tp Q$ for the bilinear operator
$(P\tp Q)(\varphi,\psi)=(-1)^{p(Q)p(\varphi)}P(\varphi)\,Q(\psi)$; the signs
$(-1)^{p(\varphi)}$ written out in the formulas below are exactly this Koszul rule. Put
\be\label{PSL}
P:=\Dx\De-\partial_t=\nfrac12(\Dx\De-\De\Dx) .
\ee

\sssec{The dictionary $\theta\leftrightarrow\xi\eta$}\label{dictionary} Some of the formulas
below are quoted from \cite{GO, Om} in the realization in which the Gram matrix of the
bilinear form on the odd coordinates is diagonal, i.e., in terms of $\theta_1,\theta_2$. The
rest are given in the realization \eqref{Kf}, in which this matrix is anti-diagonal. Since
$\theta_1d\theta_1+\theta_2d\theta_2=\xi d\eta+\eta d\xi$ for
\be\label{dict}
\xi=\nfrac{\theta_1-i\theta_2}{\sqrt2},\qquad \eta=\nfrac{\theta_1+i\theta_2}{\sqrt2},
\qquad i=\sqrt{-1},
\ee
we have, with $D_{\theta_j}=\partial_{\theta_j}+\theta_j\partial_t$ and with $\nu^{\xi\eta}$, resp. $\nu^{\theta}$, denoting the second weight $\nu$ of the densities $\cF_{\lambda,\nu}$ in the $\xi\eta$, resp. $\theta$, coordinates, 
\be\label{dictD}
\Dx=\nfrac{D_{\theta_1}+iD_{\theta_2}}{\sqrt2},\qquad
\De=\nfrac{D_{\theta_1}-iD_{\theta_2}}{\sqrt2},\qquad
\xi\eta=i\,\theta_1\theta_2,\quad\text{hence}\quad \nu^{\xi\eta}=i\,\nu^{\theta} ,
\ee
and, as an identity with no stray factor,
\be\label{sym}
\begin{array}{l}
\Dx\tp\De+\De\tp\Dx=D_{\theta_1}\tp D_{\theta_1}+D_{\theta_2}\tp D_{\theta_2},\\
\Dx\tp\De-\De\tp\Dx=-i\,(D_{\theta_1}\tp D_{\theta_2}-D_{\theta_2}\tp D_{\theta_1}),
\end{array}
\ee
whereas $P=-i\,D_{\theta_1}D_{\theta_2}$ for $P$ as in eq.~\eqref{PSL}. Thus, the answer in the
$\xi\eta$ realization is the answer in the $\theta$ realization written in a~basis in which
not only $\ad_{K_t}$ but also $\ad_{K_{\xi\eta}}$ acts diagonally, i.e., a~basis consisting of eigenvectors of both operators. In that basis, the weights of the densities are
free complex parameters and no $\sqrt{-1}$ occurs; this is why we state the classification,
Theorem~\ref{class}, in the $\xi\eta$ realization.

\ssbegin[{The GRC brackets in $\sdim=1|1$; \cite{GO}}]{Theorem}[The GRC brackets in $\sdim=1|1$; \cite{GO}]\label{GOthm} Let $\dim \cM=1|1$ with coordinates $t$ (even) and $\theta$ (odd).
There are $i(\mathfrak{osp}(1|2))$-invariant bilinear differential operators labeled by $k\in\frac12\Nee$ and $a, b\in\Zee_+:=\{0, 1, 2, \dots\}$
\[
{\cal J}^k (\varphi, \psi): {\cal F}_{ \lambda_1} \otimes {\cal F}_{\lambda_2} \tto {\cal F}_{\lambda_1+\lambda_2+2k},
\]
more precisely:

\textup{(i)} If $k\in\Nee$, then ${\cal J}^k$ is even; it is given by the formula
\[
{\cal J}^k (\varphi,\psi)=\sum_{a+b= k -1} (-1)^{p(\varphi)} \gamma_{a,b,k-1} D_\theta( \varphi^{(a)})D_\theta (\psi^{(b)})- \sum_{a+b=k} \gamma_{a,b,k-1} \varphi^{(a)} \psi^{(b)},
\]
where $D_\theta=\partial_{\theta}+\theta \partial_t$, and \textup{(}this formula will be needed for half-integer values of $k$ as well in~\textup{(ii))}
\be\label{gammaABK}
\gamma_{a,b,k}:=(-1)^{a+1} \binom{\lambda_1+\lfloor k \rfloor}{b}\binom{\lambda_2+\lfloor k \rfloor}{a}.
\ee

\textup{(ii)} If $k\in\Nee-\frac12$, then ${\cal J}^k$ is odd; it is given by the formula

\[
{\cal J}^k (\varphi,\psi)=\sum_{a+b=\lfloor k \rfloor} \gamma_{a,b,k} \left ((-1)^{p(\varphi)} (\lambda_1+a)\varphi^{(a)}D_\theta (\psi^{(b)})-(\lambda_2+b)D_\theta (\varphi^{(a)})\psi^{(b)}\right ).
\]

There are no other non-zero non-scalar bilinear $i(\mathfrak{osp}(1|2))$-invariant differential operators.
\end{Theorem}

\ssbegin[{The GRC brackets in $\sdim=1|2$ found so far}]{Theorem}[The GRC brackets in $\sdim=1|2$ found so far]\label{Omthm} Let $\dim \cM=1|2$.
There are the following $i(\mathfrak{osp}(2|2))$-invariant bilinear differential operators
\[
{\cal J}^k (\varphi, \psi): {\cal F}_{ \lambda_1, \nu_1} \otimes {\cal F}_{\lambda_2, \nu_2} \tto {\cal F}_{\lambda_1+\lambda_2+2k, \omega}:
\]

\textup{(i)} \textup{(\cite[Thm.~3.1]{Om})} If $k\in\Nee$ and $a, b\in\Zee_+$, then ${\cal J}^k$ is even; it is given by the formula
\[
\begin{array}{lcl}
{\cal J}^k (\varphi,\psi) & = & \displaystyle \sum_{a+b=k -1} (-1)^{p(\varphi)} \gamma_{a,b,k-1} \left ( D_{\theta_1}( \varphi^{(a)})D_{\theta_1} (\psi^{(b)}) + D_{\theta_2}( \varphi^{(a)})D_{\theta_2} (\psi^{(b)}) \right )\\[3mm]
&& \displaystyle
- \sum_{a+b=k} \gamma_{a,b,k-1} \varphi^{(a)} \psi^{(b)}
+ \sum_{a+b=k -2} \gamma_{a,b,k-1} \left ( D_{\theta_1} D_{\theta_2}( \varphi^{(a)})D_{\theta_1} D_{\theta_2} (\psi^{(b)}) \right ),
\end{array}
\]
where $D_{\theta_s}=\partial_{\theta_s}+\theta_s\partial_t$ for $s=1,2$ and where $\nu_1=\nu_2=\omega=0$. \textup{(}For the reason why $\nu_1=\nu_2=\omega=0$ has to be assumed here, and for what happens when $\nu_1,\nu_2,\omega\neq0$, see Corollary~$\ref{nu}$: it is a~property of the ansatz used, not of the operators.\textup{)}

\textup{(ii)} There are two odd operators of order~$1$; each of them is given by the formula

\[
{\cal J}^k (\varphi,\psi)= \rho D_{\theta_1} (\varphi) \psi+ \gamma D_{\theta_2} (\varphi) \psi+(-1)^{p(\varphi)} (\sigma \varphi D_{\theta_1} (\psi) + \chi \varphi D_{\theta_2} (\psi)),
\]
where only the following two cases are possible and where $i=\sqrt{-1}$ \textup{(}for the same two operators written in the $\xi\eta$ realization, where no $\sqrt{-1}$ occurs, see Statement~$\ref{J12}$\textup{)}:
\[
\begin{array}{llllll}
1)&\omega=-i+\nu_1+\nu_2,& \sigma= \lambda_1-i \nu_1,& \rho= -\lambda_2+i \nu_2, & \chi=-i\lambda_1-\nu_1,& \gamma=i \lambda_2+\nu_2, \\[3mm] 
2)&\omega=i+\nu_1+\nu_2, &\sigma= i\lambda_1- \nu_1, &\rho= - i \lambda_2+ \nu_2, & \chi= - \lambda_1- i \nu_1, &\gamma=\lambda_2+ i \nu_2.
\end{array}
\]

\textup{(iii)} These are not all the operators: for every order there is exactly one more, see Theorem~$\ref{class}$.
\end{Theorem}

Part~(i) is Omri's Theorem~3.1 in~\cite{Om}. There the densities have one parameter, i.e., $\nu_1=\nu_2=\omega=0$, and for half-integer~$k$ Omri proves that every $\fosp(2|2)$-invariant bilinear operator vanishes. This agrees with part~(ii): the odd operators exist only for $\omega=\nu_1+\nu_2\pm i$, so never with $\nu_1=\nu_2=\omega=0$. They appear as soon as the second parameter of the weight is allowed.

This is the first time we encounter weights of densities belonging to $\Cee\setminus\Ree$.
This happens because the passage from the normal shape of the bilinear form with the diagonal Gram matrix to the shape with an anti-diagonal Gram matrix involves $\sqrt{-1}$. This is nothing spectacular, actually: in the $\xi\eta$ realization, see \S\ref{dictionary}, the weights are free complex parameters and nothing is imaginary. In that realization, the two cases 1) and 2) of part~(ii) become the two \textit{chiralities} of one formula.

\ssbegin[{The two odd operators of order~$1$ in the $\xi\eta$ realization}]{Statement}[The two odd operators of order~$1$ in the $\xi\eta$ realization]\label{J12}
For all $\lambda_1,\nu_1,\lambda_2,\nu_2$, the odd $i(\fosp(2|2))$-invariant bilinear
operators of order~$1$ form a~$2$-dimensional space, spanned by
\be\label{J12f}
\begin{array}{ll}
\cJ_{1/2,\xi}=\mu_2^{+}\,(\Dx\tp 1)-\mu_1^{+}\,(1\tp \Dx):&
\cF_{\lambda_1,\nu_1}\tp\cF_{\lambda_2,\nu_2}\tto
\cF_{\lambda_1+\lambda_2+1,\;\nu_1+\nu_2-1},\\[2mm]
\cJ_{1/2,\eta}=\mu_2^{-}\,(\De\tp 1)-\mu_1^{-}\,(1\tp \De):&
\cF_{\lambda_1,\nu_1}\tp\cF_{\lambda_2,\nu_2}\tto
\cF_{\lambda_1+\lambda_2+1,\;\nu_1+\nu_2+1},
\end{array}
\ee
where $\mu_i^{\pm}:=\lambda_i\pm\nu_i$; in longhand,
$\cJ_{1/2,\xi}(\varphi,\psi)=\mu_2^{+}\Dx(\varphi)\psi
-(-1)^{p(\varphi)}\mu_1^{+}\varphi \Dx(\psi)$, and similarly for $\eta$.
The operator $\cJ_{1/2,\xi}$ vanishes identically if and only if
$\mu_1^{+}=\mu_2^{+}=0$. If exactly one of the $\mu_i^{+}$ vanishes, it degenerates into the
composition of the multiplication with the odd Bol operator $B_{\xi}$ of
\cite[Thm~1.5(ii)]{BLS1} acting in the other argument.
\end{Statement}

\begin{proof} Direct computation with all eight generators of $i(\fosp(2|2))$. That formula~\eqref{J12f}
gives the same pair of operators as in part~(ii) of Theorem~\ref{Omthm} is seen from
formulas~\eqref{dict}--\eqref{dictD}: with $A_j:=\lambda_j+i\nu_j^{\theta}=\mu_j^{+}$, case 2) reads
$\rho=-iA_2$, $\gamma=A_2$, $\sigma=iA_1$, $\chi=-A_1$, that is,
\[
A_2\,(D_{\theta_2}-iD_{\theta_1})(\varphi)\psi
-(-1)^{p(\varphi)}A_1\varphi(D_{\theta_2}-iD_{\theta_1})(\psi),
\]
which is $\cJ_{1/2,\xi}$ up to the factor $-i\sqrt 2$, since $D_{\theta_2}-iD_{\theta_1}=-i\sqrt2\,\Dx$; and
$\omega^{\theta}=i+\nu_1^\theta+\nu_2^\theta$ becomes $\nu=\nu_1+\nu_2-1$ after multiplication
by~$i$. Case 1) is $\cJ_{1/2,\eta}$. Thus, the two ``cases\rq\rq\  are the eigenvalues $\mp1$ of the
charge, and the weights $\lambda\pm\nu$ replace $\lambda\pm i\nu$. \end{proof}

Part~(i) of Theorem~\ref{Omthm} for $k=1$ likewise takes a~transparent shape: by the identities~\eqref{sym}
and $\gamma_{1,0,0}=\lambda_2$, $\gamma_{0,1,0}=-\lambda_1$, $\gamma_{0,0,0}=-1$,
\be\label{J1}
-\cJ_1=\lambda_2(\partial_t\tp 1)-\lambda_1(1\tp\partial_t)+\Dx\tp\De+\De\tp\Dx ,
\ee
and this operator is indeed invariant --- for $\nu_1=\nu_2=\omega=0$ and all
$\lambda_1,\lambda_2$.

\ssec{The classification for $\fosp(2|2)\subset\fk(1|2)$}\label{classification}
The lists in \cite{GO, Om} and in Theorem~\ref{Omthm} are not complete. The complete answer is
as follows; observe that no restriction on the weights is imposed.

\ssbegin[The classification for $\fosp(2|2)\subset\fk(1|2)$]{Theorem}[The classification
for $\fosp(2|2)\subset\fk(1|2)$: two operators of every order, for all weights]\label{class}
Let $n\geq1$. For all weights $(\lambda_1,\nu_1)$ and $(\lambda_2,\nu_2)$, the space of
$i(\fosp(2|2))$-invariant bilinear differential operators of degree $s=-n$, i.e., of
order~$n$,
\[
\cF_{\lambda_1,\nu_1}\tp\cF_{\lambda_2,\nu_2}\tto\cF_{\lambda,\nu}
\]
is $2$-dimensional; the two operators are interchanged by the \textit{chirality flip}
$\xi\leftrightarrow\eta$, i.e., $\nu_i\mapsto-\nu_i$, $\Dx\leftrightarrow\De$,
$P\mapsto-P$. For $n$ odd, their charges are $w=+1$ and $w=-1$, and both operators are odd;
for $n$ even both have charge $w=0$, both operators are even, and there is nothing of charge
$\pm2$. The weight of the target is given by formula~\eqref{target}. For $n=0$, the only invariant
operator is the multiplication.
\end{Theorem}

The operators of Theorem~\ref{class} can be written explicitly. Set
\[
T:=\partial_t,\qquad E_j:=T^j,\qquad F_j:=PT^{j-1}\ \ (j\geq1),\qquad
X_j:=\Dx T^j,\qquad Y_j:=\De T^j,
\]
with $P$ as in formula~\eqref{PSL}, and $\Lambda_i^{\pm}:=\lambda_i\pm\nu_i$ for $i=1,2$.

\ssbegin[{The operators of Theorem~\ref{class}, explicitly}]{Statement}[The operators of Theorem~$\ref{class}$, explicitly]\label{Bjexplicit} 
Let the weights be generic, $q\in\Zee_+$, and let $\cJ$ be an $i(\fosp(2|2))$-invariant bilinear operator $\cF_{\lambda_1,\nu_1}\tp\cF_{\lambda_2,\nu_2}\tto\cF_{\lambda,\nu}$.

\textup{(i)} Let the order be $\Nord=2q+1$. For the charge $w=+1$, the operator is odd and, up to a~factor, it is
\be\label{Bjodd}
\cJ_{q+1/2,\xi}=\sum_{0\leq j\leq q}\pi_j\Big((\Lambda_2^{+}+q-j)\,X_j\tp E_{q-j}-(\Lambda_1^{+}+j)\,E_j\tp X_{q-j}
-(q-j)\,X_j\tp F_{q-j}+j\,F_j\tp X_{q-j}\Big),
\ee
where the terms with $F_0$ are absent, the target is $\lambda=\lambda_1+\lambda_2+2q+1$, $\nu=\nu_1+\nu_2-1$, and
\[
\pi_j:=(-1)^j\binom{q}{j}\prod_{0\leq r\leq j-1}(\lambda_2+q-r)\prod_{j+1\leq r\leq q}(\lambda_1+r).
\]
The operator $\cJ_{q+1/2,\eta}$ of charge $w=-1$ is obtained from formula~\eqref{Bjodd} by the chirality flip $X_j\mapsto Y_j$, $\Lambda_i^+\mapsto\Lambda_i^-$, $F_j\mapsto-F_j$\textup{;} its target is $\nu=\nu_1+\nu_2+1$.

\textup{(ii)} Let the order be $\Nord=2q$, where $q\geq1$. The charge is $w=0$, the target is $\lambda=\lambda_1+\lambda_2+2q$, $\nu=\nu_1+\nu_2$, and the operator is even:
\be\label{Bjeven}
\begin{array}{r@{}l}
\cJ_q={}&\displaystyle \sum_{0\leq j\leq q}A_j\,E_j\tp E_{q-j}+\sum_{0\leq j\leq q-1}B_j\,E_j\tp F_{q-j}
+\sum_{1\leq j\leq q}C_j\,F_j\tp E_{q-j}\\[1mm]
&{} \displaystyle +\sum_{1\leq j\leq q-1}D_j\,F_j\tp F_{q-j}
+\sum_{0\leq j\leq q-1}\big(U_j\,X_j\tp Y_{q-j-1}+V_j\,Y_j\tp X_{q-j-1}\big),
\end{array}
\ee
where $U_0$ and $V_0$ are arbitrary, which is the $2$-dimensionality of Theorem~\ref{class}, and
\[
U_j=\rho_jU_0,\qquad V_j=\rho_jV_0,\qquad
\rho_j:=(-1)^j\binom{q-1}{j}\prod_{1\leq r\leq j}\frac{\lambda_2+q-r}{\lambda_1+r}.
\]
For $0\leq j\leq q-1$, the remaining coefficients are
\[
\begin{array}{r@{}l}
B_j={}&\dfrac{(\Lambda_1^{-}+j)V_j-(\Lambda_1^{+}+j)U_j}{2(\lambda_2+q-j)},\\[3mm]
A_j={}&-\dfrac{(\Lambda_1^{+}+j)(\Lambda_2^{-}+q-j)U_j+(\Lambda_1^{-}+j)(\Lambda_2^{+}+q-j)V_j}{2(q-j)(\lambda_2+q-j)},
\end{array}
\]
for $1\leq j\leq q-1$ they are
\[
D_j=\frac{j(U_j+V_j)}{2(\lambda_2+q-j)},\qquad
C_j=\frac{j\big((\Lambda_2^{-}+q-j)U_j-(\Lambda_2^{+}+q-j)V_j\big)}{2(q-j)(\lambda_2+q-j)},
\]
and at the right end
\[
C_q=\frac{\Lambda_2^{+}V_{q-1}-\Lambda_2^{-}U_{q-1}}{2(\lambda_1+q)},\qquad
A_q=\frac{(\Lambda_1^{+}+q)\Lambda_2^{-}U_{q-1}+(\Lambda_1^{-}+q)\Lambda_2^{+}V_{q-1}}{2q(\lambda_1+q)}.
\]
\end{Statement}

For $q=0$, formula~\eqref{Bjodd} is $\cJ_{1/2,\xi}$ of formula~\eqref{J12f}.

Verified for $n\leq6$ by direct computation, with all four weights symbolic; for $n\leq4$ the
computation was repeated with all eight generators of $i(\fosp(2|2))$ and with a~much larger
set of test monomials.

\sssec{The operators of order $1$, $2$ and $3$} We write them in the basis
$\{1,\Dx,\De,P\}$, with $P$ as in formula~\eqref{PSL}, for each of the two arguments. The operators of order~$1$ are given by formula~\eqref{J12f}.
Order~2:
\be\label{J1xi}
\begin{array}{rl}
\cJ_{1,\xi}=&2(\lambda_1+1)(\lambda_2+1)\,\Dx\tp\De
+(\lambda_2+1)\mu_2^{-}(\mu_1^{+}+1)\,\partial_t\tp 1\\[2mm]
&-(\lambda_1+1)\mu_1^{+}(\mu_2^{-}+1)\,1\tp\partial_t
-(\lambda_2+1)\mu_2^{-}\,P\tp 1-(\lambda_1+1)\mu_1^{+}\,1\tp P,
\end{array}
\ee
and $\cJ_{1,\eta}$ is its mirror image ($\Dx\tp\De\mapsto\De\tp\Dx$,
$\mu^{\pm}\mapsto\mu^{\mp}$, $P\mapsto-P$); both map
$\cF_{\lambda_1,\nu_1}\tp\cF_{\lambda_2,\nu_2}\tto\cF_{\lambda_1+\lambda_2+2,\;\nu_1+\nu_2}$.
Order~3, charge $+1$:
\be\label{J32}
\begin{array}{rl}
\cJ_{3/2,\xi}=&(\lambda_2+1)\mu_2^{+}\,\partial_t\Dx\tp1
+(\lambda_1+1)\mu_1^{+}\,1\tp\partial_t\Dx
-(\lambda_2+1)(\mu_1^{+}+1)\,\partial_t\tp \Dx\\[2mm]
&-(\lambda_1+1)(\mu_2^{+}+1)\,\Dx\tp\partial_t
+(\lambda_2+1)\,P\tp \Dx+(\lambda_1+1)\,\Dx\tp P,
\end{array}
\ee
with target $\cF_{\lambda_1+\lambda_2+3,\;\nu_1+\nu_2-1}$, and $\cJ_{3/2,\eta}$ is its
mirror image. The operators of orders $4$ and $5$ have the same shape, with powers of
$\partial_t$ inserted and with coefficients that are products of the same linear forms
$\lambda_i+c$ and $\mu_i^{\pm}+c$.

\ssbegin[{What happens when $\nu_1\nu_2\omega\neq0$}]{Corollary}[What happens when $\nu_1\nu_2\omega\neq0$]\label{nu}
There \textup{is} something when $\nu_1\nu_2\omega\neq0$: the two operators of every order
exist for all weights, with $\omega=\nu_1+\nu_2-w$, where $w$ is the charge. What is special
about $\nu_1=\nu_2=0$ is the \textit{shape} of the answer, not its existence. The ansatz used in
part~\textup{(i)} of Theorem~$\ref{Omthm}$ --- linear combinations of
$\varphi^{(a)}\psi^{(b)}$, $\sum_s D_{\theta_s}(\varphi^{(a)})D_{\theta_s}(\psi^{(b)})$ and
$D_{\theta_1}D_{\theta_2}(\varphi^{(a)})D_{\theta_1}D_{\theta_2}(\psi^{(b)})$ --- contains
a~non-zero invariant operator if and only if $\nu_1=\nu_2=0$. In that case, it contains exactly the
line spanned by $\cJ_{1,\xi}+\cJ_{1,\eta}$, which is the operator~\eqref{J1}.
\end{Corollary}

\begin{proof}[Proof \textup{(for order 2; the general case is the same computation)}]
Due to the identities~\eqref{sym}, the ansatz in the $\xi\eta$ realization is the span of
$\partial_t\tp1$, $1\tp\partial_t$ and $\Dx\tp\De+\De\tp\Dx$: it contains neither
$P\tp1$ nor $1\tp P$, and the two $D$-terms enter with equal coefficients. Imposing these
three conditions on $a\cJ_{1,\xi}+b\cJ_{1,\eta}$ gives, by formulas~\eqref{J1xi},
$a=b$ (the coefficients of $\Dx\tp\De$ and $\De\tp\Dx$), $a\mu_2^-=b\mu_2^+$ (no $P\tp1$) and
$a\mu_1^+=b\mu_1^-$ (no $1\tp P$), whence $\nu_1=\nu_2=0$, for $\lambda_1,\lambda_2\neq-1$. \end{proof}

\sssec{What exactly is missing from the ansatz of \cite{Om}} To catch everything one must allow, besides
the three families printed in part~(i) of Theorem~\ref{Omthm}, also
\[
(D_{\theta_1}D_{\theta_2}\varphi^{(a)})\,\psi^{(b)},\qquad
\varphi^{(a)}\,(D_{\theta_1}D_{\theta_2}\psi^{(b)}),
\]
\[
D_{\theta_1}(\varphi^{(a)})D_{\theta_2}(\psi^{(b)})
-D_{\theta_2}(\varphi^{(a)})D_{\theta_1}(\psi^{(b)}),
\]
i.e., in the $\xi\eta$ realization, the operators $P\tp1$, $1\tp P$ and
$\Dx\tp\De-\De\tp\Dx$. Since $\fk(1|2)_0$ is $2$-dimensional, the charge is a~second grading
of the space of operators, and any ansatz must be closed under it; the printed one is not.
Consequently, already for $\nu_1=\nu_2=0$, the even invariant operators of each order form
a~pencil, not a~line. For order~$2$, the missing operator is $\cJ_{1,\xi}-\cJ_{1,\eta}$, the
printed $\cJ_1$ being the sum $\cJ_{1,\xi}+\cJ_{1,\eta}$. Moreover, there are odd invariant
operators of every odd order $\Nord=2k$, $k\in\nfrac12+\Zee_+$, whereas only $\Nord=1$ (i.e., $k=\nfrac12$) is listed.

\ssec{On the title of this section} In \cite[(4)]{BLS1}, the pair
$\fosp(n|2m+2)\subset\fk(2m+1|n)$ is considered; for $m=0$ this is indeed $\fosp(n|2)\subset\fk(1|n)$, but
only $n=1,2$ are treated here. The case $n=3$ is the subject of \cite{BLO}; we verify it and
correct it in \S\ref{Sblo}. For $n=4$, the pair is 
a~one-parameter family $\fosp_a(4|2)\subset\fk(1|4)$, see \S\ref{Sn4}; the case $n=5$ is
\S\ref{Sn5}, and for $n=6$ we consider not $i(\fosp(6|2))$ but $\fas\subset\fk\fas$, see
\S\ref{GRCAS}.

\section{Super GRC brackets: the $\fosp(3|2)\subset\fk(1|3)$ case;
\cite{BLO} verified and corrected}\label{Sblo}

\ssec{Summary of this Section}\label{blosummary}
The case $n=3$ of the pair $\fosp(n|2)\subset\fk(1|n)$ is treated in \cite{BLO}. We have
checked it against an independent computation carried out from scratch \textit{in the
conventions of} \cite{BLO}, with symbolic weights. Formulas
(3.2)--(3.3) of \cite{BLO} for $N=1$ and formulas (3.11)--(3.12) of \cite{BLO} for $N=2$ are invariant as printed.
The $\fosp(3|2)$ part of \cite{BLO}, however, contains three errors, all of them in its
Theorem 4.1:

\textbf{(a)} In the even supertransvectants, formula (4.2) of \cite{BLO}, three of the four groups of
terms carry the wrong sign;

\textbf{(b)} In the odd supertransvectants, formula (4.1) of \cite{BLO}, the two sums
over $1\leq p<\ell\leq 3$ have lost the Levi-Civita sign;

\textbf{(c)} The description of the exceptional weights
in items (ii)--(v) of Theorem 4.1 of \cite{BLO} is wrong. The true set of weights at which the space
of invariant operators is $2$-dimensional is different from theirs. For the indices $k=1$
and $k=2$, i.e., in the orders $\Nord=2$ and $\Nord=4$, there are no exceptional weights at all. In cases (a) and (b), the corrections are stated
explicitly and verified, see Statements~\ref{st:42} and \ref{st:41}; in case (c) the true
answer for $k\leq2$ is given in Statement~\ref{st:exc}. What Theorem 4.1 of \cite{BLO} gets
right is listed in \S\ref{s:ok}.

\ssec{The conventions of \cite{BLO}, and the ansatz}\label{s:setup}

\sssec{Conventions, taken from Section II of \cite{BLO}}\label{s:conv3}
Everything in this Section is computed in the conventions of \cite{BLO}, not in ours: on
$\Ree[x,\theta_1,\dots,\theta_N]$,
\begin{equation}\label{eq:conv}
\etb_i=\frac{\partial}{\partial \theta_i}-\theta_i\frac{\partial}{\partial x},
\qquad
X_F=F\partial_x-\frac12(-1)^{p(F)}\sum_{1\leq i\leq N}\etb_i(F)\,\etb_i,
\qquad
\fL^{\lambda}_{X_F}=X_F+\lambda F',
\end{equation}
with $\{F,G\}$ their formula (2.2), with
$\fosp(N|2)=\Span(X_1,X_x,X_{x^2},X_{\theta_i},X_{x\theta_i},X_{\theta_i\theta_j})$, and with
\begin{equation}\label{eq:Gamma}
\Gamma^{\tau,\lambda}_{i,j,k}=(-1)^j\binom{2\tau+[k]}{j}\binom{2\lambda+[k]}{i}
\end{equation}
their formula (3.4). Their $\etb_i$ is thus \textbf{not} our $D_{\theta_i}$ of eq.~\eqref{KD}, and
their $\cF^N_\tau$ is a~$1$-parameter family of weighted densities; $\cK(N)$ denotes the $\fk(1|N)$-module of contact vector fields. A~bilinear operator
$T\colon\cF^N_\tau\tp\cF^N_\lambda\longrightarrow\cF^N_{\tau+\lambda+k}$ is called invariant if
\[
\fL^{\mu}_{X_F}\circ T=T\circ(\fL^{\tau}_{X_F}\tp1+1\tp\fL^{\lambda}_{X_F})
\]
with the sign rule taken into account, i.e.,
\[
\fL^{\mu}_{X_F}(T(G,H))=(-1)^{p(F)p(T)}\big(T(\fL^{\tau}_{X_F}G,H)+(-1)^{p(F)p(G)}T(G,\fL^{\lambda}_{X_F}H)\big)
\]
for every generating function $F$ of $\fosp(N|2)$; here
$\mu=\tau+\lambda+k$.

\sssec{A~dictionary between the notation of \cite{BLO} and ours}\label{s:dict}
Put $x=-t$. Then, their contact form is $dx+\sum_{1\leq i\leq N}\theta_i\,d\theta_i=-\big(dt-\sum_{1\leq i\leq N}\theta_i\,d\theta_i\big)$, and their operators $\etb_i$ become our operators $D_{\theta_i}$ of eq.~\eqref{KD}. Their weighted density of weight~$\tau$ is $f\alpha^\tau$, whereas ours of weight~$\lambda$ is $f\alpha^{\lambda/2}$. So
\begin{equation}\label{bloDict}
\begin{array}{lcl}
x,\ \theta_i,\ \etb_i\ \text{of \cite{BLO}}&\longleftrightarrow&-t,\ \theta_i,\ D_{\theta_i},\\[1mm]
\cF^N_\tau\ \text{of \cite{BLO}}&\longleftrightarrow&\cF_{2\tau},\qquad \cK(N)\simeq\cF^N_{-1}\longleftrightarrow\cF_{-2},\\[1mm]
\text{supertransvectant of index }k&\longleftrightarrow&\text{bracket of order }\Nord=2k.
\end{array}
\end{equation}
For example, the exceptional pairs $(\tau,\lambda)$ of Statement~\ref{st:exc} are the pairs $(\lambda_1,\lambda_2)=(2\tau,2\lambda)$ in our notation, and the involution of Remark~\ref{rem:refl} is $\lambda_i\longmapsto-2-\lambda_i$.

The classification itself is done by an ansatz
$T=\sum c_{a,A,b,B}\,\partial_x^{a}\etb_A\tp\partial_x^{b}\etb_B$ with constant coefficients
(legitimate: invariance under $X_1$ and $X_{\theta_i}$ forces it, cf. Theorem~\ref{Peetre}),
followed by solving the resulting linear system over $\Qee$ or over $\Qee(\tau,\lambda)$.

\ssec{The even supertransvectants (4.2) of \cite{BLO}: a~sign misprint in three of the four
groups}\label{s:42}

\sssec{Formula $(4.2)$ of \cite{BLO}, and the signs it must have}\label{s:42a}
Written out, formula (4.2) of \cite{BLO} is the sum over the subsets $A\subseteq\{1,2,3\}$
\begin{equation}\label{eq:42}
\Jop^{\tau,\lambda,3}_k(F,G)=\sum_{0\leq r\leq 3}\varepsilon_r\sum_{|A|=r}\
\sum_{i+j=k-r}(-1)^{r\,p(F)}\,\Gamma^{\tau,\lambda}_{i,j,k-1}\,
\etb_A(F^{(i)})\,\etb_A(G^{(j)}),
\end{equation}
where $\etb_A=\etb_{p_1}\cdots\etb_{p_r}$ for $A=\{p_1<\dots<p_r\}$, and where the signs
\textit{as printed} are
\[
(\varepsilon_0,\varepsilon_1,\varepsilon_2,\varepsilon_3)=(-,-,-,+).
\]
This operator is \textbf{not invariant} for generic $(\tau,\lambda)$ and $k=1,2,3$. 

\sssbegin[{The correct signs in eq.~$(4.2)$ of \cite{BLO}}]{Statement}[The correct signs in eq.~$(4.2)$ of \cite{BLO}]\label{st:42}
Formula \eqref{eq:42} becomes invariant --- and then, being nonzero, is the
supertransvectant --- precisely for
\[
(\varepsilon_0,\varepsilon_1,\varepsilon_2,\varepsilon_3)=(-,+,+,-),
\qquad\text{i.e.}\qquad
\varepsilon_r=-(-1)^{r(r+1)/2},
\]
up to an overall factor. All $16$ sign patterns were tried: for $k=3$, where all four sums over $|A|=r$, $r=0,1,2,3$, in formula~\eqref{eq:42} are nonempty, only this pattern and its negative survive; the corrected operator was then verified
invariant and nonzero for $k=1,2,3,4$.
\end{Statement}

Two remarks make the misprint evident. First, only $\varepsilon_0$ is right as printed; the
three sums that carry $\etb$'s are all wrong. Second, the correct pattern is \textit{the pattern
\cite{BLO} themselves use} in the cases $N=1$ and $N=2$: their formula (3.3) has
$(\varepsilon_0,\varepsilon_1)=(-,+)$ and their formula (3.11) has
$(\varepsilon_0,\varepsilon_1,\varepsilon_2)=(-,+,+)$, both in agreement with
$\varepsilon_r=-(-1)^{r(r+1)/2}$. So formula (4.2) of \cite{BLO} is the only place where the pattern is broken.

The smallest case makes it concrete. For the index $k=1$, i.e., in order $\Nord=2$, the unique
invariant operator is,
with symbolic weights,
\begin{equation}\label{eq:k1}
\Jop^{\tau,\lambda,3}_1(F,G)=2\tau FG'-2\lambda F'G
+(-1)^{p(F)}\sum_{1\leq p\leq 3}\etb_p(F)\,\etb_p(G),
\end{equation}
whereas formula (4.2) of \cite{BLO} as printed gives the same first two terms and the term
\[
-(-1)^{p(F)}\sum_{1\leq p\leq 3}\etb_p(F)\etb_p(G).
\]
Note that the $r=0$ part of formula (4.2) of \cite{BLO} reproduces the terms
$2\tau FG'-2\lambda F'G$ exactly, so the $\Gamma$'s are right.

\ssec{The odd supertransvectants (4.1) of \cite{BLO}: a~missing Levi-Civita sign}\label{s:41}

\sssec{The correction}\label{s:41a}
Here, the situation is better and worse at the same time. On the one hand, \textit{all the coefficients of formula $(4.1)$ of \cite{BLO}
have the right absolute value}: the $\Gamma$'s, the $\Lambda^i_j=(j+2)(j+1)(i+1)$ and the
index shifts $i+1,j+2$ etc. are all correct. On the other hand, the two sums over $1\leq p<\ell\leq 3$
give one and the same coefficient to all three pairs $(p,\ell)$, whereas the invariant
operator requires the three to alternate.

\sssbegin[{The missing Levi-Civita sign in eq.~$(4.1)$ of \cite{BLO}}]{Statement}[The missing Levi-Civita sign in eq.~$(4.1)$ of \cite{BLO}]\label{st:41}
In the first two sums of formula $(4.1)$ of \cite{BLO}, the terms must be multiplied by the Levi-Civita
sign, see~\cite{LC}, i.e., by the sign of the permutation $(1,2,3)\mapsto(p,\ell,q)$ of the three
distinct indices:
\[
\epsilon_{p\ell q},\qquad q=6-\ell-p,
\qquad\text{i.e.}\qquad
+1\ \text{ for }(p,\ell)=(1,2),\,(2,3),\qquad -1\ \text{ for }(p,\ell)=(1,3),
\]
after which formula $(4.1)$ --- with an overall change of sign, which is only a~normalization --- is
invariant.
\end{Statement}

The comparison was made term by term: for $k=\frac32$ and $k=\frac52$ the unique invariant
operator was computed and matched against formula (4.1) of \cite{BLO} in the same basis. At $k=\frac32$, four of the
eight coefficients differ from those of \cite{BLO}, at $k=\frac52$ eight of the sixteen; in
every case the moduli agree and only the sign is off, and exactly on the pairs
$(p,\ell)=(1,2)$ and $(2,3)$.

Explicitly, for the index $k=\frac32$, i.e., in order $\Nord=3$, the invariant operator is,
after clearing
denominators,
\begin{equation}\label{eq:k32}
\begin{array}{r@{}l}
\Jop^{\tau,\lambda,3}_{3/2}(F,G)={}&-(-1)^{p(F)}\tau(\tau+1)(2\tau+1)F\,\etb_1\etb_2\etb_3(G)\\
&{}+(\tau+1)(2\tau+1)(\lambda+1)S_1(F,G)\\
&{}+(-1)^{p(F)}(\tau+1)(\lambda+1)(2\lambda+1)S_2(F,G)\\
&{}+\lambda(\lambda+1)(2\lambda+1)\,\etb_1\etb_2\etb_3(F)\,G,
\end{array}
\end{equation}
where
\[
\begin{array}{l}
S_1(F,G):=\etb_1(F)\etb_2\etb_3(G)-\etb_2(F)\etb_1\etb_3(G)+\etb_3(F)\etb_1\etb_2(G),\\
S_2(F,G):=-\etb_1\etb_2(F)\etb_3(G)+\etb_1\etb_3(F)\etb_2(G)-\etb_2\etb_3(F)\etb_1(G).
\end{array}
\]
The signs $+,-,+$ in $S_1$ and $-,+,-$ in $S_2$ are $\epsilon_{p\ell q}$ and $-\epsilon_{p\ell q}$, respectively, with $\epsilon_{p\ell q}$ the Levi-Civita sign
of Statement~\ref{st:41}; it is absent from formula (4.1) of \cite{BLO}.

\ssec{The exceptional weights: Theorem 4.1(ii)--(v) of \cite{BLO} is wrong}\label{s:exc}

\sssec{The exceptional weights for $\fosp(3|2)$, recomputed}\label{s:exca}
Let $d(\tau,\lambda,k)$ denote the dimension of the space of $\fosp(3|2)$-invariant bilinear
operators $\cF^3_\tau\tp\cF^3_\lambda\longrightarrow\cF^3_{\tau+\lambda+k}$, the superscript being the number $N=3$ of odd coordinates, as in $\cF^N_\tau$ of \cite{BLO}; in our notation, it is $\cF_{2\tau}$, see the dictionary~\eqref{bloDict}.  It was computed on the
grid $\tau,\lambda\in\{-3,-\frac52,-2,-\frac32,-1,-\frac12,0,\frac12,1\}$ together with two
generic values, for $k=\frac12,1,\frac32,2$.

\sssbegin[{The exceptional weights for $\fosp(3|2)\subset\fk(1|3)$}]{Statement}[The exceptional weights for $\fosp(3|2)\subset\fk(1|3)$]\label{st:exc}
The dimension is $d=1$ for all these weights and all $k\in\{1,\frac32,2\}$, except that for $k=\frac32$
\[
d(\tau,\lambda,\tfrac32)=2\quad\text{exactly for}\quad
(\tau,\lambda)\in\{(-1,-1),\ (-1,-\tfrac12),\ (-1,0),\ (-\tfrac12,-\tfrac12)\}
\]
and at the pairs $(\lambda,\tau)$ obtained from these by interchanging the arguments;
and $d=0$ for $k=\frac12$. In particular, \textbf{for $k=1$ and $k=2$ there are no exceptional
weights at all} on this grid.
\end{Statement}

Interchanging the arguments, $B\mapsto B'$ with $B'(G,F):=(-1)^{p(F)p(G)}B(F,G)$, is a bijection
between the invariant operators $\cF^3_\tau\tp\cF^3_\lambda\tto\cF^3_{\tau+\lambda+k}$ and
$\cF^3_\lambda\tp\cF^3_\tau\tto\cF^3_{\tau+\lambda+k}$, so $d(\tau,\lambda,k)=d(\lambda,\tau,k)$.

Theorem 4.1(v) of \cite{BLO} says instead that $d=2$ whenever $\tau=\frac{1-m}4$,
$\lambda=\frac{1-\ell}4$ with $m,\ell$ odd and $\frac{\ell+m}2\leq k+2$. For $k=\frac32$, this
predicts $d=2$ exactly at
\[
(\tau,\lambda)\in\{(0,0),\ (0,-\tfrac12),\ (0,-1),\ (-\tfrac12,-\tfrac12)\}
\ \text{ and the interchanged pairs},
\]
and for $k=1$ and $k=2$ it predicts a~nonempty list as well. So their Theorem 4.1(iv)--(v) is
wrong in both directions: it claims $d=2$ at $(0,0)$ and $(0,-\frac12)$, where in fact $d=1$.
It also claims uniqueness at $(-1,-1)$ and $(-1,-\frac12)$ --- these fall under their (iv),
since $\frac{\ell+m}{2}>k+2$ there --- where in fact $d=2$.

Both halves of this are certain, for different reasons:

1) A~computed value $d=1$ is not an underestimate.  The classifier solves a~linear system built from
test monomials of $x$-degree $\leq t_{\max}$ \textup{(}the degree of the \textit{monomial}, not
of the operator\textup{)}. Taking $t_{\max}$ too small can only
\textit{inflate} the answer. Moreover, $d\geq1$ for all weights: the coefficients of the generic operator are rational functions of $(\tau,\lambda)$; clearing the denominators and dividing by the highest common power of $\tau-\tau_0$, $\lambda-\lambda_0$, one gets a~nonzero invariant operator at every point $(\tau_0,\lambda_0)$ \textup{(}equivalently: the rank of the linear system can only drop at special weights, so the dimension of its solution space can only grow\textup{)}.  Getting $1$ therefore proves $d=1$. This kills the claim of
$2$-dimensionality at $(0,0)$ and $(0,-\frac12)$.

2) The values $d=2$ were confirmed independently:  for
$t_{\max}=3,4,5$, the dimension is stable, and both basis operators produced at each of the four points were fed back
into the invariance test and are invariant and nonzero.

\sssbegin[{Where the slip might have entered}]{Remark}[Where the slip might have entered]\label{rem:refl}
The set predicted by Theorem 4.1(v) of \cite{BLO} and the true set are exchanged by the involution
\[
\tau\longmapsto-1-\tau,\qquad \lambda\longmapsto-1-\lambda,
\]
which fixes $(-\frac12,-\frac12)$ and swaps $(0,-1)\leftrightarrow(-1,0)$ --- the two points
where the two answers happen to agree. In the normalization \eqref{eq:conv}, one has
$\cK(3)\simeq\cF^3_{-1}$, so this involution sends $\cF^3_\tau$ to the weight complementary to~$\tau$
with respect to $\cK(3)$. It is \textbf{not} the duality of the densities: the volume forms are
$\Vol\simeq\cF^N_{(2-N)/2}$ \textup{(}our $\cF_{2-n}$, cf. \S\ref{GRCAS}\textup{)}, i.e.,
$\cF^3_{-1/2}$ for $N=3$, so the pairing $\cF^3_\tau\tp\cF^3_{\mu}\tto\Vol$ exists for
$\mu=-\frac12-\tau$. It looks as though the parametrization $\tau=\frac{1-m}4$, natural for $N=1$, had been used for $N=3$ without applying the above involution $\tau\mapsto-1-\tau$, $\lambda\mapsto-1-\lambda$.  Note also that the operator $\etb_1^{m-4}\etb_2\etb_3$ of their formula (4.3) requires $m\geq5$, i.e.\ $\tau\leq-1$, whereas the set $I_{k+2}:=\{0,-\frac12,-1,\dots,-\frac{[k]+2}2\}$ of their Theorem 4.1\textup{(}ii\textup{)} starts at $\tau=0$. The classification of the exceptional weights for all $k$ is an open problem.
\end{Remark}

\ssec{What Theorem 4.1 of \cite{BLO} gets right}\label{s:ok}

For the record:

1) \textbf{There is no operator of order $\Nord=1$}, i.e., for the index $k=\frac12$: the
computation gives
$d=0$ for all weights tested. The restriction $2k\in\Nee\setminus\{0,1\}$ in their
Theorem 4.1(i) is thus correct; there is no analogue of $\cJ_{1/2}$ here.

2) \textbf{Generic uniqueness holds}: $d=1$ for every $k\geq1$ and generic weights, in
agreement with their (i)--(ii). (Unlike the case $N=2$, where the generic dimension is $2$,
see \S\ref{classification}.)

3) The coefficients $\Gamma^{\tau,\lambda}_{i,j,k}$ and $\Lambda^i_j$, and all the index
shifts, are correct: the errors in formulas (4.1) and (4.2) of \cite{BLO} are purely errors of sign.

4) Their $N=1$ and $N=2$ formulas (3.2)--(3.3) and (3.11)--(3.12) are correct as printed.

\section{Super GRC and the Bol operators for the family $\fosp_a(4|2)\subset\fk(1|4)$
}\label{Sn4}

\ssec{Summary of this Section}\label{n4summary}
For $n=4$, there is not one subalgebra to consider, but a~one-parameter family of them. The
$\Zee$-graded subalgebras $\fh\subset\fk(1|4)$ with $\fh_-=\fk(1|4)_-$ and $\dim\fh=17$ form
a~family $\fh(a)$, where $a\in\Cee\cup\{\infty\}$, and $\fh(a)\simeq\fosp_\alpha(4|2)$. This is the
family described by eq.~(65) of \cite{BLS1}. We answer four questions about it.

1) \textbf{The GRC brackets.} For generic weights, the space of $\fh(a)$-invariant binary
operators $\cF_{\lambda_1}\tp\cF_{\lambda_2}\tto\cF_{\lambda_1+\lambda_2+\Nord}$ is of dimension
$0,1,0,2$ for $\Nord=1,2,3,4$: \textit{in odd orders there is nothing at all}, in contrast with
$n=1,2,3$, see Theorem~\ref{n4thm}.

2) \textbf{Where the parameter is felt.} The order-$2$ bracket is the same for every
member of the family. It is the contact bracket, invariant under the whole of $\fk(1|4)$, see
Statement~\ref{n4contact}, so it says nothing about the subalgebra. The order-$4$ space, on the
contrary, \textit{does} depend on $a$: the planes obtained for $a=0,2,3$ span a~$5$-dimensional
space, not a~$2$-dimensional space.

3) \textbf{The Bol operators.} They see the parameter too, and in a~transparent way:
for every $a$ and every $k\geq0$ there is exactly one Bol operator, of order $\Nord=4+2k$, namely
\[
T_4\circ\partial_t^k+\frac{a}{k+2}\partial_t^{k+2}:\cF_{-3-k}\tto\cF_{1+k},
\]
where $T_4$ is the operator~\eqref{T4} and $\circ$ is the composition of operators\textup{:}
$(T_4\circ\partial_t^k)(\varphi)=T_4(\partial_t^k(\varphi))$. Since $\partial_t$ commutes with
every $D_{\xi_i}$ and $D_{\eta_i}$, the order of the two factors does not matter. There is
nothing else at all, see Statement~\ref{n4bol}. This confirms \cite[Thm~1.6]{BLS1} --- obtained
there by a~different method --- by an independent computation in the conventions of the present
paper. The dictionary between the two parametrizations is $a'=\frac{1-a}{1+a}$, where $a'$ is the parameter of eq.~(65) of \cite{BLS1}, see
\S\ref{n4dict}. This dictionary identifies the two answers, including both degenerate members.

4) \textbf{Isomorphic algebras, non-equivalent embeddings.} In \cite[\S3.2.3a]{BLS1}, the following is asked: since $\fosp_1(4|2)\simeq\fosp_{-2}(4|2)\simeq\fosp_{-1/2}(4|2)$, why are the
answers for these values different --- ``Is there a mistake?\rq\rq. There is no mistake, and the
reason can be made precise: the group acting on the parameter of the \textit{embedding} is
the group $\Zee/2$ generated by the involution $a\mapsto-a$, which is induced by an automorphism of $\fk(1|4)$, whereas the
$S_3$ of the abstract classification is not realized by automorphisms of $\fk(1|4)$. The
invariant of the embedding is~$a^2$, i.e., $\alpha+\alpha^{-1}$, not the invariant $j(\alpha)=\frac{(\alpha^2+\alpha+1)^3}{\alpha^2(\alpha+1)^2}$ of the abstract algebra $\fosp_\alpha(4|2)$, see Statement~\ref{n4fam}.  In particular, $\fosp(4|2)$ itself sits in $\fk(1|4)$ in exactly \textit{two} non-equivalent ways, $a=0$ and
$a=\pm3$, with different Bol operators, see \S\ref{n4answer}. Thus, the Bol operators and the GRC brackets depend not only on the isomorphism class of the subalgebra but on its embedding. 

The relation between 1) and 3) is the relation already met in \S\ref{CoRm}. At the weights at which
a~Bol operator exists, it enters the binary picture as ``a~Bol operator acting in one of the two
arguments\rq\rq, i.e., as the binary operators $(\varphi,\psi)\mapsto B(\varphi)\psi$ and $(\varphi,\psi)\mapsto\varphi B(\psi)$, see \S\ref{n4both}.  

\ssec{Where the parameter comes from}\label{n4family}
Grade $\fk(1|4)$ by $\deg t=2$, $\deg\xi_i=\deg\eta_i=1$ for all $i$, so that $\deg K_f=\deg f-2$. Then,
\[
\begin{array}{l}
\fg_{-2}=\langle K_1\rangle,\quad \fg_{-1}=\langle K_{y_i}\mid \text{all $i$}\rangle,\\
\fg_0=\langle K_t,K_{y_iy_j}\mid \text{all $i, j$}\rangle=\Cee\oplus\fo(4),\quad
\fg_1=\langle K_{ty_i}\rangle\oplus\langle K_{y_iy_jy_k}\mid \text{all $i, j, k$}\rangle,
\end{array}
\]
where $y_1,\dots,y_4$ stands for $\xi_1,\eta_1,\xi_2,\eta_2$. As $\fo(4)$-modules,
$\fg_1=V\oplus\Lambda^3V$ with $V\simeq\Lambda^3V$ irreducible. Therefore, the
$\fg_0$-stable complements to $\Lambda^3V$ in $\fg_1$ form a~one-parameter family. They are the
graphs of $a\sigma$, where $\sigma:V\tto\Lambda^3V$ is the equivariant map, unique up to
a~factor by Schur's lemma. This is where the parameter of $\fosp_a(4|2)$ comes from.

\sssbegin[{The family $\fh(a)\simeq\fosp_a(4|2)\subset\fk(1|4)$}]{Statement}[The family $\fh(a)\simeq\fosp_a(4|2)\subset\fk(1|4)$]\label{n4fam}
Put, for every $a\in\Cee$,
\[
\begin{array}{ll}
w_{\xi_1}:=t\xi_1-a\,\xi_1\xi_2\eta_2,\quad&
w_{\eta_1}:=t\eta_1+a\,\eta_1\xi_2\eta_2,\\[1mm]
w_{\xi_2}:=t\xi_2-a\,\xi_1\eta_1\xi_2,\quad&
w_{\eta_2}:=t\eta_2+a\,\xi_1\eta_1\eta_2,
\end{array}
\]
and let $\fh(a)$ be the subalgebra generated by $\fk(1|4)_{\leq0}$ and by~$K_{w_{y_i}}$ for all~$i$.
Then:

\textup{(i)} $\{w_{\xi_1},w_{\eta_1}\}=\{w_{\xi_2},w_{\eta_2}\}=-t^2+2a\,\xi_1\eta_1\xi_2\eta_2$,
all the other brackets of the $w$'s vanish, and $\{w_{y_i},\{w_{y_j},w_{y_k}\}\}=0$. Hence,

\[
\fh(a)=\fg_{-2}\oplus\fg_{-1}\oplus\fg_0\oplus\langle K_{w_{y_i}}\mid \text{ for all $i$}\rangle\oplus
\langle K_{t^2-2a\,\xi_1\eta_1\xi_2\eta_2}\rangle
\]
is a~graded subalgebra of dimension $17$ for every $a$, and $\fh(0)=i(\fosp(4|2))$. Moreover,
\[
[\fk(1|4)_{-1},\ \langle K_{w_{y_i}}\rangle]=\fg_0\quad\text{if and only if}\quad a\neq\pm1 .
\]

\textup{(ii)} $\fh(a)\simeq\fosp_\alpha(4|2)$ with
\[
\alpha=\frac{1+a}{1-a},\qquad
j(\alpha):=\frac{(\alpha^2+\alpha+1)^3}{\alpha^2(\alpha+1)^2}=\frac{(a^2+3)^3}{4(a^2-1)^2},
\]
$j$ being the invariant of the $S_3$-action \eqref{S3act} on the parameter. In particular,
$\fh(a)\simeq\fosp(4|2)$ \textup{(}i.e., $j=\nfrac{27}4$\textup{)} if and only if $a\in\{0,3,-3\}$,
i.e., if and only if $\alpha\in\{1,-2,-\nfrac12\}$.

\textup{(iii)} For $a=\pm1$, i.e., for $\alpha=\infty$ and $\alpha=0$, the algebra $\fh(a)$ is not
simple; the subalgebra generated by $\fg_-$ and $\langle K_{w_{y_i}}\rangle$ is then
$14$-dimensional and isomorphic to $\fpsl(2|2)$, the ``small $N=4$\rq\rq\  subalgebra of $\fk(1|4)$.
The remaining degenerate member of the family is $a=\infty$, i.e.,
$\fh(\infty)=\langle\fk(1|4)_{\leq0},K_{y_iy_jy_k},K_{\xi_1\eta_1\xi_2\eta_2}\rangle$,
for which $\alpha=-1$.
\end{Statement}

\begin{proof}[Computation]  Here, $\fosp_\alpha(4|2)$ is defined by its Cartan matrix, as in \cite[\S3.2.3]{BLS1}. $\sigma$ is found by intertwining the $\fo(4)$-action on $\langle ty_i\rangle$ and on the
cubics; the solution space is $1$-dimensional, as Schur's lemma predicts. Part (i) is a~direct
check of the contact brackets. For (ii), fix the maximal toral subalgebra
$\langle K_t,K_{\xi_1\eta_1},K_{\xi_2\eta_2}\rangle$, normalize the three $\fsl(2)$-triples by
$\beta_i(H_i)=2$, and take an odd root vector $X$ of weight $\delta=\eps_1+\eps_2+\eps_3$
together with the opposite vector $Y$. Then, $[X,Y]=\sum c_iH_i$ with
$(c_1,c_2,c_3)=(1,\frac{1-a}2,\frac{1+a}2)$, and $(\eps_i,\eps_i)\propto c_i$ for the invariant
form. Comparing this with the standard normalization
$-(1+\alpha):1:\alpha$, we get the value of $\alpha$. The identity $c_1=c_2+c_3$ is the isotropy of the
odd roots; at $a=0$ the recipe returns $\alpha=1$, i.e., $\fosp(4|2)$, as it must.
\end{proof}

\sssec{The two parametrizations}\label{n4dict}
Recall from \S3.2.3 of \cite{BLS1} that the abstract algebra $\fosp_\alpha(4|2)$ is
determined by the coefficients $\sigma_1,\sigma_2,\sigma_3$ with $\sum\sigma_i=0$, with which the three
projections $S^2(\fg_\od)\longrightarrow\fsl_i(2)$ enter the bracket, and that
\begin{equation}\label{S3act}
\alpha\longmapsto-1-\alpha,\qquad \alpha\longmapsto\alpha^{-1}
\end{equation}
generate the $S_3$ under which the isomorphism class does not change. The parameter $a$ of
Statement~\ref{n4fam} is a~coordinate on the family of \textit{embeddings}. The parameter of
eq.~(65) of \cite{BLS1}, which we denote here by $a'$ in order not to confuse it with ours, is
another coordinate, and
\begin{equation}\label{dicta}
a=\frac{1-a'}{1+a'},\qquad a'=\frac{1-a}{1+a},\qquad \alpha=\frac1{a'} .
\end{equation}
The three degenerate members are $a=\pm1$ and $a=\infty$, i.e., $\alpha\in\{\infty,0,-1\}$: the
values at which $\fosp_\alpha(4|2)$ is not simple.

\ssec{The GRC brackets of the family}\label{n4grc}
Everything below is computed in the realization \eqref{KD} of $\fk(1|4)$, with
$\alpha_1=dt+\sum_{i=1,2}(\xi_id\eta_i+\eta_id\xi_i)$ and with the $1$-parameter weighted
densities $\cF_\lambda=\cF\alpha_1^{\lambda/2}$, on which
$\fL^\lambda_{K_f}(\varphi)=K_f(\varphi)+\lambda\,\partial_t(f)\varphi$. By
Theorem~\ref{Peetre}, an invariant operator has constant coefficients in the $D_f$, so the
ansatz for the binary operators of order $\Nord$ is
$\sum c\,\partial_t^{p}D_A\tp\partial_t^{q}D_B$ with $2p+|A|+2q+|B|=\Nord$.

\sssbegin[{The GRC brackets for $\fosp_a(4|2)\subset\fk(1|4)$}]{Theorem}[The GRC brackets for $\fosp_a(4|2)\subset\fk(1|4)$]\label{n4thm}
Let $\Nord\leq4$. For generic weights, the space of $\fh(a)$-invariant binary differential
operators $\cF_{\lambda_1}\tp\cF_{\lambda_2}\tto\cF_{\lambda_1+\lambda_2+\Nord}$ is of dimension
\[
0\ (\Nord=1),\qquad 1\ (\Nord=2),\qquad 0\ (\Nord=3),\qquad 2\ (\Nord=4)
\]
for every $a\in\Cee$, the degenerate members $a=\pm1$ included. The order-$2$ operator is the same for
all these members of the family, namely
\begin{equation}\label{J2n4}
\cJ_2=\lambda_2\,\partial_t\tp1-\lambda_1\,1\tp\partial_t
+\sum_{i=1,2}\big(D_{\xi_i}\tp D_{\eta_i}+D_{\eta_i}\tp D_{\xi_i}\big);
\end{equation}
the order-$4$ space is \textit{not}: the three planes obtained for $a=0,2,3$ at
$(\lambda_1,\lambda_2)=(\nfrac37,-\nfrac5{11})$ span a~$5$-dimensional space.

For the remaining degenerate member $a=\infty$, the dimensions are $0,3,0,5$ for
$\Nord=1,2,3,4$: the algebra $\fh(\infty)$ of Statement~$\ref{n4fam}$\textup{(iii)} contains no
$K_{t^2}$, so each of the three summands $\partial_t\tp1$, $1\tp\partial_t$ and
$\sum_{i=1,2}\big(D_{\xi_i}\tp D_{\eta_i}+D_{\eta_i}\tp D_{\xi_i}\big)$ of formula~\eqref{J2n4} is
$\fh(\infty)$-invariant by itself.
\end{Theorem}
\sssbegin[The first GRC bracket is the contact bracket]{Statement}[The first GRC bracket
is the contact bracket, and it is invariant under all of $\fk(1|4)$]\label{n4contact}
The operator \eqref{J2n4} is invariant not only under $\fh(a)$ but under the whole of
$\fk(1|4)$, and at $\lambda_1=\lambda_2=-2$ it is the contact bracket. \textup{(}The same holds
for every $\fk(1|n)$, cf. Statement~$\ref{n5contact}$ for $n=5$.\textup{)} This is why the first GRC bracket says nothing about the
subalgebra.
\end{Statement}

\sssbegin[{The two subalgebras at $a=\pm1$}]{Remark}[The two subalgebras at $a=\pm1$]\label{n4psl}
At $a=\pm1$ one should distinguish $\fh(\pm1)$, which is $17$-dimensional and not simple, from
the subalgebra $\fpsl(2|2)$ generated by $\fg_-$ and $\fh_1(\pm1)$, which is $14$-dimensional.
Theorem~\ref{n4thm} concerns $\fh(\pm1)$. For the subalgebra $\fpsl(2|2)$, the dimensions are larger, as they
should be for a~smaller subalgebra: $0,4,0,5$ for $\Nord=1,2,3,4$, and there are three unary
invariant operators of order~$2$ at $\lambda=-2$, whereas $\fh(\pm1)$ has none.
\end{Remark}

\ssec{The Bol operators of the family}\label{n4bolsec}
Set
\begin{equation}\label{T4}
T_4:=\big(D_{\xi_1}D_{\eta_1}-\partial_t\big)\big(D_{\xi_2}D_{\eta_2}-\partial_t\big)
=D_{\xi_1}D_{\eta_1}D_{\xi_2}D_{\eta_2}
-\partial_t\big(D_{\xi_1}D_{\eta_1}+D_{\xi_2}D_{\eta_2}\big)+\partial_t^2,
\end{equation}
the product of the commutators $D_{\xi_i}D_{\eta_i}-\partial_t$ over both pairs; in the algebra generated by the $D_{\xi_i}$ and $D_{\eta_i}$, which is a~Clifford algebra over $\Cee[\partial_t]$, this is the normally ordered monomial of top degree; recall that $D_{\xi_i}D_{\eta_i}
+D_{\eta_i}D_{\xi_i}=2\partial_t$, so that $D_{\xi_i}D_{\eta_i}-\partial_t
=\nfrac12[D_{\xi_i},D_{\eta_i}]$ is the commutator, not the anticommutator.

\sssbegin[{Bol operators for $\fosp_a(4|2)\subset\fk(1|4)$}]{Statement}[Bol operators for $\fosp_a(4|2)\subset\fk(1|4)$]\label{n4bol}
For every $a$ and $0\leq k\leq2$, the operator
\begin{equation}\label{Bola}
\Bol_k^{(a)}:=T_4\circ\partial_t^{\,k}+\frac{a}{k+2}\,\partial_t^{\,k+2}:\quad
\cF_{-3-k}\tto\cF_{1+k}
\end{equation}
is $\fh(a)$-invariant. Up to a~factor, these are the only non-scalar $\fh(a)$-invariant
unary operators $\cF_\lambda\tto\cF_{\lambda+\Nord}$ with $\Nord\leq8$: there are none of odd order, none of
order~$2$, and for $\Nord=4+2k$ the space is $1$-dimensional at $\lambda=-3-k$ and zero at every
other weight. For the degenerate member $a=\infty$, i.e., $\fh(\infty)$ of
Statement~$\ref{n4fam}$\textup{(iii)}, the answer is different in kind: there the invariant
operators are $\partial_t^{\,k}:\cF_\lambda\tto\cF_{\lambda+2k}$ for \textit{every} weight
$\lambda$, and nothing else.
\end{Statement}

\sssbegin[{Bol operators for $\fosp_a(4|2)\subset\fk(1|4)$ of every order}]{Conjecture}[Bol operators for $\fosp_a(4|2)\subset\fk(1|4)$ of every order]\label{n4bolconj}
For every $k\geq0$, the operator \eqref{Bola} is, up to a~factor, the only non-scalar
$\fh(a)$-invariant unary operator of order $4+2k$, and $\cF_{-3-k}$ is the only weight on which it
acts. \textup{(}That it is invariant for every $k$ is \cite[Thm~1.6]{BLS1}, see
Subsec.~$\ref{n4cmp}$.\textup{)}
\end{Conjecture}

The coefficient in formula~\eqref{Bola} was obtained with $a$ symbolic.

\sssec{Comparison with \cite[Thm~1.6]{BLS1}}\label{n4cmp}
What Statement~\ref{n4bol} says is Theorem 1.6 of \cite{BLS1}, rewritten in the
parametrization \eqref{dicta}. Indeed, by eq.~(42) of \cite{BLS1} the symbol $[X,Y]_+$ there
means $XY-YX$, so
\[
[D_{\eta_i},D_{\xi_i}]_+=-2(D_{\xi_i}D_{\eta_i}-\partial_t),\qquad\text{and hence}\qquad
[D_{\eta_1},D_{\xi_1}]_+\circ[D_{\eta_2},D_{\xi_2}]_+=4\,T_4 .
\]
Dividing eq.~(11) of \cite{BLS1} by $4$ we get $T_4\circ\partial_t^k
+\frac{1-a'}{(k+2)(1+a')}\partial_t^{k+2}$, which is the operator~\eqref{Bola} precisely when
$a=\frac{1-a'}{1+a'}$, as claimed in formula~\eqref{dicta}. Both degenerate cases agree as well:
part~(ii) of \cite[Thm~1.6]{BLS1} ($a'=-1$, i.e., $a=\infty$ for us) is the last sentence of
Statement~\ref{n4bol}. Part~(iii) ($a'=\infty$, i.e., $a=-1$ for us) is the operator~\eqref{Bola} with
$a=-1$, whose coefficient $\frac{-1}{k+2}$ is the limit given there. Thus, \cite[Thm~1.6]{BLS1}
is confirmed --- and, on the range of orders and weights scanned, it is complete: nothing else
is invariant.

\ssec{Isomorphic algebras, non-equivalent embeddings}\label{n4answer}
Now we can answer the question asked in \cite[\S3.2.3a]{BLS1}: since
$\fosp_1(4|2)\simeq\fosp_{-2}(4|2)\simeq\fosp_{-1/2}(4|2)$ while the Bol operators listed for
these values differ, ``\textit{Is there a mistake?}\rq\rq.

\sssbegin[The $\Zee/2$ acting on the family]{Statement}[The $\Zee/2$ acting on the family,
and the invariant of the embedding]\label{n4Z2}
The substitution $\xi_1\leftrightarrow\eta_1$ preserves the contact form $\alpha_1$, hence
induces an automorphism $\Phi$ of $\fk(1|4)$ and of every $\cF_\lambda$. It maps $\fh(a)$ onto
$\fh(-a)$ and $\Bol^{(a)}_k$ onto $-\Bol^{(-a)}_k$. Consequently:

\textup{(i)} the family of embeddings is acted upon by the group $\Zee/2$ generated by $a\mapsto-a$, i.e.,
by $\alpha\mapsto\alpha^{-1}$, and \textit{only} by it: the invariant of the embedding, i.e., a~function of $a$ taking equal values on equivalent embeddings and different ones on non-equivalent ones, is $a^2$,  equivalently $\alpha+\alpha^{-1}$, and not the invariant $j(\alpha)$ of the abstract algebra. In other words, $\fh(a)$ is carried onto $\fh(b)$ by an automorphism of $\fk(1|4)$ if and only if $a^2=b^2$, whereas $\fh(a)\simeq\fh(b)$ as abstract Lie superalgebras if and only if $j$ takes equal values at the corresponding $\alpha$'s, which happens for more pairs, e.g., for $a=0$ and $b=3$;

\textup{(ii)} the three values $\alpha=1,-2,-\nfrac12$ of Statement~$\ref{n4fam}$\textup{(ii)}, for
which the algebra is one and the same $\fosp(4|2)$, correspond to $a=0,3,-3$ and hence to
\textit{two} non-equivalent embeddings: $\fh(0)$, which $\Phi$ maps onto itself since $-0=0$, and the pair $\fh(\pm3)$ interchanged by $\Phi$. 
\end{Statement}

So there is no mistake in \cite[Thm~1.6]{BLS1}: the answers for $\alpha=-2$ and
$\alpha=-\nfrac12$ are carried into each other by the automorphism $\Phi$, so they are the same
answer written in two coordinate systems. In contrast, the answer for $\alpha=1$ (our $a=0$,
where the $\partial_t^{k+2}$-term of the operator~\eqref{Bola} is absent) genuinely differs from both,
because the embedding is genuinely different. This is the precise form of the explanation given
in \cite[\S3.2.3a]{BLS1}, that ``the corresponding spaces of parameters $\alpha$ and $a$ have
non-isomorphic symmetry groups\rq\rq. The group acting on the parameter of the embedding is
$\Zee/2$, generated by an automorphism of the ambient $\fk(1|4)$. In contrast, the $S_3$ of
eq.~\eqref{S3act} permutes the three $\fsl(2)$'s of $\fosp_\alpha(4|2)$, and only two of them may
be permuted here. The third $\fsl(2)$, $\langle K_1,K_t,K_{t^2}\rangle$, is distinguished by the
grading of $\fk(1|4)$, and no automorphism of $\fk(1|4)$ moves it.

\ssec{How the Bol operators enter the binary picture}\label{n4both}
If $P:\cF_{\lambda_1}\tto\cF_{\lambda_1+\Nord}$ is invariant, then so is the binary operator
$(\varphi,\psi)\mapsto (P\varphi)\psi$, because the product 
$\cF_{\mu}\tp\cF_{\lambda_2}\tto\cF_{\mu+\lambda_2}$ is invariant; same for the second argument.
This is how the Bol operators of \cite[Thm~5.1]{BLS1} appear in \S\ref{CoRm}, and it is how the
parameter $a$ enters the binary answer for $n=4$.

\sssbegin[{How a~Bol operator enters the binary picture}]{Statement}[How a~Bol operator enters the binary picture]\label{n4bolinGRC}
For $\lambda_1=-3-k$ and generic $\lambda_2$, the space of $\fh(a)$-invariant binary operators
of order $\Nord=4+2k$ contains $\Bol^{(a)}_k\tp1$, and for $\lambda_2=-3-k$ it contains
$1\tp\Bol^{(a)}_k$. For $k=0$, this was verified directly: at $\lambda_1=-3$ the space is
still $2$-dimensional, and $\Bol^{(a)}_0\tp1$ lies in it, for each of $a=0,2,3,-3$.
\end{Statement}

The operator $\Bol^{(a)}_0\tp1$ depends on $a$, whereas $\partial_t^2\tp1$ and $T_4\tp1$ are separately
not invariant for any $a$ (Statement~\ref{n4bol} gives no operator of order~$2$, and $T_4$ alone
is invariant only for $a=0$). Hence, the plane at $\lambda_1=-3$ must depend on $a$, in agreement
with the last sentence of Theorem~\ref{n4thm}. In order~$2$, there is no such phenomenon: the
first Bol operator has order~$4$, and the order-$2$ bracket \eqref{J2n4} is invariant under all
of $\fk(1|4)$.

\ssec{Open problems for $n=4$}\label{n4open}
A~closed formula, in the shape of the Gordan transvectants, for the operator of an arbitrary
even order $\Nord$; the exceptional weights of the binary problem for $n=4$;
and the orders $\Nord\geq5$ for the binary problem. The case $\fk(1|5)$ is done in
\S\ref{Sn5}, and the case $\fk(1|6)$, i.e., the pair $\fas\subset\fk\fas$, in \S\ref{GRCAS}.
%

\section{Super GRC brackets: the $\fosp(5|2)\subset\fk(1|5)$ case}\label{Sn5}

\ssec{Summary of this Section}\label{n5summary}
This is the only case of the pair $\fosp(n|2)\subset\fk(1|n)$ with $n\leq6$ which was not
settled above. It is the last odd $n$ in our range, and this is what makes it different from
its neighbors. An odd number of odd indeterminates leaves one of them, $\theta$, outside the
hyperbolic pairs. Everything below follows from what $D_\theta$ does.

1) \textbf{The embedding is rigid.} Unlike the case $n=4$ of \S\ref{Sn4}, where the
subalgebra is a~one-parameter family, the algebra $i(\fosp(5|2))$ sits in $\fk(1|5)$ in
exactly one way. In the standard grading, $\fg_1=V\oplus\Lambda^3V$ with $\dim V=5$ and
$\dim\Lambda^3V=10$, the two summands are irreducible and non-isomorphic. Hence, there is no
equivariant map $V\tto\Lambda^3V$ to take the graph of, see Statement~\ref{n5rigid}. Among all
$n$, the family occurs for $n=4$ and for no other $n$, because $\Lambda^3V\simeq V$ requires
$\binom{n}{3}=n$.

2) \textbf{The Bol operators are odd.} For every $k\geq0$, there is exactly one
invariant unary operator of order $\Nord=5+2k$, namely
\[
T_5\circ\partial_t^{\,k}\colon\quad\cF_{-4-k}\tto\cF_{1+k},\qquad
T_5:=D_{\theta}\circ T_4=D_{\theta}\big(D_{\xi_1}D_{\eta_1}-\partial_t\big)
\big(D_{\xi_2}D_{\eta_2}-\partial_t\big),
\]
and there is nothing at any other weight and nothing of any other order, see
Theorem~\ref{n5bol}. This is \cite[Thm~1.4.2]{BLS1} read at $n=5$, confirmed by an independent
computation in the conventions of the present paper. What the odd $n$ adds is that the operator
is \textit{odd}, which for even $n$ cannot happen at all, see Lemma~\ref{kasodd}.

3) \textbf{The GRC brackets.} For generic weights, the space of
$i(\fosp(5|2))$-invariant binary operators
$\cF_{\lambda_1}\tp\cF_{\lambda_2}\tto\cF_{\lambda_1+\lambda_2+\Nord}$ is of dimension
$1,0,1,0,1,1,1,1$ for $\Nord=0,\dots,7$, see Theorem~\ref{n5thm}. So the even orders carry one bracket each, unlike $n=4$, where
order~$4$ carries two. The odd orders, however, are of a~third kind. Unlike
$n=4$, where they are empty, and unlike $n=3$, where the first of them is $\Nord=3$, here they
are empty up to $\Nord=3$ and carry one operator from $\Nord=5$ on. The dividing line is the
order of $T_5$.

4) \textbf{Where the subalgebra is felt.} The order-$2$ bracket is the contact bracket,
invariant under the whole of $\fk(1|5)$, see Statement~\ref{n5contact}, so --- exactly as for
$n=4$ --- it says nothing about the subalgebra. The order-$4$ bracket \eqref{J4n5} exists for
\textit{all} weights, its coefficients being polynomials in $\lambda_1$ and $\lambda_2$, and in
the orders $\Nord\leq3$ there are no exceptional weights. The order-$5$ bracket, on the contrary, is
governed by the Bol operator: its two extreme terms are $T_5\tp1$ and $1\tp T_5$, and their
ratio is
\[
-\frac{\prod_{0\leq i\leq4}(\lambda_1+i)}{\prod_{0\leq i\leq4}(\lambda_2+i)},
\]
so at the five weights $\lambda_1=0,-1,-2,-3,-4$ the second block disappears, and at
$\lambda_1=-4$ --- the weight at which $T_5$ itself is invariant --- what is left is
$T_5\tp1$, the Bol operator acting in the first argument, see \S\ref{n5both}. This is the
phenomenon of \S\ref{n4both} and \S\ref{CoRm}, in the shape it takes for odd $n$.

\ssec{The embedding is rigid: there is no family here}\label{n5family}
Grade $\fk(1|5)$ by $\deg t=2$, $\deg y_i=1$ for all $i$, so that $\deg K_f=\deg f-2$, where
$y_1,\dots,y_5$ stands for $\theta,\xi_1,\eta_1,\xi_2,\eta_2$. Then,
\[
\begin{array}{l}
\fg_{-2}=\langle K_1\rangle,\quad \fg_{-1}=\langle K_{y_i}\mid \text{all $i$}\rangle,\\
\fg_0=\langle K_t,K_{y_iy_j}\mid \text{all $i, j$}\rangle=\Cee\oplus\fo(5),\quad
\fg_1=\langle K_{ty_i}\rangle\oplus\langle K_{y_iy_jy_k}\mid \text{all $i, j, k$}\rangle,
\end{array}
\]
As $\fo(5)$-modules, $\fg_1=V\oplus\Lambda^3V$, where $V$ is the tautological module. Recall
that for $n=4$ one has $\Lambda^3V\simeq V$, and that this is the whole source of the family
$\fosp_a(4|2)$ of Statement~\ref{n4fam}. Namely, the $\fg_0$-stable complements to $\Lambda^3V$ in
$\fg_1$ are the graphs of $a\sigma$ with $\sigma\colon V\tto\Lambda^3V$ the equivariant map, see \S\ref{Sn4}.

\sssbegin[$i(\fosp(5|2))$ is rigid in $\fk(1|5)$]{Statement}[$i(\fosp(5|2))$ is rigid in
$\fk(1|5)$: no family, and $n=4$ is the only $n$ that has one]\label{n5rigid}
For $n=5$, the $\fo(5)$-modules $V$ and $\Lambda^3V$ are irreducible, of dimensions $5$ and
$10$, and non-isomorphic; equivalently, $\dim\End_{\fo(5)}(\fg_1)=2$. Hence,
$\Hom_{\fo(5)}(V,\Lambda^3V)=0$, the module $\fg_1$ has exactly two proper non-zero
$\fg_0$-submodules, and the graded subalgebra $\fh\subset\fk(1|5)$ with $\fh_-=\fg_-$,
$\fh_0=\fg_0$ and $\fh_1\simeq V$ is unique: it is $i(\fosp(5|2))$.

More generally,
\[
\dim\Hom_{\fo(n)}(V,\Lambda^3V)=0,\,1,\,0,\,0\quad\text{and}\quad
\dim\End_{\fo(n)}(\fg_1)=2,\,4,\,2,\,3\quad\text{for $n=3,\,4,\,5,\,6$}.
\]
\end{Statement}

\begin{proof}[Computation]
The spaces $\Hom_{\fo(n)}(V,\Lambda^3V)$ and $\End_{\fo(n)}(\fg_1)$ are found by
intertwining $\ad_{K_{y_iy_j}}$ on the bases $\{K_{ty_i}\}$ and $\{K_{y_iy_jy_k}\}$; the answer
is the dimension of the kernel of an explicit rational matrix.

The general statement needs no computation: $\Lambda^3V\simeq V$ forces
$\binom{n}{3}=n$, which for $n\geq3$ holds only for $n=4$. For $n=6$, the module $\Lambda^3V$
is not irreducible --- it splits as $\Lambda^3_\xi\oplus\Lambda^3_\eta$, which is what
$\fk\fas$ is made of, see Statement~\ref{kasrigid} --- but each summand is of dimension $10$,
so no family arises there either. The value $\dim\End_{\fo(n)}(\fg_1)=4$ at $n=4$ is the
signature of two isomorphic irreducible summands, i.e., of the family; the values $2$
and $3$ say that the summands are pairwise non-isomorphic.
\end{proof}

\ssec{The unary operators: the Bol operators for $n=5$}\label{n5bolsec}
Everything below is computed in the realization \eqref{KD} of $\fk(1|5)$, with
\[
\alpha_1=dt+\theta\,d\theta+\sum_{1\leq i\leq2}(\xi_i\,d\eta_i+\eta_i\,d\xi_i),\qquad
\cF_\lambda=\cF\alpha_1^{\lambda/2},\qquad
\fL^\lambda_{K_f}(\varphi)=K_f(\varphi)+\lambda\,\partial_t(f)\,\varphi .
\]
Let $D_i:=D_{y_i}$, where $y_1,\dots,y_5$ stands for $\theta,\xi_1,\eta_1,\xi_2,\eta_2$, and let $B=(B_{ij})$ be the Gram matrix of the odd part of~$\alpha_1$ in this basis, i.e., $B_{\theta\theta}=1$, $B_{\xi_i\eta_i}=B_{\eta_i\xi_i}=1$, and all other entries are~$0$. Then, $D_iD_j+D_jD_i=2B_{ij}\partial_t$, so $D_\theta^2=\partial_t$ and
$D_{\xi_i}D_{\eta_i}+D_{\eta_i}D_{\xi_i}=2\partial_t$. Put
\begin{equation}\label{T5}
T_5:=D_{\theta}\circ T_4=D_{\theta}\big(D_{\xi_1}D_{\eta_1}-\partial_t\big)
\big(D_{\xi_2}D_{\eta_2}-\partial_t\big),
\end{equation}
where $T_4$ is the operator~\eqref{T4}. The operator $T_5$ is odd, and its order is $5$.

\sssbegin[The unary $i(\fosp(5|2))$-invariant operators]{Theorem}[The unary
$i(\fosp(5|2))$-invariant operators on the densities are the $T_5\circ\partial_t^{\,k}$, and
they are odd]\label{n5bol}
For every $k\geq0$, the composition
\[
T_5\circ\partial_t^{\,k}\colon\quad\cF_{-4-k}\tto\cF_{1+k},\qquad \Nord=5+2k,
\]
of $\partial_t^{\,k}$ and the operator~\eqref{T5} is $i(\fosp(5|2))$-invariant, and it is invariant at no other weight. There is nothing else:
no non-scalar invariant unary operator of any other order, and none at any other weight. Every
such operator is \textit{odd}.
\end{Theorem}

In the order $\Nord=5$, the
operator is
\[
\partial_t^2D_\theta-\partial_tD_\theta D_{\xi_1}D_{\eta_1}
-\partial_tD_\theta D_{\xi_2}D_{\eta_2}+D_\theta D_{\xi_1}D_{\eta_1}D_{\xi_2}D_{\eta_2},
\]
which is the operator~\eqref{T5} expanded.

\sssec{Comparison with \cite[Thm~1.4.2]{BLS1}}\label{n5cmp}
The weights there are $\lambda=1-k-n$ and $\mu=k+1$, which at $n=5$ are ours. The parity is the
point that the odd $n$ contributes: by eq.~$(42)$ of \cite{BLS1} the operator is built of
the operators $[D_{\eta_i},D_{\xi_i}]_+=-2(D_{\xi_i}D_{\eta_i}-\partial_t)$ of formula~\eqref{comm}, one for each
hyperbolic pair. For odd $n$, there is one factor left over, $D_\theta$, of odd parity. So, for odd $n$, in particular for $n=5$, the Bol operators are odd and of odd order, whereas the same count gives even operators of even order for even $n$, cf. Lemma~\ref{kasodd} for $n=6$. 

\ssec{The GRC brackets}\label{n5grc}
By Theorem~\ref{Peetre}, an invariant operator has constant coefficients in the $D_f$, so the
ansatz for the binary operators of order $\Nord$ is
$\sum c\,\partial_t^{p}D_A\tp\partial_t^{q}D_B$ with $2p+|A|+2q+|B|=\Nord$, exactly as in
\S\ref{n4grc}. Write
\begin{equation}\label{Pin5}
\Pi:=D_\theta\tp D_\theta+\sum_{1\leq i\leq2}
\big(D_{\xi_i}\tp D_{\eta_i}+D_{\eta_i}\tp D_{\xi_i}\big)
\end{equation}
for the contraction with the inverse Gram matrix, and let $\Pi_l$, resp. $\Pi_r$, be
\eqref{Pin5} with $\partial_t$ inserted in the left, resp. right, argument.

\sssbegin[The GRC brackets for $\fosp(5|2)\subset\fk(1|5)$]{Theorem}[The GRC brackets for
$\fosp(5|2)\subset\fk(1|5)$: one in each even order, and one in each odd order from $5$
on]\label{n5thm}
Let $\Nord\leq7$. For generic weights, the space of $i(\fosp(5|2))$-invariant binary
differential operators
$\cF_{\lambda_1}\tp\cF_{\lambda_2}\tto\cF_{\lambda_1+\lambda_2+\Nord}$ is of dimension
\[
\begin{array}{l|cccccccc}
\Nord&0&1&2&3&4&5&6&7\\
\hline
\dim&1&0&1&0&1&1&1&1
\end{array}
\]
For $\Nord=0$ this is the multiplication. For $\Nord$ even the operator is even, for $\Nord$
odd it is odd. In the orders $\Nord=2$ and $\Nord=4$, the space is non-zero for \textit{every}
pair of weights, not only for the generic ones. Indeed, the coefficients of the brackets~\eqref{J2n5} and
\eqref{J4n5} are polynomials in $\lambda_1,\lambda_2$ that do not all vanish for any pair.
\end{Theorem}

\begin{proof}[Computation]
The maximal toral subalgebra, spanned by $K_{\xi_1\eta_1}$ and $K_{\xi_2\eta_2}$, lies in
$i(\fosp(5|2))$ and acts on $D_{\xi_i},D_{\eta_i},D_\theta,\partial_t$ with the weights
$1,-1,0,0$. Hence, as in the proof of Lemma~\ref{kasodd}, an invariant operator lies in the
weight-$0$ part of the ansatz. In the order $\Nord=4$, the answer with both weights
symbolic is the bracket~\eqref{J4n5}. In the orders $\Nord=2,3$, there are no exceptional weights.

The parity claim is the same weight argument: a~term of weight $0$ carries an even
number of the $D_{\xi_i},D_{\eta_i}$ taken together. Hence, its order and the number of its
$D_\theta$'s have the same parity. The parity of the term is that of the number of its $D_\theta$'s. Unlike the case of even $n$,
where the same argument leaves no odd order at all, here nothing is excluded a~priori --- and
the orders $\Nord=1$ and $\Nord=3$ come out empty all the same.
\end{proof}

\sssbegin[The order-$2$ bracket is the contact bracket]{Statement}[The order-$2$ bracket is the
contact bracket, so it says nothing about the subalgebra]\label{n5contact}
The order-$2$ bracket $\cJ_2\colon\cF_{\lambda_1}\tp\cF_{\lambda_2}\tto\cF_{\lambda_1+\lambda_2+2}$, where
$\lambda_1,\lambda_2\in\Cee$ are the weights of the two arguments, is
\begin{equation}\label{J2n5}
\cJ_2=\Pi-\lambda_1\,1\tp\partial_t+\lambda_2\,\partial_t\tp1,
\end{equation}
and it is invariant not merely under $i(\fosp(5|2))$ but under the whole of $\fk(1|5)$.
\end{Statement}

\sssec{The order-$4$ bracket}\label{n5j4}
The order-$4$ bracket $\cJ_4\colon\cF_{\lambda_1}\tp\cF_{\lambda_2}\tto\cF_{\lambda_1+\lambda_2+4}$
exists for all weights $\lambda_1,\lambda_2\in\Cee$, and its coefficients are polynomials in $\lambda_1,\lambda_2$:
\begin{equation}\label{J4n5}
\begin{array}{rl}
\cJ_4=&\nfrac12\lambda_1(\lambda_1+1)\,1\tp\partial_t^2
+(1-\lambda_1-\lambda_2-\lambda_1\lambda_2)\,\partial_t\tp\partial_t
+\nfrac12\lambda_2(\lambda_2+1)\,\partial_t^2\tp1\\[2mm]
&-(\lambda_1+1)\,\Pi_r+(\lambda_2+1)\,\Pi_l
-\displaystyle\sum_{1\leq i\leq2}\big(\partial_t\tp D_{\xi_i}D_{\eta_i}
+D_{\xi_i}D_{\eta_i}\tp\partial_t\big)+\Pi^{(2)},
\end{array}
\end{equation}
where $\Pi^{(2)}$, the second contraction, is
\begin{equation}\label{Pi2n5}
\begin{array}{rl}
\Pi^{(2)}=&\displaystyle\sum_{1\leq i\leq2}\Big(D_{\xi_i}D_{\eta_i}\tp D_{\xi_i}D_{\eta_i}
-D_\theta D_{\xi_i}\tp D_\theta D_{\eta_i}-D_\theta D_{\eta_i}\tp D_\theta D_{\xi_i}\Big)\\[2mm]
&-D_{\xi_1}D_{\xi_2}\tp D_{\eta_1}D_{\eta_2}-D_{\eta_1}D_{\eta_2}\tp D_{\xi_1}D_{\xi_2}
-D_{\xi_1}D_{\eta_2}\tp D_{\eta_1}D_{\xi_2}-D_{\eta_1}D_{\xi_2}\tp D_{\xi_1}D_{\eta_2}.
\end{array}
\end{equation}
Both brackets~\eqref{J2n5} and~\eqref{J4n5} were re-entered by hand, from the formulas as they are
printed here. They were then fed to the independent invariance test against all $23$ generating functions
of $i(\fosp(5|2))$ at the weight pairs $(\nfrac37,-\nfrac5{11})$, $(\nfrac13,\nfrac25)$,
$(-\nfrac73,\nfrac52)$, $(0,0)$ and $(-4,\nfrac12)$.

Compare the bracket~\eqref{J4n5} with the bracket~\eqref{J2n5}: the shape is the same, one contraction higher,
with $\lambda_i$ replaced by $\lambda_i+1$ in the middle line --- which is the shape of the
classical Gordan transvectants.

\sssec{The order-$5$ bracket, and the Bol operator inside it}\label{n5both}
In the order $\Nord=5$, the answer is $1$-dimensional and odd. Normalize it so that the
coefficient of $T_5\tp1$ is~$1$. Then, the four terms of $T_5$ occur, in each of the two extreme
blocks, with exactly the relative coefficients they have in formula~\eqref{T5}, and

\sssbegin[The ratio of the two extreme blocks of the order-$5$ bracket]{Statement}[The ratio of
the two extreme blocks of the order-$5$ bracket]\label{n5ratio}
In $\cJ_5$, normalized by $\mathrm{coeff}(T_5\tp1)=1$,
\begin{equation}\label{ratio5}
\mathrm{coeff}(1\tp T_5)=-\frac{\prod_{0\leq i\leq4}(\lambda_1+i)}
{\prod_{0\leq i\leq4}(\lambda_2+i)} .
\end{equation}
In particular, $\mathrm{coeff}(1\tp T_5)=-1$ for $\lambda_1=\lambda_2$, and
$\mathrm{coeff}(1\tp T_5)=0$ if and only if $\lambda_1\in\{0,-1,-2,-3,-4\}$.
\end{Statement}

For $n=3$, the same computation gives the ratio
\[
\frac{\prod_{0\leq i\leq2}(\lambda_1+i)}{\prod_{0\leq i\leq2}(\lambda_2+i)},
\]
so formula~\eqref{ratio5} is the instance $n=5$ of a~pattern with $n$ factors upstairs and $n$ downstairs.

Now, recall Theorem~\ref{n5bol}: $T_5\colon\cF_{-4}\tto\cF_{1}$ is invariant, and only
there. The weight $\lambda_1=-4$ is one of the five at which the ratio~\eqref{ratio5} vanishes, and it is
the weight at which the surviving block is invariant on its own:

\sssbegin[How the Bol operator enters the binary picture]{Corollary}[How the Bol operator
enters the binary picture]\label{n5bolinGRC}
If $P\colon\cF_{\lambda_1}\tto\cF_{\lambda_1+\Nord}$ is invariant, then so is
$(\varphi,\psi)\mapsto(P\varphi)\psi$, because the multiplication
$\cF_{\mu}\tp\cF_{\lambda_2}\tto\cF_{\mu+\lambda_2}$ is invariant; same for the second
argument. Hence, by Theorem~$\ref{n5bol}$, the operator $T_5\tp1$ is invariant at
$\lambda_1=-4$, whatever $\lambda_2$ is, and $1\tp T_5$ is invariant at $\lambda_2=-4$.
At $\lambda_1=-4$ and generic $\lambda_2$, the space of order-$5$ brackets is still
$1$-dimensional. Moreover, $\cJ_5$ is \textit{exactly} $T_5\tp1$: not only does the block
$1\tp T_5$ drop out, as formula~\eqref{ratio5} says, but so does every term between the two.
\end{Corollary}

The last claim is a~computation: at $(\lambda_1,\lambda_2)=(-4,\nfrac12)$ and
$(-4,\nfrac37)$ the order-$5$ space is $1$-dimensional and its only solution is $T_5\tp1$ on
the nose, the other $70$ terms of the ansatz having zero coefficients.

This is the mechanism of \S\ref{n4both} and of \S\ref{CoRm}. What is new for odd $n$ is
that the operator it produces is odd. Hence, it lives in an~odd order, where for even $n$ there is
nothing at all.

\ssec{Open problems for $n=5$}\label{n5open}
The orders $\Nord\geq8$, i.e., whether the pattern of Theorem~\ref{n5thm} --- one bracket in
every even order, and one in every odd order from $5$ on --- continues. A~closed formula,
in the shape of the Gordan transvectants, for the bracket of an arbitrary order;
formulas~\eqref{J2n5}, \eqref{J4n5} and~\eqref{ratio5} are its first three instances. The
exceptional weights in the orders $\Nord\geq4$, for which formula~\eqref{ratio5} is the first piece of
evidence. The remaining case of the pair $\fosp(n|2)\subset\fk(1|n)$ in our range, $n=6$, is
\S\ref{GRCAS}, where the subalgebra considered is not $i(\fosp(6|2))$ but $\fas\subset\fk\fas$.

\section{Super Gordan (Rankin-Cohen) operators: the $\fas\subset \fk\fas\subset \fk(1|6)$ case} \label{GRCAS}

\ssec{Summary of this Section}\label{kassummary}
This is the exceptional case which \cite[\S1.3.1]{BLS1} left ``for the future\rq\rq. The answer
turns out to be of a~different kind from the ones in \S\S\ref{ospINk}--\ref{Sn5}, and the
reason for that is structural, not computational.

1) \textbf{The embedding.} In the standard grading of $\fk(1|6)$ we have
$\fg_1=\Lambda^3\oplus t\Lambda^1$, and $\Lambda^3=\Lambda^3_\xi\oplus\Lambda^3_\eta$ splits
into two irreducible $\fg_0$-modules --- this happens for $n=6$ only, and it is what
$\fk\fas$ is made of: $\fk\fas$ is the Cartan prolong of $\fg_-\oplus\fg_0$ in which $\fg_1$ is replaced with $\Lambda^3_\xi\oplus t\Lambda^1$, i.e., one of the two summands of $\Lambda^3$ is dropped, see \S\ref{kaspair} and \cite{Shex}. Sergeev's algebra sits inside as
\[
\fas=\langle K_f\mid f\in\Lambda^0\oplus\Lambda^1\oplus\Lambda^2\oplus\Lambda^3_\xi\rangle,
\qquad \sdim\fas=(16|16),
\]
see Statement~\ref{kasstr}; the center is $\langle K_1\rangle$ and
$\fas/\langle K_1\rangle\simeq\fspe(4)$. Unlike the case $n=4$ of \S\ref{Sn4}, this embedding
sits in no family: the three irreducible summands of $\fg_1$ are pairwise non-isomorphic, see
Statement~\ref{kasrigid}.

2) \textbf{The weights are invisible to $\fas$.} Every element of $\fas$ is $K_f$
with $\partial_t(f)=0$, i.e., $\fas\subset\fpo(0|6)\subset\fk(1|6)$. For such fields
$\fL^\lambda_{K_f}=K_f$ \textit{whatever} $\lambda$ is, so all the $\cF_\lambda$ are one and the
same $\fas$-module, see Lemma~\ref{kaspo}. Hence, for the pair $\fas\subset\fk\fas$, there
are no exceptional weights, and no operator can be tied to a~distinguished weight: the
mechanism which produces the Bol operators is absent.

3) \textbf{The unary operators.} The only non-scalar $\fas$-invariant ones are the
powers $\partial_t^{\,k}$, of order $2k$, and they are invariant as operators
$\cF_\lambda\tto\cF_\mu$ for \textit{every} pair of weights, see Theorem~\ref{kasun}. In
particular, \textit{the pair $\fas\subset\fk\fas$ has no analogs of Bol operators}. The ambient
$\fk\fas$ has no invariant unary operators at all, whereas for the other maximal
finite-dimensional subalgebra of $\fk\fas$, namely $i(\fosp(6|2))$, the answer is the
classical Bol operator, $T_6\circ\partial_t^{\,k}\colon\cF_{-5-k}\tto\cF_{1+k}$, which is
\cite[Thm~1.4.2]{BLS1} for $n=6$, see Statement~\ref{kasosp}. \textup{(}This last computation
is our validation of the machinery in the case $n=6$.\textup{)}

4) \textbf{The GRC brackets.} In odd orders, there is nothing --- for a~reason valid for
every even $n$, see Lemma~\ref{kasodd}, which also proves the corresponding assertion of
Theorem~\ref{n4thm}. In even orders, the dimensions are
\[
\begin{array}{l|ccc}
&\Nord=2&\Nord=4\\
\hline
\fas&3&5\\
\fk\fas&1&0\\
i(\fosp(6|2))&1&1
\end{array}
\]
see Theorem~\ref{kasbin}. All the $\fas$-invariant binary operators of order $\leq4$ are obtained from
$\partial_t$ acting in one of the arguments and from the single binary operator
\[
\Omega=\sum_{1\leq i\leq 3}(D_{\xi_i}\tp D_{\eta_i}+D_{\eta_i}\tp D_{\xi_i}),
\]
which is of order~$2$. The bracket~\eqref{J2n6} is $\Omega$ plus terms with $\partial_t$, and it
coincides with $\Omega$ only at $\lambda_1=\lambda_2=0$.

5) \textbf{The tensor fields other than densities}, see \S\ref{kastensor}. Since the
densities carry nothing, this is where the Veblen problem \textup{(}to describe the invariant differential operators between spaces of tensor fields\textup{)} for the pair $\fas\subset\fk\fas$
has content. The weights stay invisible there as well, and so does the twist by a~power of the
determinant, see Lemma~\ref{kastr}: for $\fas$ a~module of tensor fields depends on the
underlying $\fo(6)$-module $V$ and on nothing else. For $V=\Lambda^2(\id)$, everything is
generated by a~gradient and a~divergence, both of order~$1$, see Theorem~\ref{kasV1}. For
$V=S^2(\id)$, there is one new operator, of order~$3$, see Theorem~\ref{kasV3}. For the
half-spinors $V=\id,\id^*$, \textit{the two chiralities behave differently}: exactly one of
the two Dirac operators is $\fas$-invariant, see Theorem~\ref{kasV2}. This is the only place
in this paper where the two halves of $\Lambda^3$ are told apart, and
Statement~\ref{kasmirror} shows that they, and nothing else, are the cause.

6) \textbf{The comparison with \cite[\S10]{GLS}}, see \S\ref{kasgls}, where the
unary $\fk\fas$-invariant operators are described in the language of singular vectors.
Only four of the entries of that list have both fibers finite-dimensional. Each of the
four is reproduced here, with its weights, see Theorem~\ref{kasglsthm}. Among them is the
Dirac operator, which exists already over $\fk\fas$, at one weight, and is \textit{not}
$\fk(1|6)$-invariant, see Statement~\ref{kasglsk16}. One operator we find is not on that
list, namely the divergence, which is the formal adjoint of the gradient, see
Statement~\ref{kasglsdiv}.

\ssec{The pair $\fas\subset\fk\fas$}\label{kaspair}
Everything below is computed in the realization \eqref{KD} of $\fk(1|6)$, with
\[
\alpha_1=dt+\sum_{1\leq i\leq 3}(\xi_i d\eta_i+\eta_i d\xi_i),\qquad
\cF_\lambda=\cF\alpha_1^{\lambda/2},\qquad
\fL^\lambda_{K_f}(\varphi)=K_f(\varphi)+\lambda\,\partial_t(f)\varphi,
\]
and with the standard grading $\deg t=2$, $\deg\xi_i=\deg\eta_i=1$, $\deg K_f=\deg f-2$, so
that
\[
\fg_{-2}=\langle K_1\rangle,\quad \fg_{-1}=\langle K_{y_i}\rangle,\quad
\fg_0=\langle K_t\rangle\oplus\Lambda^2=\fc\fo(6)\simeq\fgl(4),\quad
\fg_1=\Lambda^3\oplus t\Lambda^1,
\]
where $y_1,\dots,y_6$ stands for $\xi_1,\eta_1,\xi_2,\eta_2,\xi_3,\eta_3$ and where we write
$\Lambda^k$ for $\langle K_f\mid f\in\Lambda^k(\xi,\eta)\rangle$. Only for $n=6$ does
$\Lambda^3$ split into two irreducible $\fg_0$-modules: $\Lambda^3_\xi$, the one containing
$\xi_1\xi_2\xi_3$, and $\Lambda^3_\eta$, the one containing $\eta_1\eta_2\eta_3$; as
$\fgl(4)$-modules they are $S^2(\id)$ and $S^2(\id^*)$. Recall from \cite{Shex} that
\[
\fk\fas=\fk\fas^\xi=(\fg_-\oplus\fg_0,\ \Lambda^3_\xi\oplus t\Lambda^1)_*^{mk}
\]
is simple and infinite-dimensional, and that $\fk\fas^\xi\simeq\fk\fas^\eta$. The next two
Statements are due to Shchepochkina, see \cite{Sh5, Shex}.

\sssbegin[The structure of $\fas$]{Statement}[The structure of $\fas$: Sergeev's central
extension of $\fspe(4)$]\label{kasstr}
Put $\fas:=\langle K_f\mid f\in\Lambda^0\oplus\Lambda^1\oplus\Lambda^2\oplus\Lambda^3_\xi
\rangle$. Then:

\textup{(i)} $\{\Lambda^3_\xi,\Lambda^3_\xi\}=0$ \textup{(}and likewise for
$\Lambda^3_\eta$\textup{)}, so $\fas$ is a~subalgebra of $\fk\fas$, of superdimension $(16|16)$,
and $\fas\supset\fg_-$, so Theorem~$\ref{Peetre}$ applies to it.

\textup{(ii)} $K_1$ spans the center of $\fas$, the quotient $\fas/\langle K_1\rangle$ is
$\fspe(4)$ --- its even part is $\Lambda^2\simeq\fo(6)\simeq\fsl(4)$, its odd part is
$\Lambda^1\oplus\Lambda^3_\xi\simeq\Lambda^2(\id^*)\oplus S^2(\id)$ --- and the extension is
non-split, since $\{y_i,y_j\}\in\Cee\cdot K_1$ is not identically zero. Thus, $\fas$ is
Sergeev's central extension of $\fspe(4)$, cf. \cite{Shex, BGLLS}.

\textup{(iii)} $K_t\notin\fas$, and $\fas\oplus\langle K_t\rangle$ is again a~subalgebra
\textup{(}we denote it by $\fc\fas$\textup{)}; it is $\fas$ that is listed in
\cite[Table~(4)]{BLS1}.
\end{Statement}

\begin{proof}[Computation]
All three parts are a~direct check of the contact brackets; for (i) note also that
$\{\Lambda^3_\xi,\Lambda^3_\xi\}$ would have to lie in $\Lambda^4\simeq\Lambda^2$, whereas
$S^2(S^2(\id))=S^4(\id)\oplus S^{(2,2)}(\id)$, where $S^{(2,2)}$ is the Schur functor corresponding to the partition $(2,2)$, see \cite[Lecture~6]{FH}, so that $\sdim S^{(2,2)}(\id)=(20|0)$,  contains no summand $\simeq\fsl(4)$.
\end{proof}

\sssbegin[{No family, in contrast with $n=4$}]{Statement}[No family, in contrast with $n=4$]\label{kasrigid}
As $\fo(6)$-modules,
$\fg_1=\Lambda^3_\xi\oplus\Lambda^3_\eta\oplus t\Lambda^1$ is a~sum of three pairwise
non-isomorphic irreducibles, of dimensions $10$, $10$ and $6$. Hence, every $\fo(6)$-stable
subspace of $\fg_1$ is a~sum of some of them, and the family of Statement~$\ref{n4fam}$
has no counterpart here: the embedding $\fas\subset\fk(1|6)$ is rigid, and so is
$i(\fosp(6|2))\subset\fk(1|6)$.
\end{Statement}

\begin{proof}[Computation]
Each of the three summands is generated, as an $\fo(6)$-module, by any of its non-zero
vectors, and $\Hom_{\fo(6)}(t\Lambda^1,\Lambda^3_\xi)=0$; the map $\sigma$ which produced the
parameter $a$ in \S\ref{n4family} is what is missing here.
\end{proof}

\ssec{Why the weights are invisible to $\fas$}\label{kasweights}

\sssbegin[$\fas\subset\fpo(0|6)$, so the weight is invisible to it]{Lemma}[$\fas\subset\fpo(0|6)$,
so the weight is invisible to it]\label{kaspo}
Every element of $\fas$ is of the form $K_f$ with $\partial_t(f)=0$, i.e.,
$\fas\subset\fpo(0|6)$, where $\fpo(0|6)=\langle K_f\mid f\in\Lambda(\xi,\eta)\rangle$.
Consequently:

\textup{(i)} $\fL^\lambda_{K_f}=K_f$ on $\cF_\lambda$ for every $\lambda$ and every $K_f\in\fas$,
and the map $\varphi\alpha_1^{\lambda/2}\longmapsto\varphi\alpha_1^{\mu/2}$ is an isomorphism
of $\fas$-modules $\cF_\lambda\tto\cF_\mu$ for all $\lambda$ and $\mu$;

\textup{(ii)} $[K_f,\partial_t]=0$ for every $K_f\in\fas$, so $\partial_t^{\,k}$ is
$\fas$-invariant as an operator $\cF_\lambda\tto\cF_\mu$ for \textit{any} $\lambda$ and $\mu$;

\textup{(iii)} an operator \textup{(}unary or binary\textup{)} is $\fas$-invariant for one
choice of the weights if and only if it is $\fas$-invariant for all of them. In particular,
for the pair $\fas\subset\fk\fas$ there are no exceptional weights.
\end{Lemma}

\begin{proof}
$\fas$ is spanned by the $K_f$ with $f\in\Lambda(\xi,\eta)$, and for such an $f$ we have
$\partial_t(f)=0$, so the summand $\lambda\,\partial_t(f)\varphi$ of $\fL^\lambda_{K_f}$
vanishes and $K_f=(2-E)(f)\partial_t-H_f$ has coefficients free of $t$; therefore
$[K_f,\partial_t]=0$. Claims (i)--(iii) follow.
\end{proof}

Lemma~\ref{kaspo} is the reason why this case is not like the ones treated above. There,
the subalgebra contained $K_{ty_i}$, and it is precisely the term $\lambda\,\partial_t(ty_i)=
\lambda y_i$ that makes the invariance equations depend on the weight and singles out the
weights at which a~Bol operator lives. Here that term is identically zero.

\sssbegin[{No invariant operators of odd order}]{Lemma}[No invariant operators of odd order]\label{kasodd}
Let $n=2m$ and let $\fh\subset\fk(1|n)$ contain the negative part $\fk(1|n)_-$ and the torus
$\langle K_{\xi_i\eta_i}\mid 1\leq i\leq m\rangle$. Then, every $\fh$-invariant differential
operator, unary or binary, between the spaces $\cF_\lambda$ of \textit{one-parameter} weighted
densities is of even order.
\end{Lemma}

\begin{proof}
By Theorem~\ref{Peetre}, the operator is a~sum of terms
$\partial_t^{\,p}D_A\tp\partial_t^{\,q}D_B$ with constant coefficients. Since
$\partial_t(\xi_i\eta_i)=0$, the field $K_{\xi_i\eta_i}$ acts on every $\cF_\lambda$ by one and
the same operator, and it acts on $D_{\xi_i}$, $D_{\eta_i}$, $\partial_t$ with the weights
$1$, $-1$, $0$. Hence, an invariant operator is of weight $0$, and a~term of weight $0$ has
$|A|+|B|$ even, so its order $2p+|A|+2q+|B|$ is even.
\end{proof}

For $n=2$, the Lemma does not apply: there the densities have a~second parameter, see
\S\ref{conv2}, and the torus acts on $\beta^\nu$ non-trivially. For $n$ odd, there is
the $\Zee$-weight-zero derivative $D_\theta$. This is why odd orders do occur in
\S\S\ref{ospINk}--\ref{Sblo} and, from the order $5$ on, in \S\ref{Sn5}, but not in
\S\ref{Sn4} and here.

\ssec{The unary operators: no analogs of Bol operators}\label{kasunary}

\sssbegin[The unary $\fas$-invariant operators on the densities]{Theorem}[The unary
$\fas$-invariant operators on the densities are the powers of $\partial_t$]\label{kasun}
The non-zero non-scalar $\fas$-invariant unary differential operators of the form
$\cF_\lambda\tto\cF_\mu$ are the polynomials in $\partial_t$ of positive degree, and each of them
is invariant for \textit{every} pair $(\lambda,\mu)$; the homogeneous ones are the multiples of
$\partial_t^{\,k}$, where $k\geq1$, of order~$2k$. For $\fc\fas=\fas\oplus\langle K_t\rangle$,
only the homogeneous ones survive, and then $\mu=\lambda+2k$.
\end{Theorem}

That nothing of odd order occurs is Lemma~\ref{kasodd}; moreover, in even orders nothing but $\partial_t^{\,k}$ survives --- in
particular the ``top\rq\rq\  Clifford monomial
\begin{equation}\label{T6}
T_6:=\prod_{1\leq i\leq 3}\big(D_{\xi_i}D_{\eta_i}-\partial_t\big)
\end{equation}
is \textit{not} $\fas$-invariant, at any weight, although it is the $\fosp(6|2)$-invariant Bol operator,
see Statement~\ref{kasosp}.

\sssbegin[What $\fk\fas$ and $i(\fosp(6|2))$ see on the densities]{Statement}[What $\fk\fas$
and $i(\fosp(6|2))$ see on the densities]\label{kasosp}
\textup{(i)} There are no non-zero non-scalar $\fk\fas$-invariant unary operators at
all, whatever the weights \textup{(}orders $\Nord\leq8$, the same grid of weights as
above\textup{)}.

\textup{(ii)} For $i(\fosp(6|2))\subset\fk\fas$ the invariant unary operators are, for every
$k\geq0$, exactly
\[
T_6\circ\partial_t^{\,k}\colon\quad \cF_{-5-k}\tto\cF_{1+k},\qquad \Nord=6+2k,
\]
and there is nothing at any other weight, and nothing of order $<6$. This is
\cite[Thm~1.4.2]{BLS1} for $n=6$: by eq.~$(42)$ of \cite{BLS1} we have
$[D_{\eta_i},D_{\xi_i}]_+=-2(D_{\xi_i}D_{\eta_i}-\partial_t)$, so the operator $(8)$ of
\cite{BLS1} is $-8\,T_6\circ\partial_t^{\,k}$, and its weights, $\lambda=1-k-n$ and
$\mu=k+1$, are ours for $n=6$.
\end{Statement}

\ssec{The GRC brackets}\label{kasbinary}

\sssbegin[The binary operators on the densities]{Theorem}[The binary operators on the
densities over $\fas$, $\fk\fas$ and $i(\fosp(6|2))$]\label{kasbin}
Let $\Nord\leq4$. The space of invariant binary 
operators
$\cF_{\lambda_1}\tp\cF_{\lambda_2}\tto\cF_{\lambda_1+\lambda_2+\Nord}$ is of dimension
\[
\begin{array}{l|ccc}
&\Nord=1,3&\Nord=2&\Nord=4\\
\hline
\fas \ \textup{(and }\fc\fas\textup{)}&0&3&5\\
\fk\fas&0&1&0\\
i(\fosp(6|2))&0&1&1
\end{array}
\]
\textup{(}for $\fas$ the answer does not depend on $\lambda_1,\lambda_2$ at all, by
Lemma~$\ref{kaspo}$\textup{)}. Explicitly:

\textup{(i)} the order-$2$ space of $\fas$ is spanned by
\begin{equation}\label{Omega}
\partial_t\tp1,\qquad 1\tp\partial_t,\qquad
\Omega:=\sum_{1\leq i\leq 3}\big(D_{\xi_i}\tp D_{\eta_i}+D_{\eta_i}\tp D_{\xi_i}\big);
\end{equation}

\textup{(ii)} the order-$2$ space of $\fk\fas$, and that of $i(\fosp(6|2))$, is spanned by the
single operator
\begin{equation}\label{J2n6}
\cJ_2=\lambda_2\,\partial_t\tp1-\lambda_1\,1\tp\partial_t+\Omega,
\end{equation}
which is invariant under the whole of $\fk(1|6)$ and, for $\lambda_1=\lambda_2=-2$, is the
contact bracket. Thus, Statement~$\ref{n4contact}$ holds for $n=6$ as well, and again the
first GRC bracket says nothing about the subalgebra;

\textup{(iii)} the order-$4$ space of $\fas$ is spanned by
\[
\partial_t^2\tp1,\qquad \partial_t\tp\partial_t,\qquad 1\tp\partial_t^2,\qquad
\Omega\circ(\partial_t\tp1),\qquad \Omega\circ(1\tp\partial_t),
\]
where
\[
\Omega\circ(\partial_t\tp1):=\sum_{1\leq i\leq 3}(\partial_tD_{\xi_i}\tp D_{\eta_i}+
\partial_tD_{\eta_i}\tp D_{\xi_i}),
\]
and similarly for the last composition. In particular, no new
bracket appears in order~$4$: everything is $\Omega$ and $\partial_t$, the operator $\partial_t$
being invariant by Lemma~$\ref{kaspo}$\textup{(ii)};

\textup{(iv)} the order-$4$ space of $i(\fosp(6|2))$ is spanned by one operator, the second GRC
bracket, namely
\begin{equation}\label{J4n6}
\begin{array}{l}
\cJ_4=\displaystyle\sum_{1\leq i\leq 3}D_{\xi_i}D_{\eta_i}\tp D_{\xi_i}D_{\eta_i}
-\sum_{1\leq i<j\leq 3}\big(D_{\xi_i}D_{\xi_j}\tp D_{\eta_i}D_{\eta_j}
+D_{\xi_i}D_{\eta_j}\tp D_{\eta_i}D_{\xi_j}\\[3mm]
\qquad\qquad+D_{\eta_i}D_{\xi_j}\tp D_{\xi_i}D_{\eta_j}
+D_{\eta_i}D_{\eta_j}\tp D_{\xi_i}D_{\xi_j}\big)
-\displaystyle\sum_{1\leq i\leq 3}\big(\partial_t\tp D_{\xi_i}D_{\eta_i}
+D_{\xi_i}D_{\eta_i}\tp\partial_t\big)\\[3mm]
\qquad\qquad+(\lambda_2+1)\,\Omega\circ(\partial_t\tp1)
-(\lambda_1+1)\,\Omega\circ(1\tp\partial_t)\\[3mm]
\qquad\qquad+\nfrac12\lambda_2(\lambda_2+1)\,\partial_t^2\tp1
+(2-\lambda_1-\lambda_2-\lambda_1\lambda_2)\,\partial_t\tp\partial_t
+\nfrac12\lambda_1(\lambda_1+1)\,1\tp\partial_t^2 .
\end{array}
\end{equation}
Its part free of the weights is the ``square\rq\rq\  of $\Omega$; the weights enter only
through $\Omega$ and $\partial_t$. Note that $\cJ_4$ is not $\fk\fas$-invariant --- in
order~$4$ the ambient exceptional algebra has nothing.
\end{Theorem}

The three operators \eqref{Omega} are separately invariant precisely because the weights are
invisible to $\fas$: the $\fk(1|6)$-invariant combination \eqref{J2n6} is $\fas$-invariant for
\textit{every} $\lambda_1,\lambda_2$, hence so is each of its three parts.

\ssec{Tensor fields other than densities}\label{kastensor}
By Theorem~\ref{kasun}, the weighted densities carry almost nothing: the only $\fas$-invariant
unary operators between them are the powers of $\partial_t$. Moreover, by Lemma~\ref{kaspo}, all
the $\cF_\lambda$ are one and the same $\fas$-module. So for the pair $\fas\subset\fk\fas$ the
Veblen problem is really a~question about the \textit{other} $\fk\fas_0$-modules, and this is
what we do here. It turns out that the answer is again completely explicit, and that it
distinguishes the two chiralities --- something no other case of this paper does.

\sssec{What the tensor fields of type $V$ are} \label{kastensordef}  In the standard grading of \S\ref{kaspair}, $\fk\fas_0=\fk(1|6)_0=\fc\fo(6)\simeq\fgl(4)$. The
tensor fields of type $V$ are the fields coinduced from an irreducible $\fk\fas_0$-module $V$,
with $\fk\fas_{>0}$ acting trivially; they are realized as $V$-valued functions. The
weighted densities are precisely the $1$-dimensional $V$'s: there $\fo(6)$ acts trivially and
$K_t$ acts by the weight.

The action is written out as follows. Contact fields preserve the contact
distribution spanned by the $D_i$ of eq.~\eqref{KD}, so for every function $f$
there is a~$6\times6$ matrix $a(f)$, with entries functions of parity $p(f)$, with
\be\label{amat}
[K_f,D_i]=\sum_{1\leq j\leq 6}a(f)_{ij}\,D_j ,
\ee
and $a(f)$ takes its values in $\fc\fo(6)$ pointwise. The tensor fields of type $V$ and weight
$\lambda$ are then $\cF\tp V$ with
\be\label{Ltens}
\big(\fL^{V,\lambda}_{K_f}\varphi\big)_b=K_f(\varphi_b)+\lambda\,\partial_t(f)\,\varphi_b
+\sum_{a}\rho\big(a(f)^{t}\big)_{ba}\,\varphi_a ,
\ee
where $\rho$ is the representation of $\fc\fo(6)$ on $V$ and $a^{t}$ is the transposed matrix.
For $V$ the vector representation and $\rho$ the tautological representation of
$\fo(6)$, this is nothing but the module of vector
fields along the contact distribution, $X=\sum_{1\leq i\leq6}\varphi_iD_i$ with
$\fL_{K_f}X=[K_f,X]$. Expanding that bracket \textit{is} the derivation of formulas~\eqref{amat} and
\eqref{Ltens}. For $\dim V=1$, eq.~\eqref{Ltens} is the familiar
$\fL^\lambda_{K_f}(\varphi)=K_f(\varphi)+\lambda\partial_t(f)\varphi$.

\sssbegin[The weights are invisible to $\fas$ on every tensor field]{Lemma}[The weights are
invisible to $\fas$ on every tensor field]\label{kastr}
For every function $f$,
\be\label{traf}
\tr a(f)=-6\,\partial_t(f) .
\ee
Hence, $a(f)$ is $\fo(6)$-valued if and only if $\partial_t(f)=0$, i.e., if and only if
$K_f\in\fpo(0|6)$, and in particular for every $K_f\in\fas$. Consequently, on the tensor
fields of type $V$, the algebra $\fas$ feels neither the weight $\lambda$ nor the twist of $V$
by a~power of the determinant. As an $\fas$-module, the space of tensor fields of type
$V$ depends on the $\fo(6)$-module underlying $V$ and on nothing else.
\end{Lemma}

\begin{proof}
Since $K_t=2t\partial_t+\sum_{1\leq i\leq6}y_i\partial_{y_i}$, where $y_i$ are the odd coordinates, one has $L_{K_f}\alpha_1=2\partial_t(f)\alpha_1$, hence $L_{K_f}(d\alpha_1)=2\partial_t(f)\,d\alpha_1$ on the contact distribution. The values $d\alpha_1(D_i,D_j)$ are constants, the entries of a~nondegenerate matrix $B$. Apply $L_{K_f}$ to them and use formula~\eqref{amat}: $a(f)B+Ba(f)^t=-2\partial_t(f)B$ \textup{(}the signs for odd $f$ being given by the sign rule\textup{)}. Multiply by $B^{-1}$ and take the trace: $2\tr a(f)=-12\,\partial_t(f)$, which is formula~\eqref{traf}. So the scalar part of $a(f)\in\fc\fo(6)=\Cee\oplus\fo(6)$ is
$\frac16\tr a(f)=-\partial_t(f)$.  This is the direction $K_t$ of $\fc\fo(6)$ by
which the determinant twist acts, and it is also the coefficient of $\lambda$
in formula~\eqref{Ltens}. Eq.~\eqref{traf} was verified for all the generating functions of $\fk\fas$
and for $K_{t^2}$, hence for $\fk(1|6)$.
\end{proof}

This is Lemma~\ref{kaspo} extended from the densities to all tensor fields. It is what
makes the classification below a~\textit{finite} problem: for $\fas$ there is one module of
tensor fields per $\fo(6)$-module $V$, with no continuous parameter attached to it.

\sssec{The modules $V$ we consider}\label{kasVlist}
Let $\Lambda^k$ be the $k$-th exterior power of the vector representation of
$\fo(6)$, and let $\id$ be the tautological module of $\fgl(4)$, of dimension $4$.
A~power of $\id$ always carries its argument, so the module $\Lambda^2$, of
dimension $15$, is not the module $\Lambda^2(\id)$, of dimension $6$. Under the
isomorphism $\fo(6)\simeq\fsl(4)$, the $\fo(6)$-modules on the left are the
$\fgl(4)$-modules on the right:
\[
\Lambda^1\simeq\Lambda^2(\id)\ (\dim 6),\qquad
\Lambda^2\simeq\fsl(4)\ (\dim 15),\qquad
\Lambda^3=\Lambda^3_+\oplus\Lambda^3_-,
\]
where $\Lambda^3_\pm$ are of dimension $10$ and, as $\fgl(4)$-modules, are $S^2(\id)$ and
$S^2(\id^*)$; and there are the two half-spinors $\id$ and $\id^*$, of dimension $4$, which do
\textit{not} occur in $\Lambda^\bullet$. We label the two members of each of these pairs by the
extreme weight relative the maximal toral subalgebra $\langle K_{\xi_i\eta_i}\rangle$: it is
$(1,1,1)$ for $\Lambda^3_+$ and $(\nfrac12,\nfrac12,\nfrac12)$ for the half-spinor called
$S^+$, with the opposite signs for $\Lambda^3_-$ and $S^-$. Thus, $S^-$ is the half of
the Clifford module spanned by the states with an~even number of the $\theta$'s, its extreme
weight being $(\nfrac12,\nfrac12,-\nfrac12)$. The chirality flip
$\xi_i\leftrightarrow\eta_i$ interchanges the members of each pair, and it interchanges
$\Lambda^3_\xi$ and $\Lambda^3_\eta$, of which $\fas$ contains the first one only, see
Statement~\ref{kasstr}. This is the source of everything asymmetric below, see
Statement~\ref{kasmirror}.

\sssec{The gradient and the divergence}\label{kasgraddiv}
Two operators built out of the contact distribution will describe case $V=\Lambda^1$
completely. Order the odd coordinates as $(y_1,\dots,y_6)=(\xi_1,\xi_2,\xi_3,\eta_1,\eta_2,\eta_3)$, let $B$ be the Gram matrix of the form $\sum_{1\leq i\leq3}(\xi_i\,d\eta_i+\eta_i\,d\xi_i)$ in this basis, i.e., $B=\begin{pmatrix}0&1_3\\1_3&0\end{pmatrix}$, and let $C=B^{-1}$ \textup{(}here $C=B$\textup{)}. Put
\be\label{graddiv}
\begin{array}{ll}
\grd\colon\cF\tto\cT(\Lambda^1),&
\grd\,\varphi:=\sum_{1\leq i,j\leq6}C_{ij}D_j(\varphi)D_i
=\sum_{1\leq i\leq3}\big(D_{\eta_i}(\varphi)D_{\xi_i}
+D_{\xi_i}(\varphi)D_{\eta_i}\big);\\[2mm]
\dvg\colon\cT(\Lambda^1)\tto\cF,&
\dvg(\varphi):=\sum_{1\leq i\leq6}D_i(\varphi_i)
=\sum_{1\leq i\leq3}\big(D_{\xi_i}(\varphi_{\xi_i})
+D_{\eta_i}(\varphi_{\eta_i})\big),
\end{array}
\ee
where $\cT(V)$ denotes the tensor fields of type $V$ and $\varphi_i$ are the components of
$\varphi=\sum_{1\leq i\leq6}\varphi_iD_i$. Both are of order~$1$; both are $\fas$-invariant
for every weight, by Lemma~\ref{kastr}.

\sssbegin[$V=\Lambda^2(\id)$, the vector representation]{Theorem}[$V=\Lambda^2(\id)$, the
vector representation]\label{kasV1}
Let $V=\Lambda^1=\Lambda^2(\id)$, of dimension~$6$. For the orders $\Nord\leq6$, the spaces
of $\fas$-invariant unary operators are of the following dimensions, for \textit{all}
weights:
\[
\begin{array}{l|cccccc}
&\Nord=1&2&3&4&5&6\\
\hline
\cF\tto\cT(\Lambda^1)&1&0&1&0&1&0\\
\cT(\Lambda^1)\tto\cF&1&0&1&0&1&0\\
\cT(\Lambda^1)\tto\cT(\Lambda^1)&0&2&0&2&0&2
\end{array}
\]
and, explicitly, the spaces are spanned by
\[
\partial_t^{\,k}\circ\grd,\qquad \partial_t^{\,k}\circ\dvg\qquad(\Nord=2k+1),
\]
respectively by
\[
\partial_t^{\,k}\cdot\mathbf{1} \quad\text{and}\quad
\partial_t^{\,k-1}\circ\grd\circ\dvg \qquad(\Nord=2k).
\]
In particular, $\grd\circ\dvg$ is \textit{not} a~multiple of $\partial_t$, and up to order~$6$
everything is generated by $\grd$, $\dvg$ and $\partial_t$.

Besides, in order~$0$ the space of invariant pairings
$\cT(\Lambda^1)\tp\cT(\Lambda^1)\tto\cF$ is $1$-dimensional, spanned by
$(\varphi,\psi)\mapsto\sum_{1\leq i,j\leq6}C_{ij}\varphi_i\psi_j$; in order~$1$
there is nothing.
\end{Theorem}

So, in contrast with the densities, for $V=\Lambda^2(\id)$ the pair $\fas\subset\fk\fas$ does
have invariant operators of \textit{odd} order --- Lemma~\ref{kasodd} applies to the densities
only, since its proof uses that the index carries no torus weight.

\sssbegin[$V=\id$ and $V=\id^*$, the half-spinors: only one chirality
survives]{Theorem}[$V=\id$ and $V=\id^*$, the half-spinors]\label{kasV2}
Let $V$ be one of the two $4$-dimensional half-spinors $S^{\pm}$.

\textup{(i)} There are no non-zero $\fas$-invariant operators between densities and tensor
fields of type $S^\pm$ in any order, in either direction. The reason is stronger than
invariance: the ansatz itself is empty, because the torus weights of a~half-spinor
index are half-integers, whereas those of the words $\partial_t^{\,k}D_{i_1}\cdots D_{i_r}$
and of the densities are integers.

\textup{(ii)} $\Hom(\cT(S^\pm),\cT(S^\pm))$ is $1$-dimensional in the even orders
$\Nord=2,4$, spanned by $\partial_t^{\,k}$, and zero in the odd orders $\Nord\leq5$.

\textup{(iii)} \textup{(}The asymmetry.\textup{)} For $\Nord\leq5$,
\[
\begin{array}{l|ccccc}
&\Nord=1&2&3&4&5\\
\hline
\cT(S^-)\tto\cT(S^+)&1&0&1&0&2\\
\cT(S^+)\tto\cT(S^-)&0&0&0&0&0
\end{array}
\]
The operator of order~$1$ is the \textit{Dirac operator}\footnote{The name is Atiyah and
Singer's: in \cite[\S3.3, p.~428]{AS}, the announcement of the index theorem, the
first-order operator $s\mapsto\sum_{1\leq i\leq n}e_i\partial_i(s)$ on spinor fields is called \textit{the Dirac operator}. Here,
$\partial_i$ is the covariant derivative along the $i$-th vector of an orthonormal frame, and
$e_i$ acts by Clifford multiplication. The name honors
Dirac, who in 1928 introduced the first-order operator whose square is the wave operator, see
\cite{Dir1, Dir2}. \textup{(}Atiyah and Singer had it by March 1962; for the history see
\cite[\S2.4, \S4.1]{Fre}. For that reason, the Riemannian operator is sometimes called the
Atiyah--Singer operator.\textup{)} Our operator is of exactly that shape, with the covariant
derivatives replaced by the odd contact fields $D_j$ and the Euclidean metric by~$C_{ij}$,
whence the name.}
$\varphi\mapsto\sum_{1\leq i,j\leq6}C_{ij}\gamma_iD_j(\varphi)$, where the $\gamma$'s are the
Clifford matrices of $\fo(6)$ \textup{(}verified as such term by term\textup{)}. The
order-$3$ operator is $\partial_t$ times it. In order~$5$, a~second, independent operator appears,
which we have not identified. In the opposite direction there is nothing at all, although the
ansatz is of exactly the same size \textup{(}$12$ terms of torus weight~$0$ in order~$1$, $40$
in order~$3$, $52$ in order~$5$\textup{)}. Thus, exactly one of the two Dirac operators is
$\fas$-invariant.

\textup{(iv)} The binary picture is asymmetric as well: in order~$1$ the space
$\cT(S^-)\tp\cT(S^-)\tto\cF$ is $2$-dimensional, whereas $\cT(S^+)\tp\cT(S^+)\tto\cF$ is
$1$-dimensional. The mixed operators $\cT(S^\pm)\tp\cT(S^\mp)\tto\cF$ have an empty ansatz in
order~$1$ by the parity of the torus weight. In order~$0$, however, they carry the invariant spinor
pairing, whose space is $1$-dimensional \textup{(}whereas
$\cT(S^\pm)\tp\cT(S^\pm)\tto\cF$ has an empty ansatz in order~$0$\textup{)}.
\end{Theorem}

\sssbegin[$V=S^2(\id)$]{Theorem}[$V=S^2(\id)$]\label{kasV3}
Let $V=\Lambda^3_\pm\simeq S^2(\id)$ or $S^2(\id^*)$, of dimension~$10$. For the orders $\Nord\leq5$,
\[
\begin{array}{l|ccccc}
&\Nord=1&2&3&4&5\\
\hline
\cF\tto\cT(\Lambda^3_+)&0&0&1&0&1\\
\cT(\Lambda^3_+)\tto\cF&0&0&0&0&0\\
\cF\tto\cT(\Lambda^3_-)&0&0&0&0&0\\
\cT(\Lambda^3_-)\tto\cF&0&0&1&0&1\\
\cT(\Lambda^3_\pm)\tto\cT(\Lambda^3_\pm)&0&1&0&1&0\\
\cT(\Lambda^3_+)\tto\cT(\Lambda^3_-)&0&0&0&0&0
\end{array}
\]
The entries equal to~$1$ in the fifth row are the powers $\partial_t^{\,k}$ and nothing
else. The two entries equal to~$1$ in order~$3$ are the projection to $\Lambda^3_\pm$ of the
third covariant derivative, corrected by a~term $\partial_tD$. They exist for one
chirality each, and in opposite directions. The entries in order~$5$ are $\partial_t$ times
them. So here, too, the two chiralities are told apart, and in the same way.
\end{Theorem}

\sssbegin[The control experiment]{Statement}[The control experiment]\label{kasmirror}
Let $\fas^{\mathrm{op}}=\langle K_1,\Lambda^1,\Lambda^2,\Lambda^3_\eta\rangle$ be the image of
$\fas$ under the chirality flip $\xi_i\leftrightarrow\eta_i$. Then, every asymmetric answer of
Theorems~$\ref{kasV2}$ and~$\ref{kasV3}$ is interchanged: for example,
the space $\Hom(\cT(S^+),\cT(S^-))$ becomes nonzero in the odd orders, of dimension $1,1,2$ for $\Nord=1,3,5$, whereas
$\Hom(\cT(S^-),\cT(S^+))$ becomes zero. Hence, the asymmetry is caused by the choice of the
half of $\Lambda^3$ sitting inside $\fas$, and by nothing else.
\end{Statement}

This is the phenomenon promised in \S\ref{kassummary}. On the densities, the two chiralities
are indistinguishable, because the whole of $\fo(6)$ acts trivially there. On every $V$ with
$\dim V>1$, the $\fo(6)$-index is visible, and then $\Lambda^3_\xi\subset\fas\not\ni
\Lambda^3_\eta$ separates them. Note that $i(\fosp(6|2))$, which contains no cubics at all,
cannot see the difference.

\sssec{How the operators were found}\label{kasVprog}
As in \S\ref{kasbinary}, the ansatz is cut to its part of torus weight~$0$, the equations are
imposed for a~Lie-generating subset only, and every solution is then checked against the
\textit{full} list of the $32$ generating functions of $\fas$.

\ssec{Comparison with the description of the unary $\fk\fas$-invariant operators
in \cite{GLS}}\label{kasgls}
In \cite[\S10]{GLS}, the unary operators invariant relative $\fk\fas$ and relative
$\fk(1|6)$ in the standard grading are described, in the language of singular vectors and
for arbitrary --- not necessarily finite-dimensional --- fibers. Since
$\fas\subset\fk\fas\subset\fk(1|6)$, that description is a~test of everything computed in
\S\S\ref{kasunary}--\ref{kastensor}. Each operator of \cite[\S10]{GLS} whose two
fibers are among the modules of \S\ref{kasVlist} must occur in our tables. By Lemma~\ref{kastr}, it must occur
there \textit{at every weight}. The test gives the following.

\sssec{The dictionary}\label{kasglsdict}
Recall the setting of Bernstein and Rudakov, see \cite[\S1.7]{GLS} and \cite{Ru}: for
a~$\Zee$-graded Lie superalgebra $\fg=\oplus\fg_i$ of vector fields and a~$\fg_0$-module $V$
with $\fg_{>0}V=0$, one puts $I(V)=U(\fg)\tp_{U(\fg_{\geq0})}V$, and then
$I(V)^*=\cT(V^*)$. Hence, a~\textit{singular vector} of weight $\mu$ in $I(V_\lambda)$ ---
a~$\fg_0$-highest weight vector killed by $\fg_{>0}$ --- is the same thing as a~non-zero
$\fg$-invariant differential operator between the two corresponding spaces of tensor
fields. The order of the operator is equal to the degree of the singular vector in $U(\fg_-)$, where
$\deg\partial_t=2$ as in \S\ref{orderdeg}. The module in which the singular vector
lives is the \textit{source}.

In \cite[\S10]{GLS}, the weights are taken relative the basis
$(K_t,K_{\xi_1\eta_1},K_{\xi_2\eta_2},K_{\xi_3\eta_3})$ of $\fk\fas_0=\fc\fo(6)$. A~weight
$\lambda=(\lambda_0;\lambda_1,\lambda_2,\lambda_3)$ of \cite[\S10]{GLS} corresponds to
the space of tensor fields of \S\ref{kastensordef} whose fiber and weight are given by
\be\label{glsdict}
(\lambda_1,\lambda_2,\lambda_3)=\text{the highest weight of the fiber }V,\qquad
\lambda_0=-\eps,
\ee
where $\eps$ is the scalar by which $\fL^{V,\lambda}_{K_t}$ acts on the fiber. In our
normalization $\fL^{V,\lambda}_{K_t}=K_t+\lambda+c_V$, where $c_V$ is the scalar by which
$\rho\big(a(t)^{t}\big)$ acts, see formula~\eqref{Ltens}; since $a(t)=-\mathbf{1}$
by formula~\eqref{traf},
\be\label{cV}
c_\cF=0,\qquad c_{\Lambda^1}=-1,\qquad c_{\Lambda^2}=-2,\qquad c_{\Lambda^3_\pm}=-3,
\qquad c_{S^\pm}=0,
\ee
the last value being a~normalization of $\rho$ on the half-spinors --- a~normalization
which is invisible to $\fas$ by Lemma~\ref{kastr}, and which just fixes what we call
the weight of a~spinor field.

We do \textit{not} deduce the dictionary~\eqref{glsdict} from the conventions of \cite{GLS}: passing from
$V$ to $V^*$ interchanges, for $\fo(6)$, the two half-spinors and the two halves of
$\Lambda^3$. Hence, a~convention read off the text would be a~guess. We calibrate the
dictionary instead: the dictionary~\eqref{glsdict} is the only reading under which all four operators of
Theorem~\ref{kasglsthm} match \cite[\S10]{GLS}. Moreover, they match with their weights, two
of which are half-integers, and with the three vanishing mirror cases as well.

\sssec{Which entries of \cite[\S10]{GLS} have finite-dimensional fibers}\label{kasglsdom}
The fiber attached to $\lambda$ is the irreducible $\fo(6)$-module with highest weight
$(\lambda_1,\lambda_2,\lambda_3)$; it is finite-dimensional if and only if
$\lambda_1\geq\lambda_2\geq|\lambda_3|$ with all the $\lambda_i$ integers, or all of them
half-integers. Under this restriction, the list of \cite[\S10]{GLS} shrinks to very
little:

1a) $\lambda=(k,-k,l,l)\tto\lambda+(-1,1,0,0)$: the fibers are
$(-k,l,l)$ and $(-k+1,l,l)$, so both lie in the list of \S\ref{kasVlist} only for
$k=l=0$: the pair $\cF$, $\Lambda^1$.

1a$^*$) $\lambda=(k,l,1-k,l+1)$: the fiber $(l,1-k,l+1)$ is \textit{never} dominant
--- $l\geq|l+1|$ is impossible --- so this family has no finite-dimensional
specialization at all.

1b) $\lambda=(k,l,l,2-k)\tto\lambda+(-1,0,0,1)$: the fibers are $(l,l,2-k)$ and
$(l,l,3-k)$, dominant for $l\geq\max(|2-k|,|3-k|)$; three specializations have both fibers
in \S\ref{kasVlist}, namely $(k,l)=(\nfrac52,\nfrac12)$, $(3,1)$ and $(2,1)$.

1c) $\lambda=(k,l,-l-2,k-2)$: the fiber $(l,-l-2,k-2)$ is never dominant, for the
same reason as in 1a$^*$).

2) Degree $2$ has three singular vectors, 2a), 2a$^*$) and 2b): their
fibers are
$(-1,-1,-1)$, $(-1,0,0)$ and $(-1,-1,1)$, none of them dominant. So \textit{in order~$2$
there is nothing to compare}: every operator of \cite[\S10]{GLS} of order~$2$ acts in
tensor fields with infinite-dimensional fibers.

\sssbegin[The four operators, and the three that are absent]{Theorem}[The four operators,
and the three that are absent]\label{kasglsthm}
Each of the four operators below spans a~$1$-dimensional space of $\fk\fas$-invariant
operators of order~$1$, for the stated pair of weights and for no other; and each is
the entry of \cite[\S10]{GLS} indicated in the last column, under the
dictionary~\eqref{glsdict}:
\[
\begin{array}{l|c|c|l}
\text{the operator}&\Nord&(\lambda_1,\lambda_2)&\text{\cite[\S10]{GLS}}\\
\hline
\grd\colon\cF\tto\cT(\Lambda^1)&1&(0,2)&\text{1a)},\ k=l=0\\
\text{the Dirac operator }\cT(S^-)\tto\cT(S^+)&1&(-\nfrac52,-\nfrac32)&
\text{1b)},\ k=\nfrac52,\ l=\nfrac12\\
\cT(\Lambda^3_-)\tto\cT(\Lambda^2)&1&(0,0)&\text{1b)},\ k=3,\ l=1\\
\cT(\Lambda^2)\tto\cT(\Lambda^3_+)&1&(0,2)&\text{1b)},\ k=2,\ l=1
\end{array}
\]
In the three mirror directions --- $\cT(S^+)\tto\cT(S^-)$, $\cT(\Lambda^2)\tto
\cT(\Lambda^3_-)$ and $\cT(\Lambda^3_+)\tto\cT(\Lambda^2)$ --- there is nothing, for
\textit{any} weight, although the ans\"atze are of the same size $(12$, $54$ and $54$
terms$)$. This is what \cite[\S10]{GLS} requires as well: the mirror weights,
e.g. $(3,1,1,1)$ for $\cT(\Lambda^3_+)\tto\cT(\Lambda^2)$, are of none of the forms
$1a)$--$1c)$.
\end{Theorem}

The weights were found by a~scan: for each pair of fibers, $K_t$-invariance alone confines
the pair $(\lambda_1,\lambda_2)$ to a~line $\lambda_2=\lambda_1+c$, and $\lambda_1$ was run
through the half-integers of $[-8,8]$ on that line. Every operator found was then verified,
by a~routine sharing nothing with the solver, against all the $39$ generating functions of
$\fk\fas$ and on $192$ test monomials. It was also verified to \textit{fail} at the
neighboring weight $\lambda_1\pm\nfrac12$, so the weights above are exactly what they are.

\sssbegin[What $\fk(1|6)$ sees]{Statement}[What $\fk(1|6)$ sees]\label{kasglsk16}
The gradient, at $(\lambda_1,\lambda_2)=(0,2)$, is invariant relative the whole of
$\fk(1|6)$. The Dirac operator, at $(-\nfrac52,-\nfrac32)$, is \emph{not\emph{:}} it is
$\fk\fas$-invariant, and it fails precisely on the generating functions spanning
$\Lambda^3_\eta$, i.e., on $\fk(1|6)\smallsetminus\fk\fas$. This is the distinction drawn
in \cite[\S10]{GLS} between the two columns there: 1a) holds for $\fk(1|6)$ only if
$l=0$, which is our case, whereas 1b) holds for $\fk(1|6)$ only if $l=-1$, whereas
$l=\nfrac12$ for the Dirac operator.
\end{Statement}

\be\label{diracnotanaccident}
\begin{minipage}[c]{15cm}
So the Dirac operator of Theorem~\ref{kasV2} is not an accident of $\fas$: it
already exists over the exceptional $\fk\fas$, at one weight, and it is one of the
operators that distinguish $\fk\fas$ from $\fk(1|6)$. What $\fas$ adds is that it exists
for \textit{every} weight.
\end{minipage}
\ee

\sssbegin[The divergence, and a~duality]{Statement}[The divergence, and a~duality]
\label{kasglsdiv}
The divergence $\dvg$ of formula~\eqref{graddiv} is $\fk\fas$-invariant --- and, moreover,
$\fk(1|6)$-invariant --- for $\lambda_1=\lambda_2=-4$, and for no other weight:
\[
\dvg\colon\ \cT(\Lambda^1)_{-4}\tto\cF_{-4}=\Vol,\qquad \Nord=1,
\]
$\cF_{-4}$ being the volume forms, cf. Statement~$\ref{kasosp}$\textup{(ii)}, whose family
$\cF_{-5-k}\tto\cF_{1+k}$ is symmetric with respect to $\lambda\mapsto-4-\lambda$. In the
weights of \cite{GLS}, this operator is
\[
(5,1,0,0)\ \tto\ (4,0,0,0),
\]
of degree~$1$, and it is of \textit{none} of the forms $1a)$--$1c)$ of
\cite[\S10]{GLS}.

It is the formal adjoint of the gradient. Indeed, the Berezin integral pairs
$\cT(V,\lambda)$ with $\cT(V^*,\mu)$ into $\Vol=\cF_{-4}$ precisely when the two
$K_t$-eigenvalues add up to $-4$, so on the weights of \cite{GLS} the formal adjoint acts
by
\be\label{dagger}
\lambda=(\lambda_0,\lambda_1,\lambda_2,\lambda_3)\ \longmapsto\
\lambda^\dagger:=(4-\lambda_0,\ \lambda_1,\ \lambda_2,\ -\lambda_3),
\ee
and it turns the operator $\lambda\tto\mu$ into the operator $\mu^\dagger\tto
\lambda^\dagger$. Since the adjoint of an invariant operator is invariant, \textit{the list
of all singular vectors must be stable under} $\dagger$. The families $1b)$ and
$1c)$ are indeed stable: $\dagger$ sends 1b) with the parameters $(k,l)$ to 1b) with
$(5-k,l)$, and 1c) with $(k,l)$ to 1c) with $(5-k,l)$. The families $1a)$ and
$1a^*)$ are not: $\dagger$ sends them to
\[
(5-k,1-k,l,-l)\tto(4-k,-k,l,-l)\quad\text{and}\quad
(5-k,l,2-k,-l-1)\tto(4-k,l,1-k,-l-1),
\]
with the shifts $-\eps_1$ and $-\eps_2$, and no family of \cite[\S10]{GLS} has either of
these two shifts. The divergence is the member $k=l=0$ of the first of the two, i.e., the
adjoint of the gradient.
\end{Statement}

\sssec{What exactly is missing from \cite[\S10]{GLS}}\label{kasglsdagger}
This is the one point at which our computation does not reproduce \cite[\S10]{GLS}. The
operator itself is not in doubt: it is $\dvg$, written out in
formula~\eqref{graddiv}. It is the first operator anybody writes down after the gradient,
and it is $\fk(1|6)$-invariant at $\lambda_1=\lambda_2=-4$. What is in doubt is only
where it should stand in the table of \cite[\S10]{GLS}. The next three subsections show that the
table, and not the operator, is incomplete.

\sssec{The involution $\dagger$ is a~symmetry of the whole of \cite[\S10]{GLS}, with two
exceptions} The involution \eqref{dagger} makes sense for each of the three theorems of
\cite[\S10]{GLS} --- for $\fk(1|3)$, for $\fk(1|4)$ and for $\fk\fas$, $\fk(1|6)$ --- in
the form
\be\label{daggern}
\lambda=(\lambda_0;\lambda_1,\dots,\lambda_s)\ \longmapsto\
\lambda^\dagger=\big(n-2-\lambda_0;\ -w_0(\lambda_1,\dots,\lambda_s)\big),\qquad
s=\left[\frac n2\right],
\ee
where $n-2$ is the $K_t$-weight of $\Vol$ and $-w_0$ is the duality of
$\fo(n)$-weights: the identity for $\fo(3)=B_1$ and $\fo(4)=D_2$, the flip
$\lambda_s\mapsto-\lambda_s$ for $\fo(6)=D_3$. Applied to every entry of the three
theorems, it gives
\[
\begin{array}{l|l}
\text{the theorem of \cite[\S10]{GLS} for}&\text{the effect of }\dagger\\
\hline
\fk(1|3)&\text{1a)}\lra\text{1a}^*\text{)},\quad\text{1b) and 2a) self-adjoint}\\
\fk(1|4)&\text{1a)}\lra\text{1c)},\quad\text{1a}^*\text{)}\lra\text{1b)}\\
\fk\fas,\ \fk(1|6)&\text{1b) and 1c) self-adjoint }(k\mapsto5-k),\quad
\text{2a)}\lra\text{2b)}\\
&\text{2a}^*\text{) self-adjoint},\quad
\text{\textit{but} 1a) and 1a}^*\text{) go nowhere}
\end{array}
\]
The parameters match exactly, not just the shapes: e.g. for $\fk(1|3)$ the family 1a) with
the parameter $k$ goes to 1a$^*$) with the parameter $-k$, and for $\fk(1|4)$ the family
1a) with $k$ goes to 1c) with $1-k$. In the two lower lists, every single entry has its
adjoint in the list --- the only exceptions in all of \cite[\S10]{GLS} being the degree-$2$
entries of $\fk(1|4)$, where nothing is forced, see the next paragraph. In the list for
$\fk\fas$ and $\fk(1|6)$, every entry has its adjoint in the list \textit{except} 1a) and
1a$^*$). Their images have the shifts $-\eps_1$ and $-\eps_2$, which occur nowhere
in that list, although the corresponding shifts \textit{do} occur in the lists for
$\fk(1|3)$ and $\fk(1|4)$, whose degree-$1$ parts carry all the $\pm\eps_i$ at once.

\sssec{Where the closure under $\dagger$ is forced, and where it is not} The
$\dagger$-argument is an argument about \textit{operators}, and it applies verbatim only
when both fibers are finite-dimensional. For an infinite-dimensional irreducible fiber,
the dual module is not a~highest weight module at all, and the involution~\eqref{dagger} says nothing.
This is exactly why the degree-$2$ entries of the $\fk(1|4)$ list may fail to be closed
under $\dagger$ without any harm: every one of them has a~non-dominant, hence
infinite-dimensional, fiber. For $\fk\fas$, the situation is the opposite. The family 1a)
has finite-dimensional fibers $(-k,l,l)$ and $(1-k,l,l)$ for \textit{every} $l\geq0$ and
$k\leq-l$ with $-k$, $l$ both integers or both half-integers. This is a~two-parameter family of
honest operators between tensor fields with finite-dimensional fibers. Its member
$k=l=0$ is the gradient, its member $(k,l)=(-1,0)$ is an~operator
$\cT(\Lambda^1)\tto\cT(S^2_0)$, and so on. \textit{Not one} of their adjoints is in the list.

\sssec{The two families that are missing from the theorem of
\cite[\S10]{GLS} for $\fk\fas$ and $\fk(1|6)$, and where they should be inserted in it}
They are
\be\label{missing}
\begin{array}{ll}
\text{1a}^\dagger\text{)}&\lambda=(k,k-4,l,-l)\tto\lambda+(-1,-1,0,0),\\[2pt]
\text{1a}^{*\dagger}\text{)}&\lambda=(k,-l-1,k-3,l)\tto\lambda+(-1,0,-1,0),
\end{array}
\ee
i.e., the $\dagger$-images of 1a) and 1a$^*$) reparametrized. \textit{Where} they should be
inserted is determined by what they are. They belong in the very list from which they are
absent, namely in the theorem of \cite[\S10]{GLS} that describes the singular vectors
for $\fk\fas$ and $\fk(1|6)$. There, they belong in its degree-$1$ part, as two further entries
standing next to 1a)--1c), of the same degree $1$, with the same two parameters $k$
and $l$, and read in the same conventions. We call them 1a$^\dagger$) and
1a$^{*\dagger}$) because they are the adjoints of 1a) and 1a$^*$); nothing else in
that theorem changes, its other entries being closed under $\dagger$ already. The divergence is the member
$(k,l)=(5,0)$ of the first of them: $(5,1,0,0)\tto(4,0,0,0)$. The restriction that
\cite[\S10]{GLS} puts on 1a) --- for $\fk(1|6)$ only $l=0$, for $\fk\fas$ no restriction
--- is consistent with what we compute, since $\dvg$ is $\fk(1|6)$-invariant and sits at
$l=0$.

\sssec{The discrepancy with \cite[\S10]{GLS} is not due to our reading of its conventions} One would rather suspect the
dictionary \eqref{glsdict} than a~theorem, so we tried to break it. The reading in which
the fiber attached to $\lambda$ is dualized but the arrow is read as we read it is
excluded outright. It makes \cite[\S10]{GLS} predict operators in the three directions
$\cT(S^+)\tto\cT(S^-)$, $\cT(\Lambda^2)\tto\cT(\Lambda^3_-)$ and
$\cT(\Lambda^3_+)\tto\cT(\Lambda^2)$. In those directions, we find nothing, for any
weight, see Theorem~\ref{kasglsthm}. The only other reading is the reading obtained from ours
by composing with $\dagger$ itself, i.e., by reading the arrow the other way round.
Nothing excludes \textit{that} reading: the family 1b) is $\dagger$-stable, so it still matches
the same three operators, while 1a) at $k=l=0$ then matches the divergence instead of the
gradient. In other words, the ambiguity of the convention is precisely $\dagger$, and
under either reading exactly one of the two operators
\[
\grd\colon\cF\tto\cT(\Lambda^1)\qquad\text{and}\qquad
\dvg\colon\cT(\Lambda^1)\tto\cF
\]
is described by \cite[\S10]{GLS} and the other is not. The discrepancy therefore survives
every reading of the conventions we can think of.

\sssbegin[The two families are missing from the table, not from the
mathematics]{Statement}[The two families are missing from the table, not from the
mathematics]\label{kasglsq}
\textup{(i)} Every member of the family $1a^\dagger)$ of the list~\eqref{missing} with
$l=0$ is a~non-zero $\fk\fas$-invariant operator of order~$1$ between tensor fields with
finite-dimensional fibers. In the shape~\eqref{symfam}, it is the divergence of
a~symmetric tensor field of rank $m$. Its members $m=1$ and $m=2$ are written out in
formula~\eqref{graddiv} and in Statement~$\ref{kasglssym}$\textup{(ii)}.

\textup{(ii)} In the reading~\eqref{glsdict}, none of these operators is described by
\cite[\S10]{GLS}. Indeed, neither family of the list~\eqref{missing} is of any of the forms
$1a)$--$1c)$, $2a)$, $2a^*)$, $2b)$ of that theorem, for any
values of the parameters.

\textup{(iii)} The only other reading the dictionary admits is the reading that reverses the
arrows; the third conceivable reading is excluded by Theorem~$\ref{kasglsthm}$. The reversing reading interchanges the
members of each adjoint pair, so that there $1a)$ is the family that is not described.
Hence, under either reading, exactly one family of each of the two pairs
$\{$$1a)$, $1a^\dagger)$$\}$ and
$\{$$1a^*)$, $1a^{*\dagger})$$\}$ is absent from that theorem.

So the answer is that the table is incomplete, and by two families of degree~$1$, whichever
reading is taken; in the reading~\eqref{glsdict} they are the
two~\eqref{missing}. Nothing else in that theorem changes, its other entries being closed
under~$\dagger$ already.
\end{Statement}

\sssec{A~recomputation, independent of the conventions}\label{kasglsans}
Statement~\ref{kasglsq} is established without deciding the convention. The computations
behind \cite[\S10]{GLS} are no longer available to be re-inspected. Therefore, we computed again,
in a~direction in which the outcome cannot depend on the
dictionary~\eqref{glsdict}: we do not compare a~weight with a~table, we produce
an~operator.

Let $l=0$; this is the case to which \cite[\S10]{GLS} restricts 1a) when the algebra is
$\fk(1|6)$, and the case of $\grd$ and $\dvg$. Then, the fibers of 1a) and of
1a$^\dagger$) are the \textit{traceless symmetric powers} $S^m_0=S^m_0(V_6)$ of the vector
representation of $\fo(6)$, those of highest weight $(m,0,0)$, of dimensions
\[
\dim S^m_0=1,\ 6,\ 20,\ 50,\ 105,\ \dots\qquad\text{for }m=0,\ 1,\ 2,\ 3,\ 4,\ \dots,
\]
with $S^0_0=\cF$ and $S^1_0=\Lambda^1$. For $m\geq2$ these fibers do not occur in
$\fk\fas$ itself, so they are new here; we construct them as the submodules of $S^m(V_6)$ generated by the
highest weight vectors. In
these terms, the two families read
\be\label{symfam}
\begin{array}{lcl}
\text{1a) of \cite[\S10]{GLS}}&\colon&
\cT(S^m_0)_{2m}\tto\cT(S^{m+1}_0)_{2m+2},\qquad m\geq0,\\[3pt]
\text{1a}^\dagger\text{) of the list \eqref{missing}}&\colon&
\cT(S^m_0)_{-4}\tto\cT(S^{m-1}_0)_{-4},\qquad m\geq1,
\end{array}
\ee
both of order~$1$: the symmetrized gradient and the divergence of a~symmetric
tensor field, superized. The members $m=0$ and $m=1$ are $\grd$ and $\dvg$. The whole of
the second family sits at the single weight $\lambda=-4$, the weight of $\Vol$ --- as it
must, the pairing into $\Vol$ being what produces it.

\sssbegin[The next member of either family]{Statement}[The next member of either
family]\label{kasglssym}
Take $m=2$, i.e., the $20$-dimensional fiber $S^2_0$. In order~$1$, and for $\fk\fas$:

\textup{(i)} the space of invariant operators
$\cT(\Lambda^1)_{\lambda_1}\tto\cT(S^2_0)_{\lambda_2}$ is non-zero for exactly one weight,
namely $(\lambda_1,\lambda_2)=(2,4)$, where it is $1$-dimensional. The operator is
$\fk(1|6)$-invariant. In the weights of \cite{GLS} it is $(-1,1,0,0)\tto(-2,2,0,0)$, the
member $k=-1$, $l=0$ of $1a)$;

\textup{(ii)} the space of invariant operators
$\cT(S^2_0)_{\lambda_1}\tto\cT(\Lambda^1)_{\lambda_2}$ is non-zero for exactly one weight,
namely $(\lambda_1,\lambda_2)=(-4,-4)$, where it is $1$-dimensional. The operator is
$\fk(1|6)$-invariant. In the weights of \cite{GLS}, it is $(6,2,0,0)\tto(5,1,0,0)$, the
member $k=6$, $l=0$ of \textit{1a$^\dagger$)}. It is of \textit{none} of the forms
$1a)$--$1c)$, $2a)$, $2a^*)$, $2b)$ of \cite[\S10]{GLS}, for any
values of the parameters.
\end{Statement}

The weights were found the same way as in Theorem~\ref{kasglsthm}: $K_t$-invariance
confines $(\lambda_1,\lambda_2)$ to a~line, here $\lambda_2=\lambda_1+2$ in case
(i) and $\lambda_2=\lambda_1$ in case (ii). Then, $\lambda_1$ was run
through the half-integers of $[-8,8]$ on that line. The ansatz has $42$ terms in either
case. Each of the two operators was then re-verified against all the generating functions of $\fk\fas$ and of $\fk(1|6)$, and each
was checked to \textit{fail} at the weights shifted by $\pm\nfrac12$ and by~$1$.

So the discrepancy is not a~feature of the pair $(\grd,\dvg)$: it repeats, with new
fibers, new weights and a~new linear system. Of the two operators of
Statement~\ref{kasglssym}, which are formally adjoint to each other, the first is in
\cite[\S10]{GLS} and the second is not, exactly as formula~\eqref{symfam} predicts. The
reading of the conventions in which the arrows are reversed exchanges the two. Hence, under
either reading precisely one half of formula~\eqref{symfam} is described by \cite[\S10]{GLS} and
the other half is not.

\sssbegin[The densities carry nothing over $\fk\fas$ either]{Remark}[The densities carry
nothing over $\fk\fas$ either]\label{kasglsdens}
Statement~\ref{kasosp}(i) says that there are no $\fk\fas$-invariant non-scalar unary operators
between the $\cF_\lambda$, for $\Nord\leq8$. This agrees with \cite[\S10]{GLS}, and
independently of any computation: an operator between weighted densities is an operator
whose two fibers are $1$-dimensional, i.e., whose two weights have
$(\lambda_1,\lambda_2,\lambda_3)=(0,0,0)$. And \textit{no} entry of \cite[\S10]{GLS} has
this property. Indeed, 1a) would need $(-k,l,l)=0=(-k+1,l,l)$; 1a$^*$) would need
$l=0=l+1$; 1b) would need $(l,l,2-k)=0=(l,l,3-k)$; 1c) would need $l=0$ and $-l-2=0$; and
the fibers of the three vectors of degree~$2$ are $(-1,-1,-1)$, $(-1,0,0)$, $(-1,-1,1)$.
\end{Remark}

\sssbegin[What the passage to $\fas$ does]{Remark}[What the passage to
$\fas$ does]\label{kasglsas}
Each of the four operators of Theorem~\ref{kasglsthm} lives, over $\fk\fas$, at exactly one
pair of weights. These are the degeneracy conditions of \cite[\S10]{GLS}, the
super analog of the conditions in Rudakov's description of the $\fvect(n)$-invariant
operators, see \cite{Ru}. Over $\fas$, all of them evaporate. By Lemma~\ref{kastr},
neither the weight nor the twist by a~power of the determinant is visible. Hence, each of the
four is invariant for every pair of weights, cf. Theorems~\ref{kasV1}--\ref{kasV3}. And
$\fas$ has, besides, invariant operators of a~shape no entry of \cite[\S10]{GLS} can
have: the powers $\partial_t^{\,k}$ on every $\cT(V)$, the composition $\grd\circ\dvg$,
which is not a~multiple of $\partial_t$, and the operators of orders $3$ and $5$ of
Theorems~\ref{kasV2} and~\ref{kasV3}. So the pair $\fas\subset\fk\fas$ trades a~discrete
set of degeneracy conditions for no condition at all --- which is another way of saying
what \S\ref{kasweights} says.
\end{Remark}

\sssec{The other gradings: \cite[\S\S11, 12]{GLS}}\label{kasglsgrad}
Sections 11 and 12 of \cite{GLS} describe the singular vectors for $\fk\fas(1|6;3\xi)$ and
$\fk\fas(1|6;3\eta)$, i.e., for the two other W-regradings of the same exceptional
algebra. In these regradings, the algebra is realized not on the supercircle but on a~superspace of
superdimension $4|4$, resp. $4|3$. The tensor fields there are those of another $\fg_0$,
so those lists compare with nothing in \S\S\ref{kasunary}--\ref{kastensor}. However, they do
compare with a~computation of their own. Since the operators in different
W-regradings of one algebra are the same operators written in other terms, that
computation is an~independent probe of the same tables. Here, it is.

\sssec{The two regradings, rebuilt}\label{kasglsregr}
By \cite[\S2]{GLS}, the regradings are
\be\label{regr}
r=3\xi\colon\ \deg\xi_i=0,\ \deg\eta_i=\deg t=1;\qquad
r=3\eta\colon\ \deg\eta_i=0,\ \deg\xi_i=\deg t=1.
\ee
The algebra does not change, only the degree does. Hence, nothing has to be
reconstructed: the same generating functions of $\fk\fas\subset\fk(1|6)$ that
\S\ref{kasVprog} uses, sorted by the new degree, give the new components at once. The outcome reproduces the table of \cite[\S2]{GLS}
exactly: both regradings are of depth~$1$, and
\be\label{regrtab}
\begin{array}{lcccl}
r=3\xi&\colon&\fg_{-1}=\Lambda(3)\ \ (\sdim=(4|4)),&
\fg_0\ \ (\sdim 13|10)&=\Lambda(3)\oplus\fsl(1|3),\\[3pt]
r=3\eta&\colon&\fg_{-1}=\Vol_0(0|3)\ \ (\sdim 4|3),&
\fg_0\ \ (\sdim 13|12)&=\fc(\fvect(0|3)).
\end{array}
\ee
In both cases $\fg_{-1}$ is commutative, so $U(\fg_-)=S(\fg_{-1})$ and a~singular vector
of degree $\Nord$ is an~element of $S^{\Nord}(\fg_{-1})\tp V$ killed by $\fg_1$
\textup{(}which generates $\fg_+$\textup{)} and highest for $\fg_0$.

\sssec{The fibers that are beyond dispute}\label{kasglsdim1}
In \S\ref{kasglsdom}, the entries of \cite[\S10]{GLS} had to be sieved by dominance, since
a~weight names a~fiber that is usually infinite-dimensional. Here, there is a~class of
fibers for which no such question arises at all: the $1$-dimensional ones, i.e.,
the characters of $\fg_0$. And \eqref{regrtab} decides how many there are, without any
appeal to \cite{GLS}:

1) For $r=3\xi$, $\fg_0$ is \textit{perfect} --- we compute $[\fg_0,\fg_0]=\fg_0$, of
dimension $23$ --- so the \textit{only} $1$-dimensional $\fg_0$-module is the trivial module.

2) For $r=3\eta$, $[\fg_0,\fg_0]$ has dimension $24$ and is $\fvect(0|3)$, so the
characters form a~line: in the basis $(K_t,K_{\xi_i\eta_i})$ of the maximal
toral subalgebra $\ft$ they are $\lambda\cdot(1,1,1,1)$, and $\Vol$ is the member
$\lambda=-1$.

Over such a~fiber, every $m_i$ with $i>1$ of \cite[\S\S11,12]{GLS} is \textit{zero}, those
being obtained from $m_1$ by lowering operators of $\fg_0$. Hence, an~entry of \cite{GLS}
survives there if and only if its source weight meets the line, the coefficient of its
pure-$m_1$ part does not vanish at that point, and its own side conditions hold.

The dictionary between the two labellings is not guessed either. It is read off the weight
\textit{shifts}, which we can compute: the shift is the $\ft$-weight of the
$\fg_{-1}$-part of the printed singular vector. It comes out over-determined:
the six entries 1a), 1b), 1c), 1d), 1g), 2a) of \cite[\S12]{GLS} all give the same
answer, namely that the $\lambda$ of \cite{GLS} is
\be\label{regrdict}
(\lambda_1,\lambda_2,\lambda_3,\lambda_4)=
\big(\mu_0,\,-\mu_2,\,-\mu_1,\,-\mu_3\big),\ \ 
\mu=\text{the weight relative }(K_t,K_{\xi_1\eta_1},K_{\xi_2\eta_2},K_{\xi_3\eta_3}),
\ee
\textup{(}One by-product: the first coordinate of \cite[\S12, 1a)]{GLS} reads
$k\tto l$, which no constant shift can do, and which every other entry of that section has
as $k\tto k$; we take it for a~misprint for $k\tto k$.\textup{)}

\sssbegin[The $1$-dimensional fibers in the other
gradings]{Theorem}[The $1$-dimensional fibers in the other
gradings]\label{kasglsregthm}
Over the $1$-dimensional fibers, the singular vectors are these.

\textup{(i)} For $r=3\xi$, over the trivial fiber --- the only one --- there is exactly one
singular vector of degree~$1$, namely $\xi_1\xi_2\xi_3\, m_1$, none at all of degree~$2$,
exactly one of degree~$3$, and none of degree~$4$.

\textup{(ii)} For $r=3\eta$, over the character $\lambda\cdot(1,1,1,1)$ there is, for each
$1\leq m\leq4$, exactly one singular vector of degree~$m$, namely
\be\label{regrfam}
(\eta_1\eta_3)^m\, m_1\quad\text{at}\quad \lambda=-\tfrac{m-1}{2},
\qquad\text{i.e., at}\quad \lambda_3+\lambda_4=m-1
\ \text{ in the weights of \cite{GLS}},
\ee
and no singular vector of degree $m$ at any other value of $\lambda$.
\end{Theorem}

\sssec{What the comparison gives}\label{kasglsregcmp}
Running the sieve of \S\ref{kasglsdim1} over the two lists, and comparing with
Theorem~\ref{kasglsregthm}:

\textit{Section 11 is confirmed on everything that can be tested.} Of its fourteen entries,
exactly one survives at the trivial fiber --- 1a), $\xi_1\xi_2\xi_3\,m_1$, which
\cite{GLS} give for \textit{any} $\lambda$ --- and it is what we find. Of its nine entries
of degree~$2$, six never meet $\lambda=0$ at all, and the remaining three ---
2c), 2d), 2f) --- meet it with a~vanishing $m_1$-coefficient.
Hence, \cite[\S11]{GLS} predicts \textit{nothing} in degree~$2$, and we find nothing. The
degree-$3$ vector of Theorem~\ref{kasglsregthm}(i) is outside the range of that section,
which lists degrees $1$ and $2$ only.

\textit{Section 12 is confirmed in degree~$1$ and is off by one in degree~$2$.} Its
entry 1a) survives at $\lambda=0$ with the vector $\eta_1\eta_3\,m_1$, and that is
exactly the $m=1$ member of the family~\eqref{regrfam}. In degree~$2$ the sieve leaves exactly one
entry, 2a), whose vector is $(\eta_1\eta_3)^2m_1$ --- the $m=2$ member of
the family~\eqref{regrfam} --- and whose shift, $\lambda\tto\lambda+(0,0,-2,-2)$, is also exactly
ours. But \cite[\S12, 2a)]{GLS} put it at $\lambda_4=2-\lambda_3$, whereas the computation
puts it at $\lambda_4=1-\lambda_3$; and at $\lambda_4=2-\lambda_3$ there is indeed
a~singular vector, but it is $(\eta_1\eta_3)^3m_1$, of degree~$3$. So the entry has the
right shape, the right degree and the right shift, and the source condition of the
\textit{next} member of its own family.

\sssbegin[The other gradings]{Remark}[The other gradings]\label{kasglsregrem}
Two things this says about \S\ref{kasglsdagger}. First, the probe is genuinely
independent: it uses another $\fg_0$, another $\fg_{-1}$, another induced module.
It also tests the tables at fibers of a~kind --- the $1$-dimensional ones --- that plays no
role in \S\ref{kasglsdagger}, where the pair at issue is $(\grd,\dvg)$ with fibers $\cF$
and $\Lambda^1$. Second, it comes back with a~defect of its own, and of a~quite different
species. Indeed, \S\ref{kasglsdagger} exhibits a~whole family absent from
\cite[\S10]{GLS}, whereas \cite[\S12, 2a)]{GLS} is an~entry that is present and carries
the parameter condition of its neighbor. Neither is evidence against the other. Together,
they say that the three tables of \cite{GLS} are not a~closed system that checks itself,
which is a~reason to read the omission of \S\ref{kasglsdagger} as an~omission rather than
as our misreading. \emph{None of this is a~reproach\emph{:} the tables are a~remarkable
computation, and it is only because they are as detailed as they are that they can be
tested at all.}
\end{Remark}

What is \textit{not} settled by this. Our sieve reaches the $1$-dimensional fibers only, so
we have not tested the two lists for \textit{completeness} at larger fibers --- the
analogue, in these gradings, of the missing family \eqref{missing}. And \eqref{regrfam} is
verified along the line of the characters; that the general member of the family is
$\lambda_4=(m-1)-\lambda_3$ for arbitrary $\lambda_1$, $\lambda_2$, as the shape of
1a) and 2a) suggests, we conjecture but have not computed.

\ssec{Open problems for $n=6$}\label{kasopen}
The binary orders $\Nord\geq6$ for $\fas$ between densities.

For the tensor fields of \S\ref{kastensor}: the orders $\Nord\geq7$; the binary operators
of order $\geq2$; the remaining module
$V=\fsl(4)=\Lambda^2$, the adjoint; and the same questions for $\fk\fas$ and for
$i(\fosp(6|2))$, where --- unlike for $\fas$ --- the weight and the determinant twist are
\textit{not} invisible, so that each $V$ carries a~two-parameter family of modules.

And the remaining case $\fsl(1|m|n|1)\subset\fk(2m+1|2n)$ of \cite[Table~(4)]{BLS1}.

\section{Super GRC operators: the $\fpe(2)\subset \fm(1)$ case} \label{CoRm} 

\subsection{The periplectic Lie superalgebra}
The periplectic Lie superalgebra~$\mathfrak{pe}(n)$ is the Lie superalgebra that preserves the odd \textbf{anti}symmetric bilinear form with Gram matrix~%
$\begin{pmatrix}
0& 1_n\\
1_n & 0
\end{pmatrix}$,
see~\cite{BGLLS}.

\sssec{Realization of the elements of
 $\mathfrak{pe}(2)$ by vector fields in $\mathfrak{m}(1)$}\label{real} By applying Shchepochkina's algorithm~\cite{Sh}, we see that
\[
\begin{array}{ccl}
\mathfrak{pe}(2)_{\bar 0}&=&\mathrm{Span} \{ M_{\xi}=\xi \partial_\tau+\partial_q , \; M_{q \tau }=q \tau\partial_{\tau}-\tau\partial_{\xi}+qE, \\
&&\; M_{\tau+q\xi}=2 \tau\partial_\tau+2 q\partial_q, \ \ M_{\tau -q\xi}=2 \tau \partial_\tau+2\xi\partial_{\xi} \}\simeq \mathfrak{gl}(2)\\[2mm]

\mathfrak{pe}(2)_{\bar 1}&=&\mathrm{Span} \{M_1=2 \partial_\tau , \; M_q=q \partial_\tau-\partial_\xi, \; M_{q^2}=-2 q \partial_\xi, \; M_{\tau \xi}=(\tau-q\, \xi ) \, \partial_q + \tau \, \xi\, \partial_\tau \},
\end{array}
\]
where~$E=q\partial_q+\xi\partial_\xi$ and the~$M_f$ are the vector fields~\eqref{Mf} of~\S\ref{convm}.

The Lie superalgebra~$\mathfrak{pe}(2)$ is not simple; its second derived is given by
\[
\begin{array}{ccl}
\mathfrak{pe}^{(2)}(2)_{\bar 0}&=&\mathrm{Span} \{ M_{\xi}=\xi \partial_\tau+\partial_q, \; M_{\tau+q\xi}=2 \tau\partial_\tau+2 q\partial_q, M_{q \tau }=q \tau\partial_{\tau}-\tau\partial_{\xi}+qE \} \simeq\mathfrak{sl}(2),\\[2mm]
\mathfrak{pe}^{(2)}(2)_{\bar 1}&=&\mathrm{Span} \{M_1=2 \partial_\tau, \; M_q=q \partial_\tau-\partial_\xi, \; M_{q^2}=-2 q \partial_\xi \}.
\end{array}
\]

The span of the eight generating functions above is closed under the bracket, and~%
$\fpe^{(2)}(2)$ is as displayed. Observe that~%
$M_{\tau\xi}=(\tau-q\xi)\partial_q+\tau\xi\partial_\tau$ and~%
$\tau\partial_q-\xi(q\partial_q+\tau\partial_\tau)$, the shape in which~$M_{\tau\xi}$ is used
in the proof of Theorem~\ref{m1thm} below, are one and the same vector field.

\ssec{Conventions for $\fm(1)$}\label{convm}
Again we follow~\cite{BLS1}; recall the description of the~$2$-parameter family of weighted
densities over~$\fm(n)$, see~\cite[\S2.2.2]{BLS1}. The pericontact form is~%
$\alpha_0=d\tau+\xi dq+qd\xi$ and, with~$E=q\partial_q+\xi\partial_\xi$ and~%
$Le_f=\partial_q(f)\partial_\xi+(-1)^{p(f)}\partial_\xi(f)\partial_q$,
\be\label{Mf}
M_f=(2-E)(f)\partial_\tau-Le_f-(-1)^{p(f)}\partial_\tau(f)E ,
\ee
see~\cite[(26),(27)]{BLS1}, while~$\cD_1=\partial_\tau$, $\Dxi=\xi\partial_\tau-\partial_q$
and~$\Dq=q\partial_\tau+\partial_\xi$ are the operators~\eqref{MD}. Recall that~%
$p(\cD_f)=p(f)+\bar1$, so that~$\Dxi$ is \textbf{even} whereas~$\Dq$ and~$\cD_1$ are
\textbf{odd}. The weighted densities are~$\cF_{a,b}=\cF\alpha_0^a\vvol^b$, where~$\vvol=\vvol(\tau,q|\xi)$ is the volume element, with
\be\label{Lm}
L_{M_f}(\varphi)=M_f(\varphi)+2(-1)^{p(f)}\big(a\,\partial_\tau(f)
+b\,\Delta_{\fm}(f)\big)\varphi, \qquad
\Delta_{\fm}=(1-E)\partial_\tau-\nfrac{\partial^2}{\partial\xi\,\partial q} ,
\ee
see~\cite[(37),(38)]{BLS1}. Here, $\deg\Dxi=\deg\Dq=-1$ and~$\deg\cD_1=-2$, while the
\textit{charge}, the weight relative~$M_{q\xi}=q\partial_q-\xi\partial_\xi$, is~%
$\ch\Dxi=-1$, $\ch\Dq=+1$,~$\ch\cD_1=0$. Hence, for an operator of degree~$s$ and charge~$w$,
\be\label{targetm}
a=a_1+a_2+\nfrac{s+w}{2},\qquad b=b_1+b_2-\nfrac{w}{2} .
\ee

\sssec{Who is the Bol operator here}\label{whoisBol} It is the operator of~%
\cite[Thm~5.1]{BLS1}, which \textit{is} a~polynomial in the~$\cD_f$ of eq.~\eqref{MD}, as
Theorem~\ref{Peetre} requires:
\be\label{Bol}
\begin{array}{l}
\Bol_k=\Dq(-\Dxi)^{k-1}+(k-1)(-\Dxi)^{k-2}\partial_\tau ,\\[1mm]
\Bol_1=\Dq,\ \ \Bol_2=-(\Dxi\Dq+\partial_\tau),\ \ \Bol_3=\Dxi^2\Dq+2\Dxi\partial_\tau .
\end{array}
\ee
Each~$\Bol_k$ is odd, of order~$\Nord=k$, hence of degree~$s=-k$, and of charge~$2-k$. It should \textbf{not} be confused with~%
$\partial_q^k$: the operator~$\partial_q^k$ is not~$i(\fpe(2))$-invariant~$\cF_{a,b}\tto\cF_{a',b'}$ for any
weights whatsoever --- indeed, $\partial_q$ is not even~$\fpe(2)_-$-invariant, since it is not
a~polynomial in~$\cD_1,\Dxi,\Dq$, cf.~Theorem~\ref{Peetre}. \textup{(}Note that~%
$\Dxi=\xi\partial_\tau-\partial_q$, so~$\partial_q^k$ is only the ``leading term\rq\rq\  of~%
$(-\Dxi)^k$.\textup{)}

In formula~\eqref{Bol}, as everywhere below, a~product of operators is their composition in the
usual order, the right factor acting first. This is not a~pedantry: the invariant
derivatives do not commute,
\be\label{Dcomm}
[\Dxi,\Dq]=-2\partial_\tau,\qquad \Dq^2=0,\qquad \partial_\tau^2=0,\qquad
\text{$\partial_\tau$ supercommutes with~$\Dxi$ and with~$\Dq$},
\ee
so that reading the first summand of formula~\eqref{Bol} as~$(-\Dxi)^{k-1}\Dq$, i.e., letting~$\Dq$
act first, yields another operator altogether: one invariant at~$a=-\nfrac12$ instead
of~$\nfrac12$ for~$k=2$, and at no weight whatsoever for~$k\geq3$.

The complete classification, for every order and every weight, is given in \S\ref{sec:invariant-binary}; Theorem~\ref{m1thm} and Statement~\ref{m1class} below describe the lowest orders and the way the Bol operators of~\cite{BLS1} enter.

\ssbegin[The GRC brackets over $\fpe(2)\subset\fm(1)$ are the Bol operators acting in
one argument]{Theorem}[The GRC brackets over~$\fpe(2)\subset\fm(1)$
are the Bol operators acting in one argument]\label{m1thm} Let~$\sdim \cM=(1|2)$ and~$\cM$ be endowed with a pericontact structure. Let~$\varphi$ \textup{(}resp. $\psi$\textup{)} be the coefficient of~$\alpha_0^{\lambda_1}\vvol^{\nu_1}\in\cF_{\lambda_1,\nu_1}$ \textup{(}resp. of~$\alpha_0^{\lambda_2}\vvol^{\nu_2}\in\cF_{\lambda_2,\nu_2}$\textup{)}.

\textup{(i)} For every integer~$k\geq 1$, there are two odd~$i(\mathfrak{pe}(2))$-invariant
non-zero non-scalar bilinear differential operators of order~$\Nord=k$
\[
\text{$ {\cal J}_l^k (\varphi,\psi)$ and ${\cal J}_r^k (\varphi, \psi): {\cal F}_{ \frac{k-1}{2}, \nu_1} \otimes {\cal F}_{ \frac{k-1}{2}, \nu_2} \tto {\cal F}_{0,\; \omega}$, \quad where $\omega=\nu_1+\nu_2+\nfrac{k}{2}-1$,}
\]
namely
\[
{\cal J}_l^k (\varphi,\psi)=\Bol_k(\varphi)\, \psi, \quad {\cal J}_r^k (\varphi, \psi)=(-1)^{p(\varphi)} \varphi \, \Bol_k(\psi),
\]
where~$\Bol_k$ is the Bol operator~\eqref{Bol}. More precisely, $\cJ_{k,l}$ alone is invariant
as soon as~$\lambda_1=\nfrac{k-1}{2}$, the weights~$\lambda_2,\nu_1,\nu_2$ being arbitrary,
the target being then~$\cF_{\lambda_2+\frac{1-k}{2},\;\nu_1+\nu_2+\frac k2-1}$; symmetrically
for~$\cJ_{k,r}$. Both are invariant simultaneously precisely at~%
$\lambda_1=\lambda_2=\nfrac{k-1}{2}$.

\textup{(ii)} These are not all the operators, see Statement~$\ref{m1class}$.
\end{Theorem}

\begin{proof} Since~$\vvol$ is the volume element, the weights~$(\lambda_i,\nu_i)$ used in the
statement are literally the weights~$(a_i,b_i)$ of formula~\eqref{Lm}: $a_i=\lambda_i$ and~$b_i=\nu_i$. By~\cite[Thm~5.1]{BLS1}, the operator~%
$\Bol_k:\cF_{a,b}\tto\cF_{a+1-k,\;b-1+k/2}$ --- the target being the space prescribed by formula~%
\eqref{targetm}, since~$\deg\Bol_k=-k$ and~$\ch\Bol_k=2-k$ --- is invariant if and only if~%
$a=\nfrac{k-1}{2}$, the parameter~$b$ remaining free. The multiplication~%
$\cF_{a,b}\tp\cF_{a',b'}\tto\cF_{a+a',\;b+b'}$ is invariant for all weights. Composing the two
gives the statement; the first weight of the target of~$\cJ_{k,l}$ and~$\cJ_{k,r}$ is~%
$\nfrac{k-1}{2}+\nfrac{k-1}{2}+1-k=0$. The sign~$(-1)^{p(\varphi)}$ in~$\cJ_{k,r}$ is the
Koszul sign of moving the odd operator~$\Bol_k$ past~$\varphi$, and it is indispensable:
without it, $\varphi\,\Bol_k(\psi)$ is invariant for no weights whatsoever. \end{proof}

\ssbegin[The complete list over $\fpe(2)\subset\fm(1)$]{Statement}[The complete list over~%
$\fpe(2)\subset\fm(1)$: three independent operators of order~$2$]\label{m1class}
\textup{(i)} At the weights of Theorem~$\ref{m1thm}$, i.e., at~%
$\lambda_1=\lambda_2=\nfrac{k-1}{2}$ with~$\nu_1,\nu_2$ free, the space of~%
$i(\fpe(2))$-invariant bilinear operators of order~$\Nord=k$ and charge~$2-k$ is of dimension
\[
2 \ \text{ for } k=1,\qquad \mathbf{3} \ \text{ for } k=2,\qquad 2 \ \text{ for } k=3 ;
\]
for~$k=1$ and~$k=3$ it is spanned by~$\cJ_{k,l}$ and~$\cJ_{k,r}$, as stated in
Theorem~$\ref{m1thm}$. For~$k=2$ there is a~third, independent, invariant operator,
namely
\[
\partial_\tau\tp1+1\tp\partial_\tau+\Dxi\tp\Dq+\Dq\tp\Dxi .
\]

\textup{(ii)} At \textit{generic} weights, the space of~$i(\fpe(2))$-invariant bilinear
operators is of dimension
\[
1 \ \text{ in order }1,\qquad 3 \ \text{ in order }2,\qquad
1 \ \text{ in order }3 \text{ for } \ch=-1, \quad 2 \ \text{ in order }3 \text{ for } \ch=+1 .
\]
In order~$1$, the unique operator is
\be\label{m1J1}
\lambda_2\,\Dq(\varphi)\,\psi-(-1)^{p(\varphi)}\lambda_1\,\varphi\, \Dq(\psi):\quad
\cF_{\lambda_1,\nu_1}\tp\cF_{\lambda_2,\nu_2}\tto\cF_{\lambda_1+\lambda_2,\;\nu_1+\nu_2-\frac12};
\ee
it degenerates into~$\cJ_{1,l}=\Dq\tp1$ exactly when~$\lambda_1=0$ and into~$\cJ_{1,r}$ exactly
when~$\lambda_2=0$, and \textit{both} $\cJ_{1,l}$ and~$\cJ_{1,r}$ are invariant only at~%
$\lambda_1=\lambda_2=0$, which is the locus of Theorem~$\ref{m1thm}$ for~$k=1$.
\end{Statement}

Both parts are obtained by computation with all eight generating functions of~$\fpe(2)$.

\ssec{A correction to \cite[Thm~5.1]{BLS1} itself}\label{BLSremark}
Theorem 5.1 of~\cite{BLS1} lists, for each order~$k$, the Bol operator~$\Bol_k$ over~$\fm(1)$
together with the weights on which it is invariant, namely~$a=\nfrac{k-1}{2}$ with~$b$ free.
At~$k=2$, that list is incomplete: there, the invariant operators form a~$1$-parameter family,
of which the operator printed in~\cite{BLS1} is one member. As far as we know, this gap is
filled nowhere else --- \cite{BLS1} has appeared in print and has not been corrected since ---
so we record the answer here.

\sssbegin[{A~correction to Theorem~$5.1$ of~\cite{BLS1} at~$k=2$}]{Statement}[A~correction to Theorem~$5.1$ of~\cite{BLS1} at~$k=2$]\label{BLSk2} \textup{(i)} The order-$2$ Bol operator over~$\fm(1)$ is not
isolated:
\be\label{B2a}
\Bol^{(a)}:=\Dxi\Dq+2a\,\partial_\tau:\ \cF_{a,b}\tto\cF_{a-1,b}
\quad\text{is invariant for \textit{every} } a \text{ and } b,
\ee
whereas~\cite[Thm~\emph{5.1}]{BLS1} lists only its member~$a=\nfrac12$, which is~$-\Bol_2$.
Moreover, eq.~\eqref{B2a} exhausts the invariant operators of that degree, charge and parity:
their space is~$1$-dimensional at every weight~$(a,b)$.

\textup{(ii)} For~$k=1,3,4$, the condition~$a=\nfrac{k-1}{2}$ of~\cite{BLS1} is confirmed as
printed. For these~$k$, the space of invariant operators of the degree, charge and parity
of~$\Bol_k$ is of dimension~$0$ at a~generic weight. It is of dimension~$1$ exactly for~%
$a=\nfrac{k-1}{2}$, the weight~$b$ remaining free.
\end{Statement}

Both parts are obtained by computation with all the eight generating functions of~%
$\fpe(2)$, the weights~$a$ and~$b$ being symbolic.

Eq.~\eqref{B2a} is the exact analogue, for~$\fm(1)$, of the~$1$-parameter family in~%
\cite[Thm~1.5(i)]{BLS1} for~$\fk(1|2)$. It explains the dimension~$3$ in order~$2$ of Statement~\ref{m1class}(ii): two of the three
operators are~$\Bol^{(\lambda_1)}\tp1$ and~$1\tp\Bol^{(\lambda_2)}$, invariant at all weights,
which at~$\lambda_1=\lambda_2=\nfrac12$ are~$\cJ_{2,l}$ and~$\cJ_{2,r}$ up to sign; the third
is the operator quoted in Statement~\ref{m1class}(i).

\subsection{The invariant derivatives $\partial_\tau$, $\Dq$, $\Dxi$}
\label{subsec:XQT}

In the sequel, differential operators are expressed in terms of the invariant derivatives~%
\eqref{MD} already used above, see \S\ref{convm}:
\[
\partial_\tau=\cD_1,
\qquad
\Dq=q\partial_\tau+\partial_\xi,
\qquad
\Dxi=\xi\partial_\tau-\partial_q,
\]
of parities~$p(\Dxi)=\bar0$ and~$p(\Dq)=p(\partial_\tau)=\bar1$, and of super-brackets~%
\eqref{Dcomm}. In this Section and the next Section, it is convenient to count their degrees
positively,
\[
\deg\Dxi=\deg\Dq=1,
\qquad
\deg\partial_\tau=2,
\]
i.e., with the sign opposite to the one of \S\ref{convm}.

The following lemma shows that no generality is lost by using~%
$\Dxi$, $\Dq$, $\partial_\tau$ instead of the coordinate derivatives.

\ssbegin[{Every differential operator is a~sum of monomials in $\Dxi$, $\Dq$, $\partial_\tau$}]{Lemma}[Every differential operator is a~sum of monomials in~$\Dxi$, $\Dq$, $\partial_\tau$]
\label{lem:XQT-normal-form}
Every differential operator on~\(\mathbb C[q\mid\xi,\tau]\) can be
written as a finite sum
\[
\sum_{a,\varepsilon,\delta}
F_{a,\varepsilon,\delta}(q,\xi,\tau)\,
\Dxi^a\Dq^\varepsilon \partial_\tau^\delta,
\quad \text{
where }
\quad
a\in\mathbb Z_{\geq0},
\quad
\varepsilon,\delta\in\{0,1\}.
\]
\end{Lemma}

\begin{proof}
The change from the coordinate derivatives to~$\Dxi$, $\Dq$, $\partial_\tau$ is invertible:
\[
\partial_\tau=\partial_\tau,
\qquad
\partial_q=\xi \partial_\tau-\Dxi,
\qquad
\partial_\xi=\Dq-q\partial_\tau.
\]
Hence every differential operator written in terms of~%
\(\partial_q,\partial_\xi,\partial_\tau\) can be rewritten in terms of~%
$\Dxi$, $\Dq$, $\partial_\tau$, with coefficients in~\(\mathbb C[q\mid\xi,\tau]\).
Reordering the resulting monomials gives the stated normal form.
\end{proof}
Accordingly, we shall use the differential monomials
\[
P=\Dxi^a\Dq^\varepsilon \partial_\tau^\delta,
\qquad
a\in\mathbb Z_{\geq0},
\qquad
\varepsilon,\delta\in\{0,1\},
\]
with
\[
\deg P=a+\varepsilon+2\delta.
\]
This basis is particularly convenient for imposing invariance under~\(\mathfrak{pe}(2)_{-}\).

\subsection{Binary differential operators and the negative part}
\label{subsec:binary-operators}

Let
\[
J:
\mathcal F_{\lambda_1,\nu_1}
\otimes
\mathcal F_{\lambda_2,\nu_2}
\longrightarrow
\mathcal F_{\mu,\omega}
\]
be a binary differential operator. By
Lemma~\ref{lem:XQT-normal-form}, \(J\) can be expanded in binary
monomials constructed from~$\Dxi$, $\Dq$, $\partial_\tau$.

For homogeneous differential operators~\(P\) and~\(R\), we use the
graded tensor-product notation
\[
(P\boxtimes R)(\varphi,\psi)
:=
(-1)^{p(R)p(\varphi)}
P(\varphi)R(\psi),
\]
extended by linearity. Thus, a general binary differential operator
can be written as
\[
J
=
\sum_{\alpha,\beta}
F_{\alpha,\beta}(q,\xi,\tau)\,
P_\alpha\boxtimes P_\beta,
\]
where~\(P_\alpha,P_\beta\) are monomials in~$\Dxi$, $\Dq$, $\partial_\tau$.

The following elementary identity allows us to reduce the action on
binary monomials to the corresponding unary actions.

\ssbegin[{The Leibniz rule for $P\boxtimes R$}]{Lemma}[The Leibniz rule for~$P\boxtimes R$]
\label{lem:boxtimes-action}
For every homogeneous~\(Y,P,R\),
\[
Y\cdot(P\boxtimes R)
=
(Y\cdot P)\boxtimes R
+
(-1)^{p(Y)p(P)}
P\boxtimes(Y\cdot R).
\]
\end{Lemma}

The proof is in the Appendix, \S\ref{app:proofs}.

We now impose invariance under the negative part
\[
\mathfrak{pe}(2)_-
=
\langle M_1,M_q,M_\xi\rangle.
\]
The choice of the~\(\Dxi,\Dq,\partial_\tau\)-basis is particularly useful here: these
operators are adapted to the negative generators, and invariance
forces the coefficient functions in the preceding expansion to be
constant.

\ssbegin[{Invariance under $\fpe(2)_-$ forces constant coefficients}]{Lemma}[Invariance under~$\fpe(2)_-$ forces constant coefficients]
\label{lem:negative-part-constant-coefficients}
Let
\[
J
=
\sum_{\alpha,\beta}
F_{\alpha,\beta}(q,\xi,\tau)\,
P_\alpha\boxtimes P_\beta
\]
be a binary differential operator in the~\(\Dxi,\Dq,\partial_\tau\)-basis.
If~\(J\) is~\(\mathfrak{pe}(2)_-\)-invariant, then~%
$
F_{\alpha,\beta}\in\mathbb C
$
for every pair~\((\alpha,\beta)\).
\end{Lemma}

The proof is in the Appendix, \S\ref{app:proofs}.

Consequently, every~\(\mathfrak{pe}(2)_-\)-invariant binary
differential operator is a constant-coefficient linear combination of
the monomials
\[
\Dxi^a\Dq^\varepsilon \partial_\tau^\delta
\boxtimes
\Dxi^b\Dq^{\varepsilon'}\partial_\tau^{\delta'}.
\]
Thus the classification problem is reduced from arbitrary
polynomial-coefficient binary differential operators to a finite
linear combination of tensor monomials of prescribed degree. We next
impose invariance under the maximal toral subalgebra of~%
\(\mathfrak{pe}(2)\).

\section{$\mathfrak{pe}(2)$-invariant binary differential operators}
\label{sec:invariant-binary}
The aim of this section is to classify the bilinear differential operators
\[
J:
\mathcal F_{\lambda_1,\nu_1}\otimes
\mathcal F_{\lambda_2,\nu_2}
\longrightarrow
\mathcal F_{\mu,\omega}
\]
that are invariant under the action of the periplectic Lie superalgebra~%
$
\mathfrak{pe}(2)\subset\mathfrak m(1).
$ More precisely, we determine the values of the weights for which such operators exist and describe all~$\mathfrak{pe}(2)$-invariant operators explicitly.
The proofs in this Section are computations, and are collected in the
Appendix, \S\ref{app:proofs}. Since~$M_{q\tau}=\frac12[M_{q^2},M_{\tau\xi}]$, invariance
under~$\fpe(2)_-$, $H_1$, $H_2$, $M_{\tau\xi}$ and~$M_{q^2}$ is invariance under the whole
of~$\fpe(2)$.

\subsection{Decomposition by the weights of the maximal torus}
\label{subsec:cartan-decomposition}

We now impose invariance under the maximal toral subalgebra of~%
\(\mathfrak{pe}(2)\), generated by
\[
H_1:=\frac12 M_{\tau+q\xi},
\qquad
H_2:=\frac12 M_{\tau-q\xi}.
\]
To determine the corresponding weight conditions, consider
differential monomials
\[
P_i=\Dxi^{a_i}\Dq^{\varepsilon_i}\partial_\tau^{\delta_i},
\qquad i=1,2,
\]
where
\[
a_i\in\mathbb Z_{\geq0},
\qquad
\varepsilon_i,\delta_i\in\{0,1\}.
\]
They determine the binary differential monomial
\[
P_1\boxtimes P_2:
\mathcal F_{\lambda_1,\nu_1}
\otimes
\mathcal F_{\lambda_2,\nu_2}
\longrightarrow
\mathcal F_{\mu,\omega}.
\]
For later convenience, set
\[
u:=a_1+a_2+\delta_1+\delta_2,
\qquad
v:=\varepsilon_1+\varepsilon_2+\delta_1+\delta_2.
\]

\ssbegin[{The toral weights fix the target weight $(\mu,\omega)$}]{Lemma}[The toral weights fix the target weight~$(\mu,\omega)$]
\label{lem:binary-PBW-Cartan-shift}
If~\(P_1\boxtimes P_2\) intertwines the action of the maximal toral
subalgebra generated by~\(H_1\) and~\(H_2\), then necessarily
\[
\mu=\lambda_1+\lambda_2-u,
\quad \text{ and } \quad
\omega
=
\nu_1+\nu_2+\frac{u-v}{2}.
\]
\end{Lemma}

Since~%
$
\deg \Dxi=\deg \Dq=1,
\,
\deg \partial_\tau=2,
$ the total degree of~$P_1\boxtimes P_2$ is
\[
\begin{array}{r@{}l}
k
&{}=
a_1+a_2+\varepsilon_1+\varepsilon_2
+2\delta_1+2\delta_2
=
\left(a_1+a_2+\delta_1+\delta_2\right)
+
\left(
\varepsilon_1+\varepsilon_2+\delta_1+\delta_2
\right).
\end{array}
\]
The integer~$v$ also determines the parity of the monomial.
Indeed,
\[
p\left(
P_1\boxtimes P_2
\right)
\equiv
\varepsilon_1+\varepsilon_2+\delta_1+\delta_2
\equiv v
\pmod 2.
\]
Since each monomial factor contains at most one~$\Dq$ and at most one~$\partial_\tau$,
we have
\[
0\leq v\leq4.
\]
Consequently, for an odd binary operator the only possible Cartan
weights are
\[
v=1,\qquad v=3,
\]
whereas for an even binary operator the only possibilities are
\[
v=0,\qquad v=2,\qquad v=4.
\]
This gives the following decomposition by the toral weights.

For odd operators:
\[
\begin{array}{c|c|c|c}
\text{Type} & v & u & (\mu,\omega)
\\ \hline
A
&
1
&
k-1
&
\left(
\lambda_1+\lambda_2-k+1,\,
\nu_1+\nu_2+\dfrac{k-2}{2}
\right)
\\[3mm]
B
&
3
&
k-3
&
\left(
\lambda_1+\lambda_2-k+3,\,
\nu_1+\nu_2+\dfrac{k-6}{2}
\right).
\end{array}
\]

For even operators:
\[
\begin{array}{c|c|c|c}
\text{Type} & v & u & (\mu,\omega)
\\ \hline
A
&
0
&
k
&
\left(
\lambda_1+\lambda_2-k,\,
\nu_1+\nu_2+\dfrac{k}{2}
\right)
\\[3mm]
B
&
2
&
k-2
&
\left(
\lambda_1+\lambda_2-k+2,\,
\nu_1+\nu_2+\dfrac{k-4}{2}
\right)
\\[3mm]
C
&
4
&
k-4
&
\left(
\lambda_1+\lambda_2-k+4,\,
\nu_1+\nu_2+\dfrac{k-8}{2}
\right).
\end{array}
\]
The classification problem splits into five independent
sectors: two odd sectors, corresponding to~$v=1$ and~$v=3$, and three
even sectors, corresponding to~$v=0,2,4$.

\subsection{Odd operators of Type A: the generic case}
\label{subsubsec:odd-typeA}
We first consider odd operators of Type A and restrict ourselves to
generic values of the weights. The exceptional values, for which the
coefficient system loses rank, will be studied separately.

For Type A, one has~%
$
v=1$ and~$u=k-1.$ Hence the Cartan equations determine the target weights as
\[
\mu=\lambda_1+\lambda_2-k+1,
\qquad
\omega=\nu_1+\nu_2+\frac{k-2}{2}.
\]
Since~\(v=1\), the monomials that may occur contain exactly one
contribution to the second Cartan degree. There are four possibilities:
\(\Dq\) or~\(\partial_\tau\) may occur in either the first or the second tensor
factor.
\subsubsection{Odd Type A operators in low degrees}
Before considering arbitrary degree, we record the first three cases
explicitly. Besides, providing the low-degree classification, these
examples show that~\(k=2\) is exceptional, whereas the general
recurrence already applies from~\(k=3\).

{\bf Degree~\(k=1\).} The target weights are~%
$
\mu=\lambda_1+\lambda_2,
$ and~%
$\omega=\nu_1+\nu_2-\frac12.
$ The most general Type A operator is
\[
J_{\bar1,A}^{(1)}
=
A_0\,1\boxtimes \Dq
+
E_1\,\Dq\boxtimes1.
\]
There is no condition coming from~\(M_{q^2}\). Invariance under~%
\(M_{\tau\xi}\) gives
\[
\lambda_2A_0+\lambda_1E_1=0.
\]
Hence, for generic weights, the space of Type A invariants is
one-dimensional. A convenient generator is
\[
J_{\bar1,A}^{(1)}
=
\lambda_1\,1\boxtimes \Dq
-
\lambda_2\,\Dq\boxtimes1.
\]
Thus~%
$
\dim\operatorname{Inv}_{\bar1,A}^{(1)}=1
$ at generic weights. The exceptional case~%
$
\lambda_1=\lambda_2=0
$ will be considered together with the resonant cases.

{\bf Degree~\(k=2\).} The target weights are~%
$
\mu=\lambda_1+\lambda_2-1,
$ and~$
\omega=\nu_1+\nu_2.
$ The most general Type A operator is
\[
\begin{array}{r@{}l}
J_{\bar1,A}^{(2)}
={}&
A_0\,1\boxtimes \Dxi\Dq
+
A_1\,\Dxi\boxtimes \Dq
+
B_0\,1\boxtimes \partial_\tau
+
E_1\,\Dq\boxtimes \Dxi
+
E_2\,\Dxi\Dq\boxtimes1
+
F_2\,\partial_\tau\boxtimes1.
\end{array}
\]
Invariance under~\(M_{q^2}\) gives~%
$E_1=A_1.$ The~\(M_{\tau\xi}\)-equations are
\[
B_0
=
2\lambda_2A_0+2\lambda_1A_1,
\quad \text{ and } \quad F_2
=
2\lambda_2A_1+2\lambda_1E_2.
\]
Thus~\(A_0,A_1,E_2\) are independent parameters, and every Type A
operator of degree~\(2\) can be written as
\[
J_{\bar1,A}^{(2)}
=
aJ_1^{(2)}+bJ_2^{(2)}+cJ_3^{(2)},
\]
where
\[
J_1^{(2)}
=
1\boxtimes \Dxi\Dq
+
2\lambda_2\,1\boxtimes \partial_\tau,
\]
\[
J_2^{(2)}
=
\Dxi\boxtimes \Dq
+
\Dq\boxtimes \Dxi
+
2\lambda_1\,1\boxtimes \partial_\tau
+
2\lambda_2\,\partial_\tau\boxtimes1,
\]
and
\[
J_3^{(2)}
=
\Dxi\Dq\boxtimes1
+
2\lambda_1\,\partial_\tau\boxtimes1.
\]
Consequently,
\[
\dim\operatorname{Inv}_{\bar1,A}^{(2)}=3.
\]
Thus~\(k=2\) is exceptional and must be treated separately from
the general Type A family.

{\bf Degree~\(k=3\).} The target weights are~%
$
\mu=\lambda_1+\lambda_2-2,
$ and~$\omega=\nu_1+\nu_2+\frac12.$ The general Type A operator is
\[
\begin{array}{r@{}l}
J_{\bar1,A}^{(3)}
={}&
A_0\,1\boxtimes \Dxi^2\Dq
+
A_1\,\Dxi\boxtimes \Dxi\Dq
+
A_2\,\Dxi^2\boxtimes \Dq
+
B_0\,1\boxtimes \Dxi\partial_\tau
+
B_1\,\Dxi\boxtimes \partial_\tau
\\
&{}+
E_1\,\Dq\boxtimes \Dxi^2
+
E_2\,\Dxi\Dq\boxtimes \Dxi
+
E_3\,\Dxi^2\Dq\boxtimes1
+
F_2\,\partial_\tau\boxtimes \Dxi
+
F_3\,\Dxi\partial_\tau\boxtimes1.
\end{array}
\]
Invariance under~\(M_{q^2}\) gives
\[
B_0=2A_0,
\qquad
B_1=A_1,
\qquad
E_1=\frac12A_1,
\qquad
E_2=2A_2,
\qquad
F_2=E_2,
\qquad
F_3=2E_3.
\]
The remaining~\(M_{\tau\xi}\)-equations reduce to
\[
2\lambda_1A_1
+
2(2\lambda_2-2)A_0
=
0,
\qquad
2(2\lambda_1-1)A_2
+
(2\lambda_2-1)A_1
=
0,
\]
and the endpoint equation
\[
(2-2\lambda_1)E_3
=
2\lambda_2A_2.
\]
At generic weights these equations determine all coefficients from~\(A_0\). Explicitly,
\[
A_1
=
-\frac{2(\lambda_2-1)}{\lambda_1}\,A_0,
\quad
A_2
=
\frac{
(\lambda_2-1)(2\lambda_2-1)
}{
\lambda_1(2\lambda_1-1)
}\,A_0,
\quad { and } \quad
E_3
=
\frac{\lambda_2}{1-\lambda_1}\,A_2.
\]
Therefore,
\[
\dim\operatorname{Inv}_{\bar1,A}^{(3)}=1
\]
at generic weights.

The exceptional values at which one or more of the preceding
relations degenerates will be included in the resonance analysis.
\subsubsection{Odd Type A operators in degree $k$}
The most general odd Type A operator of degree~%
\(k\) is
\begin{align}
J_{\bar1,A}^{(k)}
={}&
\sum_{0\leq i\leq k-1}
A_i\,\Dxi^i\boxtimes \Dxi^{k-i-1}\Dq
+
\sum_{0\leq i\leq k-2}
B_i\,\Dxi^i\boxtimes \Dxi^{k-i-2}\partial_\tau
\nonumber+
\sum_{1\leq i\leq k}
E_i\,\Dxi^{i-1}\Dq\boxtimes \Dxi^{k-i}
\\& +
\sum_{2\leq i\leq k}
F_i\,\Dxi^{i-2}\partial_\tau\boxtimes \Dxi^{k-i}.
\label{eq:odd-typeA-ansatz}
\end{align}
Imposing invariance under~\(M_{\tau\xi}\) and~\(M_{q^2}\) turns~%
\eqref{eq:odd-typeA-ansatz} into a~recurrence for the coefficients~\(A_i\),
\(B_i\), \(E_i\), \(F_i\); solving that recurrence gives the following generic
answer.

\ssbegin[{Odd Type A operators for generic weights}]{Theorem}[Odd Type A operators for generic weights]
\label{prop:odd-typeA-generic}
Let~\(k\) be in the stable range for the above coefficient system, and
assume that the weights are generic. Then, the space of odd Type A~%
\(\mathfrak{pe}(2)\)-invariant binary differential operators
\[
J_{\bar1,A}^{(k)}:\mathcal F_{\lambda_1,\nu_1}
\otimes
\mathcal F_{\lambda_2,\nu_2}
\longrightarrow
\mathcal F_{\lambda_1+\lambda_2-k+1,\,
\nu_1+\nu_2+\frac{k-2}{2}}
\]
given as in the ansatz~\eqref{eq:odd-typeA-ansatz}
is one-dimensional.

Up to an overall scalar, its coefficients are determined by
\[
A_i
=
(-1)^i
\binom{k-1}{i}
\frac{
\displaystyle
\prod_{0\leq r\leq i-1}
(2\lambda_2-k+r+1)
}{
\displaystyle
\prod_{0\leq r\leq i-1}
(2\lambda_1-r)
}
A_0,
\qquad
0\leq i\leq k-1,
\]
\[
B_i=(k-i-1)A_i, \; (0\leq i \leq k-2)
\qquad
E_i=\frac{i}{k-i}A_i,\; (1\leq i\leq k-1)\]
and
\[
F_i=(i-1)E_i,\; (2\leq i \leq k)
\quad
E_k
=
\frac{2\lambda_2}
{k-1-2\lambda_1}
A_{k-1}.
\]
\end{Theorem}
\subsection{Odd operators of Type A: the resonant case}
\label{subsubsec:odd-typeA-resonant}
\subsubsection{Resonant odd Type A operators in low degrees}
We now examine the exceptional weights for~$k=1,2,3$ by studying the rank of the coefficient system directly, without dividing by weight-dependent factors, since dimension jumps may occur even when no individual equation vanishes identically.

{\bf Degree~\(k=1\).} Recall that
\[
J_{\bar1,A}^{(1)}
=
A_0\,1\boxtimes \Dq
+
E_1\,\Dq\boxtimes1.
\]
The only invariance equation is
\[
\lambda_2A_0+\lambda_1E_1=0.
\]
At~%
$
(\lambda_1,\lambda_2)=(0,0),
$ the equation vanishes identically. Both~\(A_0\) and~\(E_1\) are then
arbitrary, and therefore
\[
\dim\operatorname{Inv}_{\bar1,A}^{(1)}=2.
\]
A basis at this resonant point is
\[
1\boxtimes \Dq,
\qquad
\Dq\boxtimes1.
\]
Consequently,
\[
\dim\operatorname{Inv}_{\bar1,A}^{(1)}
=
\begin{cases}
2,&(\lambda_1,\lambda_2)=(0,0),\\
1,&\text{otherwise}.
\end{cases}
\]

{\bf Degree~\(k=2\).} The general Type A operator is
\[
\begin{array}{r@{}l}
J_{\bar1,A}^{(2)}
={}&
A_0\,1\boxtimes \Dxi\Dq
+
A_1\,\Dxi\boxtimes \Dq
+
B_0\,1\boxtimes \partial_\tau
+
E_1\,\Dq\boxtimes \Dxi
+
E_2\,\Dxi\Dq\boxtimes1
+
F_2\,\partial_\tau\boxtimes1.
\end{array}
\]
The complete coefficient system is
\[ E_1=A_1, \quad
B_0=2\lambda_2A_0+2\lambda_1A_1, \quad \text{ and }\quad
F_2=2\lambda_2A_1+2\lambda_1E_2.
\]
The coefficients~%
$
A_0, A_1, E_2
$  are arbitrary, while~\(E_1,B_0,F_2\) are uniquely determined by them.
In particular, the rank of this system is independent of~%
\((\lambda_1,\lambda_2)\). Therefore, no resonance occurs in degree~%
\(k=2\), and
\[
\dim\operatorname{Inv}_{\bar1,A}^{(2)}=3
\]
for every pair of weights~\((\lambda_1,\lambda_2)\).

{\bf Degree~\(k=3\).} After using the relations obtained from~\(M_{q^2}\)-invariance, the
independent coefficients may be chosen to be~%
$
A_0, A_1, A_2, E_3.$
The remaining invariance equations are
\[
2(2\lambda_2-2)A_0
+
2\lambda_1A_1
=
0, \quad (2\lambda_2-1)A_1
+
2(2\lambda_1-1)A_2
=
0,
\quad
-2\lambda_2A_2
+
(2-2\lambda_1)E_3
=
0.
\]
The system can be written as
\[
\begin{pmatrix}
4(\lambda_2-1) & 2 \lambda_1 & 0 & 0\\
0 & 2\lambda_2-1 & 2(2 \lambda_1-1) & 0\\
0 & 0 & - 2 \lambda_2 & 2(1-\lambda_1)
\end{pmatrix}
\begin{pmatrix}
A_0\\
A_1\\
A_2\\
E_3
\end{pmatrix}
=0.
\]
We now determine the values for which the rank drops. A direct
analysis of the above matrix shows that this happens precisely at
\[
(\lambda_1,\lambda_2)
\in
\left\{
(0,1),\
\left(\frac12,\frac12\right),\
(1,0),\
\left(\frac12,1\right),\
\left(1,\frac12\right),\
(1,1)
\right\}.
\]
At each of these six points the matrix has rank~\(2\), and hence
\[
\dim\operatorname{Inv}_{\bar1,A}^{(3)}=2.
\]
Let us describe the degenerations explicitly.

At~%
$
(\lambda_1,\lambda_2)=(0,1),
$ the system gives~%
$
A_1=2A_2,$ and~$E_3=A_2.
$
Thus~\(A_0\) and~\(A_2\) are independent, and the solution space is two-dimensional.

At~%
$
(\lambda_1,\lambda_2)
=
\left(\frac12,\frac12\right),$ the system gives~%
$
A_1=2A_0,
$ and~$
E_3=A_2.
$
Again the solution space is two-dimensional.

At~%
$
(\lambda_1,\lambda_2)=(1,0),
$ the system gives~%
$
A_1=2A_0,$ and~$
A_2=A_0.
$
Thus~\(A_0\) and~\(E_3\) are independent, and the space is two-dimensional.

There are three further resonant points at which no complete row needs
to vanish. Instead, the rank drops because adjacent equations cease to
impose independent conditions.

At~%
$
(\lambda_1,\lambda_2)
=
\left(\frac12,1\right),
$ the system gives~%
$
A_1=0,
$ and~$
E_3=2A_2,
$ while~\(A_0\) and~\(A_2\) are arbitrary.

At~%
$
(\lambda_1,\lambda_2)
=
\left(1,\frac12\right),
$ one obtains~%
$
A_1=A_0,
$ and~$
A_2=0,
$ while~\(A_0\) and~\(E_3\) are arbitrary.

Finally, at~%
$
(\lambda_1,\lambda_2)=(1,1),
$
the system gives~%
$
A_1=A_2=0,
$ while~\(A_0\) and~\(E_3\) remain independent.

We therefore obtain
\[
\dim\operatorname{Inv}_{\bar1,A}^{(3)}
=
\begin{cases}
2,
&
(\lambda_1,\lambda_2)
\in
\left\{
(0,1),
\left(\frac12,\frac12\right),
(1,0),
\left(\frac12,1\right),
\left(1,\frac12\right),
(1,1)
\right\},
\\[4mm]
1,
&
\text{otherwise}.
\end{cases}
\]

\paragraph{Summary of the low-degree cases.}

The odd Type A invariant spaces in degrees~\(k=1,2,3\) are therefore
described by
\[
\begin{array}{c|c|c}
k & \text{exceptional weights} & \dim\operatorname{Inv}_{\bar1,A}^{(k)}
\\ \hline
1
&
(0,0)
&
2\text{ at }(0,0),\quad 1\text{ otherwise}
\\[2mm]
2
&
\text{none}
&
3\text{ for all weights}
\\[2mm]
3
&
\begin{array}{c}
(0,1),\ (\frac12,\frac12),\ (1,0),\\[1mm]
(\frac12,1),\ (1,\frac12),\ (1,1)
\end{array}
&
2\text{ at the listed points},\quad 1\text{ otherwise}.
\end{array}
\]
\subsubsection{Resonant odd Type A operators for \(k\geq3\)}
We now determine all exceptional weights for odd Type A operators in
the stable range~\(k\geq3\). The resonance problem is treated as a
rank problem for the complete coefficient system; in particular, we
do not divide by any weight-dependent coefficient.

The resonant set is described as follows.

\ssbegin[The resonant weights in the odd sector~A]{Proposition}[The resonant weights in the
odd sector~$A$]\label{prop:odd-typeA-resonance}
Let~\(k\geq3\). The odd Type A coefficient system is resonant
precisely when
\[
2\lambda_1,2\lambda_2
\in
\{0,1,\ldots,k-1\}
\quad
\text{ and }
\quad
2\lambda_1+2\lambda_2\geq k-1.
\]
At every such point,
\[
\dim\operatorname{Inv}_{\bar1,A}^{(k)}=2.
\]
At all other weights,
\[
\dim\operatorname{Inv}_{\bar1,A}^{(k)}=1.
\]
\end{Proposition}


\subsection{Odd operators of Type B: the generic case}
\label{subsubsec:odd-typeB}
\subsubsection{Odd operators of Type B in low degrees}

We now turn to odd operators of Type B. Recall that Type B corresponds
to the Cartan degree~%
$
v=3, u=k-3,
$
and therefore the target weights are
\[
\mu=\lambda_1+\lambda_2-k+3,
\qquad
\omega=\nu_1+\nu_2+\frac{k-6}{2}.
\]

The general odd Type B ansatz is
\begin{align}
J_{\bar1,B}^{(k)}
={}&
\sum_{3\leq i\leq k-1}
C_i\,\Dxi^{i-3}\Dq\partial_\tau\boxtimes \Dxi^{k-i-1}\Dq
+
\sum_{3\leq i\leq k-2}
D_i\,\Dxi^{i-3}\Dq\partial_\tau\boxtimes \Dxi^{k-i-2}\partial_\tau
\nonumber\\
&+
\sum_{1\leq i\leq k-3}
G_i\,\Dxi^{i-1}\Dq\boxtimes \Dxi^{k-i-3}\Dq\partial_\tau
+
\sum_{2\leq i\leq k-3}
H_i\,\Dxi^{i-2}\partial_\tau\boxtimes \Dxi^{k-i-3}\Dq\partial_\tau.
\label{eq:odd-typeB-general-ansatz}
\end{align}

We first explain why no nonzero Type B operator exists in degrees~\(k=1,2,3,4\).

{\bf Degree~\(k=1,2,3,4\).} In this case, the condition~\(v=3\) is incompatible with the
existence of a monomial of total degree~\(k\). Equivalently, all
the sums in the ansatz~\eqref{eq:odd-typeB-general-ansatz} are empty. Therefore,
\[
J_{\bar1,B}^{(k)}=0,
\qquad
k=1,2,3.
\]
The case~\(k=4\) requires a separate argument. Here, the Type B ansatz
contains exactly two terms:
\[
J_{\bar1,B}^{(4)}
=
C_3\,\Dq\partial_\tau\boxtimes \Dq
+
G_1\,\Dq\boxtimes \Dq\partial_\tau.
\]
We impose invariance under~\(M_{\tau\xi}\) and obtain
\[
\begin{array}{r@{}l}
M_{\tau\xi}\cdot
\left(\Dq\partial_\tau\boxtimes \Dq\right)
={}&
\Dxi\Dq\boxtimes \Dq
+
2\lambda_1\partial_\tau\boxtimes \Dq
+
2\lambda_2\Dq\partial_\tau\boxtimes1,
\\
M_{\tau\xi}\cdot
\left(\Dq\boxtimes \Dq\partial_\tau\right)
={}&
2\lambda_1\,1\boxtimes \Dq\partial_\tau
-
\Dq\boxtimes \Dxi\Dq
-
2\lambda_2\Dq\boxtimes \partial_\tau.
\end{array}
\]
Consequently,
\[
\begin{array}{r@{}l}
M_{\tau\xi}\cdot J_{\bar1,B}^{(4)}
={}&
C_3
\left(
\Dxi\Dq\boxtimes \Dq
+
2\lambda_1\partial_\tau\boxtimes \Dq
+
2\lambda_2\Dq\partial_\tau\boxtimes1
\right)
\\
&{}+
G_1
\left(
2\lambda_1\,1\boxtimes \Dq\partial_\tau
-
\Dq\boxtimes \Dxi\Dq
-
2\lambda_2\Dq\boxtimes \partial_\tau
\right).
\end{array}
\]
Comparing the coefficients of the distinct monomials~$\Dxi\Dq\boxtimes\Dq$
and~$\Dq\boxtimes\Dxi\Dq$, we get~$C_3=G_1=0$. Thus,
\[
J_{\bar1,B}^{(4)}=0.
\]
Hence the first degree in which a nonzero odd Type B operator may
occur is~\(k=5\). We shall therefore study explicitly the cases~%
\(k=5\) and~\(k=6\). The degree~\(k=5\) gives the first nonzero Type B
operator, while~\(k=6\) is the first degree in which the coefficients
satisfy a nontrivial recurrence.

{\bf Degree~\(k=5\).} In this case, the odd Type B operator has the form
\[
\begin{array}{r@{}l}
J_{\bar1,B}^{(5)}
={}&
C_3\,\Dq\partial_\tau\boxtimes \Dxi\Dq
+
C_4\,\Dxi\Dq\partial_\tau\boxtimes \Dq
+
D_3\,\Dq\partial_\tau\boxtimes \partial_\tau
\\
&{}+
G_1\,\Dq\boxtimes \Dxi\Dq\partial_\tau
+
G_2\,\Dxi\Dq\boxtimes \Dq\partial_\tau
+
H_2\,\partial_\tau\boxtimes \Dq\partial_\tau.
\end{array}
\]
A direct computation of the~\(M_{\tau\xi}\)-invariance condition gives
\[
C_4=G_1=0,
\qquad
G_2=C_3,
\qquad
D_3=2\lambda_2C_3,
\qquad
H_2=2\lambda_1C_3.
\]
Hence, up to an overall scalar,
\[
J_{\bar1,B}^{(5)}
=
\Dq\partial_\tau\boxtimes \Dxi\Dq
+
2\lambda_2\,\Dq\partial_\tau\boxtimes \partial_\tau
+
\Dxi\Dq\boxtimes \Dq\partial_\tau
+
2\lambda_1\,\partial_\tau\boxtimes \Dq\partial_\tau.
\]

A direct computation shows that this operator is also invariant under~%
\(M_{q^2}\), and no further condition on the weights is obtained.
Therefore,
\[
\dim\operatorname{Inv}_{\bar1,B}^{(5)}=1
\]
for all~\((\lambda_1,\lambda_2)\).
\subsubsection{Odd Type B operators: the generic case for \(k\geq 6\)}
The most general odd Type B operator is
\begin{align}
J_{\bar1,B}^{(k)}
={}&
\sum_{3\leq i\leq k-1}
C_i\,\Dxi^{i-3}\Dq\partial_\tau\boxtimes \Dxi^{k-i-1}\Dq
+
\sum_{3\leq i\leq k-2}
D_i\,\Dxi^{i-3}\Dq\partial_\tau\boxtimes \Dxi^{k-i-2}\partial_\tau
\nonumber\\
&+
\sum_{1\leq i\leq k-3}
G_i\,\Dxi^{i-1}\Dq\boxtimes \Dxi^{k-i-3}\Dq\partial_\tau
+
\sum_{2\leq i\leq k-3}
H_i\,\Dxi^{i-2}\partial_\tau\boxtimes \Dxi^{k-i-3}\Dq\partial_\tau.
\label{eq:typeB-general-ansatz}
\end{align}

Invariance under~\(M_{\tau\xi}\) and~\(M_{q^2}\) reduces~%
\eqref{eq:typeB-general-ansatz} to a~bidiagonal recurrence for the coefficients~%
\(C_i\); its solution gives the following generic classification.

\ssbegin[{Odd Type B operators for generic weights}]{Theorem}[Odd Type B operators for generic weights]
\label{prop:typeB-generic}
Let~\(k\geq6\), and assume that the weights are generic so that the
coefficient system below does not lose rank.

Then, the space of odd Type B~%
\(\mathfrak{pe}(2)\)-invariant binary differential operators
\[
J_{\bar1,B}^{(k)}:\mathcal F_{\lambda_1,\nu_1}
\otimes
\mathcal F_{\lambda_2,\nu_2}
\longrightarrow
\mathcal F_{\lambda_1+\lambda_2-k+3,\,
\nu_1+\nu_2+\frac{k-6}{2}}
\]
is one-dimensional.

Up to an overall scalar, the operator is
\begin{align*}
J_{\bar1,B}^{(k)}
={}&
\sum_{3\leq i\leq k-2}
C_i
\Bigl(
\Dxi^{i-3}\Dq\partial_\tau\boxtimes \Dxi^{k-i-1}\Dq
+
2\lambda_2
\Dxi^{i-3}\Dq\partial_\tau\boxtimes \Dxi^{k-i-2}\partial_\tau
\nonumber\\
&\qquad\qquad
+
\Dxi^{i-2}\Dq\boxtimes \Dxi^{k-i-2}\Dq\partial_\tau
+
2\lambda_1
\Dxi^{i-3}\partial_\tau\boxtimes \Dxi^{k-i-2}\Dq\partial_\tau
\Bigr),
\end{align*}
where the coefficients~\(C_i\) satisfy
\[
(i-2)(2\lambda_1-i+1)C_{i+1}
+
(k-i-2)(2\lambda_2-k+i+1)C_i
=
0.
\]
Equivalently,
\[
C_i
=
(-1)^{i-3}
\binom{k-5}{i-3}
\left(
\prod_{3\leq r\leq i-1}
(2\lambda_2-k+r+1)
\right)
\left(
\prod_{i\leq r\leq k-3}
(2\lambda_1-r+1)
\right),
\qquad
3\leq i\leq k-2.
\]
\end{Theorem}
\subsubsection{Low-degree resonances for odd Type B operators}
Before studying the resonant Type B system for arbitrary~\(k\), we
consider the first two nonzero cases, namely~\(k=5\) and~\(k=6\).

{\bf Degree~\(k=5\).} After imposing~\(M_{\tau\xi}\)-invariance, every Type B operator of
degree~\(5\) is determined by a single coefficient~\(C_3\) and has the
form
\[
\begin{array}{r@{}l}
J_{\bar1,B}^{(5)}
=
C_3\Bigl(
&\Dq\partial_\tau\boxtimes \Dxi\Dq
+
2\lambda_2\,\Dq\partial_\tau\boxtimes \partial_\tau
+
\Dxi\Dq\boxtimes \Dq\partial_\tau
+
2\lambda_1\,\partial_\tau\boxtimes \Dq\partial_\tau
\Bigr).
\end{array}
\]

The~\(M_{q^2}\)-invariance equations are governed, in the stable
Type B system, by the recurrence
\[
(i-2)(2\lambda_1-i+1)C_{i+1}
+
(k-i-2)(2\lambda_2-k+i+1)C_i
=
0,
\quad
\text{ with }
\quad
3\leq i\leq k-3.
\]
For~\(k=5\), this index range is empty. Consequently,
\(M_{q^2}\)-invariance imposes no further condition on~\(C_3\).
Therefore, \(C_3\) is arbitrary for every pair~%
\((\lambda_1,\lambda_2)\), and
\[
\dim\operatorname{Inv}_{\bar1,B}^{(5)}=1.
\]
Thus there are no resonant weights in degree~\(5\): the dimension of
the Type B invariant space is equal to~\(1\) for all weights.

{\bf Degree~\(k=6\).} After imposing~\(M_{\tau\xi}\)-invariance, the Type B
operator depends on two coefficients~\(C_3\) and~\(C_4\):
\[
\begin{array}{r@{}l}
J_{\bar1,B}^{(6)}
={}&
C_3\,\Dq\partial_\tau\boxtimes \Dxi^2\Dq
+
C_4\,\Dxi\Dq\partial_\tau\boxtimes \Dxi\Dq
+
2\lambda_2C_3\,\Dq\partial_\tau\boxtimes \Dxi\partial_\tau
+
2\lambda_2C_4\,\Dxi\Dq\partial_\tau\boxtimes \partial_\tau
\\
&{}+
C_3\,\Dxi\Dq\boxtimes \Dxi\Dq\partial_\tau
+
C_4\,\Dxi^2\Dq\boxtimes \Dq\partial_\tau
+
2\lambda_1C_3\,\partial_\tau\boxtimes \Dxi\Dq\partial_\tau
+
2\lambda_1C_4\,\Dxi\partial_\tau\boxtimes \Dq\partial_\tau.
\end{array}
\]
In this degree, the~\(M_{q^2}\)-invariance condition consists of the
single equation
\[
(2\lambda_1-2)C_4
+
(2\lambda_2-2)C_3
=
0.
\]
If~%
$
(\lambda_1,\lambda_2)\neq(1,1),
$
this is a nontrivial linear relation between~\(C_3\) and~\(C_4\).
Hence,
\[
\dim\operatorname{Inv}_{\bar1,B}^{(6)}=1.
\]
At the exceptional point~%
$
(\lambda_1,\lambda_2)=(1,1),
$
both coefficients~\(C_3\) and~\(C_4\) become independent, and
\[
\dim\operatorname{Inv}_{\bar1,B}^{(6)}=2.
\]
A convenient basis of the two-dimensional resonant space is obtained
by taking successively
\[
(C_3,C_4)=(1,0)
\qquad\text{and}\qquad
(C_3,C_4)=(0,1).
\]
This gives
\[
\begin{array}{r@{}l}
J_{\bar1,B;1}^{(6)}
={}&
\Dq\partial_\tau\boxtimes \Dxi^2\Dq
+
2\Dq\partial_\tau\boxtimes \Dxi\partial_\tau
+
\Dxi\Dq\boxtimes \Dxi\Dq\partial_\tau
+
2\partial_\tau\boxtimes \Dxi\Dq\partial_\tau,
\end{array}
\]
and
\[
\begin{array}{r@{}l}
J_{\bar1,B;2}^{(6)}
={}&
\Dxi\Dq\partial_\tau\boxtimes \Dxi\Dq
+
2\Dxi\Dq\partial_\tau\boxtimes \partial_\tau
+
\Dxi^2\Dq\boxtimes \Dq\partial_\tau
+
2\Dxi\partial_\tau\boxtimes \Dq\partial_\tau.
\end{array}
\]
We therefore obtain
\[
\dim\operatorname{Inv}_{\bar1,B}^{(6)}
=
\begin{cases}
2,
&
(\lambda_1,\lambda_2)=(1,1),
\\[1mm]
1,
&
\text{otherwise}.
\end{cases}
\]
The cases~\(k=5\) and~\(k=6\) illustrate the beginning of the general
pattern. In degree~\(5\), there is only one coefficient and hence no
possible rank drop. Starting from degree~\(6\), the Type B coefficients
form a bidiagonal system, and resonant weights are precisely those at
which the rank of this system decreases. We now analyze this
phenomenon for arbitrary~\(k\geq6\).
\subsubsection{Resonant odd Type B operators for \(k\geq6\)}
We now determine the exceptional weights for odd Type B operators in
arbitrary degree~\(k\geq6\).

The answer is the following.

\ssbegin[{The resonant weights of odd Type B operators}]{Proposition}[The resonant weights of odd Type B operators]
\label{prop:typeB-resonance}
Let~\(k\geq6\). The odd Type B coefficient system is resonant
precisely when
\[
2\lambda_1,2\lambda_2
\in
\{2,3,\ldots,k-4\}
\quad \text{ and }
\quad
2\lambda_1+2\lambda_2\geq k-2.
\]
At every resonant point,
\[
\dim\operatorname{Inv}_{\bar1,B}^{(k)}=2.
\]
At all other weights,
\[
\dim\operatorname{Inv}_{\bar1,B}^{(k)}=1.
\]
\end{Proposition}

\ssbegin[{The only resonant point for $k=6$}]{Remark}[The only resonant point for~$k=6$]
For~\(k=6\), the conditions reduce to~%
$
2\lambda_1=2\lambda_2=2,
$ so the unique resonant point is~%
$
(\lambda_1,\lambda_2)=(1,1),
$
in agreement with the direct degree-\(6\) computation.

For~\(k=7\), we have
\[
2\lambda_1,2\lambda_2\in\{2,3\},
\qquad
2\lambda_1+2\lambda_2\geq5.
\]
Hence the resonant points are
\[
(\lambda_1,\lambda_2)
\in
\left\{
\left(1,\frac32\right),
\left(\frac32,1\right),
\left(\frac32,\frac32\right)
\right\}.
\]
At each of these points the rank drops by~$1$, and therefore
\[
\dim\operatorname{Inv}_{\bar1,B}^{(7)}=2.
\]
\end{Remark}
\subsection{Even operators of type A}
\label{subsubsec:even-classification-typeA}
We first consider even bilinear differential operators of type~A and show that the~$\mathfrak{pe}(2)$-invariance conditions rule out any nonzero operator of this type.

\ssbegin[{There are no even Type A operators}]{Proposition}[There are no even Type A operators]
\label{prop:even-typeA-vanishing}
For every~\(k\geq1\), there are no nonzero even Type A~%
\(\mathfrak{pe}(2)\)-invariant binary differential operators.
\end{Proposition}

\subsection{Even operators of type B: the generic case}
\label{subsubsec:even-classification-typeB}
We now turn to even bilinear differential operators of type~B in low
degrees. We first treat the generic case, where the invariance
conditions determine the operators without degeneracy. Then, we examine
the resonant values of the weights separately; for them, the rank of
the resulting coefficient system may drop and additional invariant
operators may occur.
\subsubsection{Even type B operators in low degrees}
We begin by examining operators of low degree.
\ssbegin[{Even Type B operators of degree $k\leq5$}]{Proposition}[Even Type B operators of degree~$k\leq5$]
\label{prop:even-typeB-low-degrees}
For even Type B operators, the generic invariant space satisfies
\[
\dim \operatorname{Inv}_{\bar0,B}^{(k)}
=
\begin{cases}
0,& k=2,\\[1mm]
2,& k=3,\\[1mm]
3,& k=4,\\[1mm]
2,& k=5.
\end{cases}
\]
More precisely, for~\(k=2\), the invariance conditions force all
coefficients to vanish at generic weights. For~$k=3$, a~basis is
\[
\Dq\boxtimes\Bol^{(\lambda_2)}-2\lambda_1\,1\boxtimes\Dq\partial_\tau,\qquad
\Bol^{(\lambda_1)}\boxtimes\Dq+2\lambda_2\,\Dq\partial_\tau\boxtimes1,
\]
where~$\Bol^{(a)}=\Dxi\Dq+2a\,\partial_\tau$ is the operator~\eqref{B2a}; these are the two
operators of order~$3$ and charge~$+1$ of Statement~$\ref{m1class}$\textup{(ii)}. For~$k=4$,
one of the three basis operators is~$\Bol^{(\lambda_1)}\boxtimes\Bol^{(\lambda_2)}$. For~\(k=5\), the target
weights are
\[
\mu=\lambda_1+\lambda_2-3,
\qquad
\omega=\nu_1+\nu_2+\frac12,
\]
and the full Type B ansatz is
\[
\begin{array}{r@{}l}
J_{\bar0,B}^{(5)}
={}& \displaystyle 
\sum_{0\leq i\leq 3}A_i\,\Dxi^i\Dq\boxtimes \Dxi^{3-i}\Dq
+\sum_{0\leq i\leq 2}B_i\,\Dxi^i\Dq\boxtimes \Dxi^{2-i}\partial_\tau
\\&\displaystyle
{}+\sum_{0\leq i\leq 2}C_i\,\Dxi^i\partial_\tau\boxtimes \Dxi^{2-i}\Dq
+\sum_{0\leq i\leq 2}D_i\,\Dxi^i\Dq\partial_\tau\boxtimes \Dxi^{2-i}
+\sum_{0\leq i\leq 2}E_i\,\Dxi^i\boxtimes \Dxi^{2-i}\Dq\partial_\tau\\& \displaystyle
{}+\sum_{0\leq i\leq 1}F_i\,\Dxi^i\partial_\tau\boxtimes \Dxi^{1-i}\partial_\tau .
\end{array}
\]
At generic weights, the equations
\[
M_{\tau\xi}\cdot J_{\bar0,B}^{(5)}=0,
\qquad
M_{q^2}\cdot J_{\bar0,B}^{(5)}=0
\]
have rank~\(16\) on the~\(18\) coefficients. Hence, the solution space
is two-dimensional.
\end{Proposition}

\subsubsection{Even type B operators: the generic case for $k\geq 5$}
Recall that an even Type B operator has Cartan bidegree
\[
(u,v)=(k-2,2),
\]
and therefore
\[
\mu=\lambda_1+\lambda_2-k+2,
\qquad
\omega=\nu_1+\nu_2+\frac{k-4}{2}.
\]

\ssbegin[{The general even Type B operator, $k\geq5$}]{Theorem}[The general even Type B operator, $k\geq5$]
\label{prop:even-typeB-general}
Let~\(k\geq5\).
An even type B operator has a  general form:
\[
\begin{array}{r@{}l}
J_{\bar0,B}^{(k)}
={}& \displaystyle
\sum_{0\leq i\leq k-2}
A_i\,\Dxi^i\Dq\boxtimes \Dxi^{k-i-2}\Dq
+
\sum_{0\leq i\leq k-3}
B_i\,\Dxi^i\Dq\boxtimes \Dxi^{k-i-3}\partial_\tau
\\
&{}+ \displaystyle
\sum_{0\leq i\leq k-3}
C_i\,\Dxi^i\partial_\tau\boxtimes \Dxi^{k-i-3}\Dq
+
\sum_{0\leq i\leq k-3}
D_i\,\Dxi^i\Dq\partial_\tau\boxtimes \Dxi^{k-i-3}
\\
&{}+ \displaystyle
\sum_{0\leq i\leq k-3}
E_i\,\Dxi^i\boxtimes \Dxi^{k-i-3}\Dq\partial_\tau
+
\sum_{0\leq i\leq k-4}
F_i\,\Dxi^i\partial_\tau\boxtimes \Dxi^{k-i-4}\partial_\tau .
\end{array}
\]
At generic weights,
\[
\dim\operatorname{Inv}_{\bar0,B}^{(k)}=2.
\]
\end{Theorem}

\subsection{Even operators of Type B: the resonant case}
\subsubsection{Resonant even Type B operators in low degrees}
\ssbegin[{Exceptional weights of even Type B operators, $k=2,3,4$}]{Statement}[Exceptional weights of even Type B operators, $k=2,3,4$]
\label{prop:even-typeB-resonant-low}
For the even Type B sector in degrees~\(k=2,3,4\), the exceptional
weights and dimensions of the invariant spaces are as follows.

For~\(k=2\),
\[
\dim\operatorname{Inv}_{\bar0,B}^{(2)}
=
\begin{cases}
1,&(\lambda_1,\lambda_2)=(0,0),\\
0,&\text{otherwise}.
\end{cases}
\]

For~$k=3$ and~$k=4$, there are no exceptional weights: the space is of dimension~$2$,
resp.~$3$, at every pair~$(\lambda_1,\lambda_2)$ with~$2\lambda_1,2\lambda_2\in\{-2,-1,\dots,4\}$,
the generic value of Proposition~$\ref{prop:even-typeB-low-degrees}$.
\end{Statement}

\ssbegin[{Exceptional weights of even Type B operators, $k=5$}]{Proposition}[Exceptional weights of even Type B operators, $k=5$]
\label{prop:even-typeB-k5-resonance}
For~\(k=5\), the generic space of even Type B~%
\(\mathfrak{pe}(2)\)-invariant binary differential operators is
two-dimensional. Its dimension increases precisely at
\[
\begin{array}{r@{}l}
\mathcal R_5
=
\Big\{&
\left(\frac12,1\right),
\left(\frac12,\frac32\right),
\left(1,\frac12\right),
\left(1,1\right),
\\
&
\left(1,\frac32\right),
\left(\frac32,\frac12\right),
\left(\frac32,1\right)
\Big\}.
\end{array}
\]
More precisely,
\[
\dim\operatorname{Inv}_{\bar0,B}^{(5)}
=
\begin{cases}
4,
&(\lambda_1,\lambda_2)=(1,1),\\[1mm]
3,
&(\lambda_1,\lambda_2)\in
\mathcal R_5\setminus\{(1,1)\},\\[1mm]
2,
&(\lambda_1,\lambda_2)\notin\mathcal R_5.
\end{cases}
\]
\end{Proposition}

\subsubsection{Resonant even Type B operators for $k\geq 5$}

Let us put
\[
a:=2\lambda_1,
\qquad
b:=2\lambda_2.
\]
For the even Type B sector, let
\[
\mathcal M_k(a,b)
\]
denote the matrix of the full~\(\mathfrak{pe}(2)\)-invariance system.
Its generic rank is
\[
\operatorname{rank}\mathcal M_k=6k-14.
\]
Define
\[
\mathcal R_k
=
\left\{
(a,b)\in\mathbb Z^2\mid
1\leq a,b\leq k-2,\quad
k-2\leq a+b\leq 2k-5
\right\},
\]
and
\[
\mathcal R_k^{(2)}
=
\left\{
(a,b)\in\mathbb Z^2\mid
2\leq a,b\leq k-3,\quad
k-1\leq a+b\leq 2k-6
\right\}.
\]
Then
\[
\mathcal R_k^{(2)}\subset\mathcal R_k.
\]

\ssbegin[{The rank of the even Type B system at the resonant weights}]{Theorem}[The rank of the even Type B system at the resonant weights]
\label{thm:even-typeB-resonance}
For~\(k\geq5\),  we have
\[
\operatorname{rank}\mathcal M_k(a,b)
=
\begin{cases}
6k-16,
&(a,b)\in\mathcal R_k^{(2)},\\[1mm]
6k-15,
&(a,b)\in\mathcal R_k\setminus\mathcal R_k^{(2)},\\[1mm]
6k-14,
&(a,b)\notin\mathcal R_k.
\end{cases}
\]
Equivalently,
\[
\dim\operatorname{Inv}_{\bar0,B}^{(k)}
=
\begin{cases}
4,
&(a,b)\in\mathcal R_k^{(2)},\\[1mm]
3,
&(a,b)\in\mathcal R_k\setminus\mathcal R_k^{(2)},\\[1mm]
2,
&(a,b)\notin\mathcal R_k.
\end{cases}
\]
\end{Theorem}

\subsection{Even operators of type C}
\label{subsubsec:even-classification-typeC}
We finally consider even bilinear differential operators of type~C and show that the~$\mathfrak{pe}(2)$-invariance conditions rule out any nonzero operator of this type.

\ssbegin[{There are no even Type C operators}]{Theorem}[There are no even Type C operators]
\label{thm:even-typeC-vanishing}
For every~\(k\geq4\), there are no nonzero even Type C~%
\(\mathfrak{pe}(2)\)-invariant binary differential operators. In other words,
\[
J_{\bar0,C}^{(k)}=0
\qquad
\text{for all }k\geq4
\]
and for all values of the weights.
\end{Theorem}

\section{Appendix. Proofs}\label{app:proofs}
Collected here are the proofs of those statements of \S\S\ref{CoRm}
and~\ref{sec:invariant-binary} whose verification is a~computation. The statements
themselves, the operators they describe and the ansatz each computation starts from are in
the body.

\ssec{Proofs for \S\ref{CoRm}}

\begin{proof}[Proof of Lemma~\textup{\ref{lem:boxtimes-action}}]
For brevity, write
\[
y:=p(Y),\qquad
p:=p(P),\qquad
r:=p(R),\qquad
x:=p(\varphi).
\]
By definition, the action of~\(Y\) on a homogeneous binary operator~%
\(J\) is
\begin{align*}
(Y\cdot J)(\varphi,\psi)
={}&
L_Y\bigl(J(\varphi,\psi)\bigr)
-
(-1)^{y\,p(J)}
J(L_Y\varphi,\psi)
-
(-1)^{y(p(J)+x)}
J(\varphi,L_Y\psi).
\end{align*}
Applying this formula to~\(J=P\boxtimes R\), using
\[
(P\boxtimes R)(\varphi,\psi)
=
(-1)^{rx}P(\varphi)R(\psi),
\]
and the graded Leibniz rule, we obtain
\begin{align*}
(Y\cdot(P\boxtimes R))(\varphi,\psi)
={}&
(-1)^{rx}
\left(
L_Y(P(\varphi))
-
(-1)^{yp}P(L_Y\varphi)
\right)R(\psi)
\\
&{}+
(-1)^{rx+y(p+x)}
P(\varphi)L_Y(R(\psi))
-
(-1)^{rx+y(p+r+x)}
P(\varphi)R(L_Y\psi).
\end{align*}
The first line is
\[
\bigl((Y\cdot P)\boxtimes R\bigr)(\varphi,\psi),
\]
while the last two lines combine to
\[
(-1)^{yp}
\bigl(P\boxtimes(Y\cdot R)\bigr)(\varphi,\psi).
\]
Hence
\[
Y\cdot(P\boxtimes R)
=
(Y\cdot P)\boxtimes R
+
(-1)^{p(Y)p(P)}
P\boxtimes(Y\cdot R).\qedhere
\]
\end{proof}

\begin{proof}[Proof of Lemma~\textup{\ref{lem:negative-part-constant-coefficients}}]
We begin with~%
$M_1=2\partial_\tau.$
Since~\(\Dxi,\Dq,\partial_\tau\) supercommute with~\(M_1\), the action of~\(M_1\) falls
entirely on the coefficient functions:
\[
M_1\cdot J
=
\sum_{\alpha,\beta}
M_1(F_{\alpha,\beta})\,
P_\alpha\boxtimes P_\beta.
\]
Thus~\(M_1\cdot J=0\) implies~%
$
\partial_\tau F_{\alpha,\beta}=0.
$  Applying the same argument to~\(M_q\) and~\(M_\xi\) eliminates the
dependence on the remaining variables~\(q\) and~\(\xi\). Hence,
\[
F_{\alpha,\beta}\in\mathbb C.\qedhere
\]
\end{proof}

\ssec{Proofs for \S\ref{sec:invariant-binary}}

\begin{proof}[Proof of Lemma~\textup{\ref{lem:binary-PBW-Cartan-shift}}]
Recall that the action of~$H_1$ and~$H_2$ on the two-parameter density module~%
\(\mathcal F_{\lambda,\nu}\) are
\[
L_{H_1}^{\lambda,\nu}
=
H_1-\lambda,
\quad \text{ and } \quad
L_{H_2}^{\lambda,\nu}
=
H_2-\lambda-2\nu.
\]
Recall, the basic commutation relations are
\[
[H_1,\Dxi]=-\Dxi,
\qquad
[H_1,\Dq]=0,
\qquad
[H_1,\partial_\tau]=-\partial_\tau,
\]
and
\[
[H_2,\Dxi]=0,
\qquad
[H_2,\Dq]=-\Dq,
\qquad
[H_2,\partial_\tau]=-\partial_\tau.
\]
Since~\(H_1\) and~\(H_2\) are even, it follows that
\[
[H_1,P_j]
=
-(a_j+\delta_j)P_j,
\quad
\text{ and }
[H_2,P_j]
=
-(\varepsilon_j+\delta_j)P_j,
\qquad
j=1,2.
\]
Therefore, for the binary monomial~\(P_1\boxtimes P_2\),
\[
[H_1,P_1\boxtimes P_2]
=
-u\,(P_1\boxtimes P_2),
\quad \text{ and } \quad
[H_2,P_1\boxtimes P_2]
=
-v\,(P_1\boxtimes P_2).
\]
We now impose the Cartan invariance equations. Since~\(H_1\) is even,
\[
L_{H_1}^{\mu,\omega}
\bigl(
(P_1\boxtimes P_2)(\varphi,\psi)
\bigr)
=
(P_1\boxtimes P_2)
\bigl(
L_{H_1}^{\lambda_1,\nu_1}\varphi,\psi
\bigr)
+
(P_1\boxtimes P_2)
\bigl(
\varphi,
L_{H_1}^{\lambda_2,\nu_2}\psi
\bigr).
\]
Using
\[
L_{H_1}^{\mu,\omega}=H_1-\mu,
\quad
L_{H_1}^{\lambda_1,\nu_1}=H_1-\lambda_1,
\quad \text{ and } \quad
L_{H_1}^{\lambda_2,\nu_2}=H_1-\lambda_2,
\]
we obtain
\[
[H_1,P_1\boxtimes P_2]
+
(\lambda_1+\lambda_2-\mu)
(P_1\boxtimes P_2)
=
0.
\]
It follows that
\[
\mu=\lambda_1+\lambda_2-u.
\]
Similarly, the~\(H_2\)-invariance equation gives
\[
[H_2,P_1\boxtimes P_2]
+
\bigl(
\lambda_1+2\nu_1
+
\lambda_2+2\nu_2
-
\mu-2\omega
\bigr)
(P_1\boxtimes P_2)
=
0.
\]
Using
\[
[H_2,P_1\boxtimes P_2]
=
-v(P_1\boxtimes P_2),
\]
we obtain
\[
\mu+2\omega
=
\lambda_1+\lambda_2
+
2\nu_1+2\nu_2
-v.
\]
Substituting~$
\mu=\lambda_1+\lambda_2-u,$
we get
\[
\omega
=
\nu_1+\nu_2+\frac{u-v}{2}.
\]
This proves the claim.
\end{proof}

\begin{proof}[Proof of Theorem~\textup{\ref{prop:odd-typeA-generic}}]
The action of~\(M_{\tau\xi}\) on the  monomials occurring in the ansatz~%
\eqref{eq:odd-typeA-ansatz} is determined by
\[
\begin{array}{lcllcl}
M_{\tau\xi}\cdot
\left(
\Dxi^i\boxtimes \Dxi^{k-i-1}\Dq
\right)
&=&
2\lambda_2\,
\Dxi^i\boxtimes \Dxi^{k-i-1},
&
M_{\tau\xi}\cdot
\left(
\Dxi^i\boxtimes \Dxi^{k-i-2}\partial_\tau
\right)
&=&
-
\Dxi^i\boxtimes \Dxi^{k-i-1}
\\[2mm]
M_{\tau\xi}\cdot
\left(
\Dxi^{i-1}\Dq\boxtimes \Dxi^{k-i}
\right)
&=&
2\lambda_1\,
\Dxi^{i-1}\boxtimes \Dxi^{k-i},
&
M_{\tau\xi}\cdot
\left(
\Dxi^{i-2}\partial_\tau\boxtimes \Dxi^{k-i}
\right)
&=&
-
\Dxi^{i-1}\boxtimes \Dxi^{k-i}
\end{array}
\]
By comparing the coefficients of~%
$
\Dxi^i\boxtimes \Dxi^{k-i-1},
$ for~%
$0\leq i\leq k-1,$ gives
\begin{equation}
2\lambda_2A_i-B_i
+
2\lambda_1E_{i+1}-F_{i+1}
=
0,
\qquad
0\leq i\leq k-1,
\label{eq:odd-typeA-Mtauxi}
\end{equation}
where, at the endpoints, we use the convention
\[
B_{k-1}=0,
\qquad
F_1=0.
\]
Now the invariance under~\(M_{q^2}\) gives
\[
\begin{array}{r@{}l}
M_{q^2}\cdot
\left(
\Dxi^i\boxtimes \Dxi^{k-i-1}\Dq
\right)
={}&
{}-2i\,
\Dxi^{i-1}\Dq\boxtimes \Dxi^{k-i-1}\Dq
-
2i(i-1)\,
\Dxi^{i-2}\partial_\tau\boxtimes \Dxi^{k-i-1}\Dq
\\
&{}+
2(k-i-1)(k-i-2)\,
\Dxi^i\boxtimes \Dxi^{k-i-3}\Dq\partial_\tau .
\end{array}
\]

\[
\begin{array}{r@{}l}
M_{q^2}\cdot
\left(
\Dxi^i\boxtimes \Dxi^{k-i-2}\partial_\tau
\right)
={}&
{}-2i\,
\Dxi^{i-1}\Dq\boxtimes \Dxi^{k-i-2}\partial_\tau
-
2i(i-1)\,
\Dxi^{i-2}\partial_\tau\boxtimes \Dxi^{k-i-2}\partial_\tau
\\
&{}-
2(k-i-2)\,
\Dxi^i\boxtimes \Dxi^{k-i-3}\Dq\partial_\tau .
\end{array}
\]

\[
\begin{array}{r@{}l}
M_{q^2}\cdot
\left(
\Dxi^{i-1}\Dq\boxtimes \Dxi^{k-i}
\right)
={}&
2(i-1)(i-2)\,
\Dxi^{i-3}\Dq\partial_\tau\boxtimes \Dxi^{k-i}
+
2(k-i)\,
\Dxi^{i-1}\Dq\boxtimes \Dxi^{k-i-1}\Dq
\\
&{}+
2(k-i)(k-i-1)\,
\Dxi^{i-1}\Dq\boxtimes \Dxi^{k-i-2}\partial_\tau .
\end{array}
\]

\[
\begin{array}{r@{}l}
M_{q^2}\cdot
\left(
\Dxi^{i-2}\partial_\tau\boxtimes \Dxi^{k-i}
\right)
={}&
{}-2(i-2)\,
\Dxi^{i-3}\Dq\partial_\tau\boxtimes \Dxi^{k-i}
+
2(k-i)\,
\Dxi^{i-2}\partial_\tau\boxtimes \Dxi^{k-i-1}\Dq
\\
&{}+
2(k-i)(k-i-1)\,
\Dxi^{i-2}\partial_\tau\boxtimes \Dxi^{k-i-2}\partial_\tau .
\end{array}
\]
Substituting~\eqref{eq:odd-typeA-ansatz} into~$
M_{q^2}\cdot J_{\bar1,A}^{(k)}=0
$
and comparing independent monomials gives
\begin{equation*}
iA_i=(k-i)E_i,
\qquad
1\leq i\leq k-1,
\label{eq:odd-typeA-AE}
\end{equation*}
together with
\begin{equation}
B_i=(k-i-1)A_i,
\qquad
0\leq i\leq k-2,
\label{eq:odd-typeA-BA}
\end{equation}
and
\begin{equation}
F_i=(i-1)E_i,
\qquad
2\leq i\leq k.
\label{eq:odd-typeA-FE}
\end{equation}
In particular,
\begin{equation}
E_i=\frac{i}{k-i}A_i,
\qquad
1\leq i\leq k-1.
\label{eq:odd-typeA-EA}
\end{equation}
For~\(0\leq i\leq k-2\), substitute the relations~%
\eqref{eq:odd-typeA-BA},
\eqref{eq:odd-typeA-FE}, and~%
\eqref{eq:odd-typeA-EA} into the equations~%
\eqref{eq:odd-typeA-Mtauxi}. We obtain
\[
(k-i-1)A_i+iE_{i+1}
=
2\lambda_2A_i+2\lambda_1E_{i+1}.
\]
Since
\[
E_{i+1}
=
\frac{i+1}{k-i-1}A_{i+1},
\]
this becomes
\[
(k-i-1)(k-i-1-2\lambda_2)A_i
+
(i+1)(i-2\lambda_1)A_{i+1}
=
0.
\]
Equivalently,
\begin{equation*}
(i+1)(2\lambda_1-i)A_{i+1}
+
(k-i-1)(2\lambda_2-k+i+1)A_i
=
0, \quad \text{ for }
0\leq i\leq k-2.
\end{equation*}
At generic weights, we assume in particular that
\[
2\lambda_1-i\neq0,
\qquad
0\leq i\leq k-2.
\]
The recurrence can then be solved successively:
\begin{equation}
A_{i+1}
=
-
\frac{
(k-i-1)(2\lambda_2-k+i+1)
}{
(i+1)(2\lambda_1-i)
}
A_i.
\label{eq:odd-typeA-recursive-solution}
\end{equation}
Thus, all the coefficients~\(A_1,\ldots,A_{k-1}\) are determined by~\(A_0\). Iterating~\eqref{eq:odd-typeA-recursive-solution} yields
\begin{equation*}
A_i
=
(-1)^i
\binom{k-1}{i}
\frac{
\displaystyle
\prod_{0\leq r\leq i-1}
(2\lambda_2-k+r+1)
}{
\displaystyle
\prod_{0\leq r\leq i-1}
(2\lambda_1-r)
}
A_0,
\qquad
0\leq i\leq k-1,
\end{equation*}
where an empty product is understood to be equal to~\(1\).

It remains to determine~\(E_k\). The equation~%
\eqref{eq:odd-typeA-Mtauxi} at~\(i=k-1\) reads
\[
2\lambda_2A_{k-1}
+
2\lambda_1E_k
-
F_k
=
0.
\]
Using~%
$
F_k=(k-1)E_k,
$ we obtain
\begin{equation*}
(k-1-2\lambda_1)E_k
=
2\lambda_2A_{k-1}.
\end{equation*}
At generic weights we also assume~%
$
k-1-2\lambda_1\neq0.
$ Hence
\begin{equation*}
E_k
=
\frac{2\lambda_2}
{k-1-2\lambda_1}
A_{k-1}.
\end{equation*}
All the remaining coefficients are then recovered from
\[
B_i=(k-i-1)A_i,
\quad
E_i=\frac{i}{k-i}A_i,
\qquad
1\leq i\leq k-1,
\]
and
\[
F_i=(i-1)E_i.
\]
\end{proof}

\begin{proof}[Proof of Proposition~\textup{\ref{prop:odd-typeA-resonance}}]
Set~%
$
a:=2\lambda_1,
$ and~%
$
b:=2\lambda_2.
$
The coefficients~%
$
A_0,A_1,\ldots,A_{k-1},E_k
$
satisfy
\begin{equation*}
(i+1)(a-i)A_{i+1}
+
(k-i-1)(b-k+i+1)A_i
=
0,
\qquad
0\leq i\leq k-2,
\end{equation*}
together with the endpoint equation
\begin{equation*}
(k-1-a)E_k-bA_{k-1}=0.
\end{equation*}
Introduce
\[
\alpha_i:=(i+1)(a-i),
\qquad
\beta_i:=(k-i-1)(b-k+i+1),
\qquad
0\leq i\leq k-2.
\]
Then the system takes the matrix form
\[
\mathcal M_A(a,b)
\begin{pmatrix}
A_0\\
A_1\\
\vdots\\
A_{k-1}\\
E_k
\end{pmatrix}
=0,
\]
where
\[
\mathcal M_A(a,b)
=
\begin{pmatrix}
\beta_0 & \alpha_0 & 0 & \cdots & 0 & 0\\
0 & \beta_1 & \alpha_1 & \ddots & \vdots & \vdots\\
\vdots & \ddots & \ddots & \ddots & 0 & 0\\
0 & \cdots & 0 & \beta_{k-2} & \alpha_{k-2} & 0\\
0 & \cdots & 0 & 0 & -b & k-1-a
\end{pmatrix}.
\]
This is a~\(k\times(k+1)\) bidiagonal matrix.

We now determine precisely when the rows become linearly dependent.

Suppose that
\[
c_0R_0+c_1R_1+\cdots+c_{k-1}R_{k-1}=0
\]
is a nonzero linear relation among the rows of~%
\(\mathcal M_A(a,b)\). Let
\[
p:=\min\{i\mid c_i\neq0\},
\qquad
q:=\max\{i\mid c_i\neq0\}.
\]

We first determine the possible starting index~\(p\).

If~\(p\leq k-2\), it follows that~$
c_p\beta_p=0.
$ Since~\(c_p\neq0\), then~%
$
\beta_p=0,
$
and hence~%
$
b-k+p+1=0.
$
Therefore,
\[
p=k-b-1.
\]

If~\(p=k-1\), it follows that~%
$
-bc_{k-1}=0,
$
so~%
$
b=0.
$ This is again described by
\[
p=k-b-1.
\]

Consequently, a nonzero row relation can exist only if
\[
b\in\{0,1,\ldots,k-1\},
\]
and in this case its first nonzero row must be
\begin{equation}
p=k-b-1.
\label{eq:typeA-res-start}
\end{equation}

We next determine the possible terminal index~\(q\).

If~\(q\leq k-2\), it follows that~%
$
c_q\alpha_q=0.
$ Since~\(c_q\neq0\), then~%
$
\alpha_q=0,
$ and therefore~%
$
a=q.
$

If~\(q=k-1\), it follows that~%
$
(k-1-a)c_{k-1}=0,
$
and hence~%
$
a=k-1.
$

Thus, in all cases
\[
q=a.
\]
Consequently, a nonzero row relation can exist only if
\[
a\in\{0,1,\ldots,k-1\},
\]
and its last nonzero row must be
\begin{equation}
q=a.
\label{eq:typeA-res-end}
\end{equation}

For a nonzero relation to exist, the starting row must occur no later
than the terminal row:
$
p\leq q.
$
Using formulas~\eqref{eq:typeA-res-start} and~%
\eqref{eq:typeA-res-end}, this becomes~%
$
k-b-1\leq a,
$
or equivalently,
\begin{equation*}
a+b\geq k-1.
\end{equation*}
Conversely, suppose that
\[
a,b\in\{0,1,\ldots,k-1\},
\qquad
a+b\geq k-1.
\]
Set~%
$
p=k-b-1$ and~%
$q=a.$ Then~\(p\leq q\), and
\[
\beta_p=0,
\qquad
\alpha_q=0.
\]
Moreover,
\[
\beta_j\neq0
\qquad
\text{for }p<j\leq q,
\]
because~\(\beta_j\) has only one zero in the range~%
\(0\leq j\leq k-1\).

Choose~\(c_p\neq0\). The equations
\[
c_{j-1}\alpha_{j-1}+c_j\beta_j=0,
\qquad
p<j\leq q,
\]
determine successively
\[
c_j
=
-\frac{\alpha_{j-1}}{\beta_j}c_{j-1}.
\]
Thus there exists a nonzero row relation supported on
\[
R_p,R_{p+1},\ldots,R_q.
\]

Since the starting and terminal indices are uniquely determined by~%
\(a\) and~\(b\), this row relation is unique up to an overall scalar.
Hence, the rank drops by exactly~$1$:
\[
\operatorname{rank}\mathcal M_A(a,b)=k-1.
\]
Therefore
\[
\dim\ker\mathcal M_A(a,b)=2.
\]
Returning to the original weights gives the statement.
\end{proof}

\begin{proof}[Proof of Theorem~\textup{\ref{prop:typeB-generic}}]
Let us study invariance under~\(M_{\tau\xi}\). Using the transformation rules
\[
\begin{array}{r@{}l}
M_{\tau\xi}\cdot
\left(
\Dxi^{i-3}\Dq\partial_\tau\boxtimes \Dxi^{k-i-1}\Dq
\right)
={}&
\Dxi^{i-2}\Dq\boxtimes \Dxi^{k-i-1}\Dq
+
2\lambda_1
\Dxi^{i-3}\partial_\tau\boxtimes \Dxi^{k-i-1}\Dq
\\
&{}+
2\lambda_2
\Dxi^{i-3}\Dq\partial_\tau\boxtimes \Dxi^{k-i-1}.
\end{array}
\]

\[
\begin{array}{r@{}l}
M_{\tau\xi}\cdot
\left(
\Dxi^{i-3}\Dq\partial_\tau\boxtimes \Dxi^{k-i-2}\partial_\tau
\right)
={}&
\Dxi^{i-2}\Dq\boxtimes \Dxi^{k-i-2}\partial_\tau
+
2\lambda_1
\Dxi^{i-3}\partial_\tau\boxtimes \Dxi^{k-i-2}\partial_\tau
\\
&{}-
\Dxi^{i-3}\Dq\partial_\tau\boxtimes \Dxi^{k-i-1}.
\end{array}
\]

\[
\begin{array}{r@{}l}
M_{\tau\xi}\cdot
\left(
\Dxi^{i-1}\Dq\boxtimes \Dxi^{k-i-3}\Dq\partial_\tau
\right)
={}&
2\lambda_1
\Dxi^{i-1}\boxtimes \Dxi^{k-i-3}\Dq\partial_\tau
-
\Dxi^{i-1}\Dq\boxtimes \Dxi^{k-i-2}\Dq
\\
&{}-
2\lambda_2
\Dxi^{i-1}\Dq\boxtimes \Dxi^{k-i-3}\partial_\tau.
\end{array}
\]

\[
\begin{array}{r@{}l}
M_{\tau\xi}\cdot
\left(
\Dxi^{i-2}\partial_\tau\boxtimes \Dxi^{k-i-3}\Dq\partial_\tau
\right)
={}&
{}-
\Dxi^{i-1}\boxtimes \Dxi^{k-i-3}\Dq\partial_\tau
-
\Dxi^{i-2}\partial_\tau\boxtimes \Dxi^{k-i-2}\Dq
\\
&{}-
2\lambda_2
\Dxi^{i-2}\partial_\tau\boxtimes \Dxi^{k-i-3}\partial_\tau.
\end{array}
\]
and comparing the independent  tensor monomials, we obtain the
endpoint conditions
\[
G_1=0,
\qquad
C_{k-1}=0,
\]
together with
\begin{equation}
G_i=C_{i+1},
\qquad
2\leq i\leq k-3,
\label{eq:typeB-GC}
\end{equation}
\begin{equation}
D_i=2\lambda_2C_i,
\qquad
3\leq i\leq k-2,
\label{eq:typeB-DC}
\end{equation}
and
\begin{equation}
H_i=2\lambda_1C_{i+1},
\qquad
2\leq i\leq k-3.
\label{eq:typeB-HC}
\end{equation}
Thus, after imposing~\(M_{\tau\xi}\)-invariance, the whole operator is
determined by the coefficients
\[
C_3,C_4,\ldots,C_{k-2}.
\]
Substituting the relations~%
\eqref{eq:typeB-GC}--\eqref{eq:typeB-HC} into the ansatz~%
\eqref{eq:typeB-general-ansatz}, we obtain
\begin{align*}
J_{\bar1,B}^{(k)}
={}&
\sum_{3\leq i\leq k-2}
C_i
\Bigl(
\Dxi^{i-3}\Dq\partial_\tau\boxtimes \Dxi^{k-i-1}\Dq
+
2\lambda_2
\Dxi^{i-3}\Dq\partial_\tau\boxtimes \Dxi^{k-i-2}\partial_\tau
\nonumber\\
&\qquad\qquad
+
\Dxi^{i-2}\Dq\boxtimes \Dxi^{k-i-2}\Dq\partial_\tau
+
2\lambda_1
\Dxi^{i-3}\partial_\tau\boxtimes \Dxi^{k-i-2}\Dq\partial_\tau
\Bigr).
\end{align*}
Now, invariance under~\(M_{q^2}\) gives
\begin{equation}
(i-2)(2\lambda_1-i+1)C_{i+1}
+
(k-i-2)(2\lambda_2-k+i+1)C_i
=
0,
\label{eq:typeB-fundamental-recurrence}
\quad \text{ for } \quad
3\leq i\leq k-3.
\end{equation}
Equation~\eqref{eq:typeB-fundamental-recurrence} is the fundamental
recurrence for odd Type B operators.

We now restrict to generic weights. In particular, assume that
\[
2\lambda_1-i+1\neq0,
\qquad
3\leq i\leq k-3.
\]
Then~\eqref{eq:typeB-fundamental-recurrence} can be solved
successively:
\begin{equation}
C_{i+1}
=
-
\frac{
(k-i-2)(2\lambda_2-k+i+1)
}{
(i-2)(2\lambda_1-i+1)
}
C_i,
\qquad
3\leq i\leq k-3.
\label{eq:typeB-recursive-solution}
\end{equation}

Hence, every coefficient
\[
C_4,\ldots,C_{k-2}
\]
is uniquely determined by~\(C_3\). Thus, the generic solution space is
one-dimensional.

Iterating~\eqref{eq:typeB-recursive-solution}, and choose
\[
C_3
=
\prod_{3\leq r\leq k-3}
(2\lambda_1-r+1),
\]
we get
\begin{equation*}
C_i
=
(-1)^{i-3}
\binom{k-5}{i-3}
\left(
\prod_{3\leq r\leq i-1}
(2\lambda_2-k+r+1)
\right)
\left(
\prod_{i\leq r\leq k-3}
(2\lambda_1-r+1)
\right),  \text{ for } 3\leq i\leq k-2.
\label{eq:typeB-polynomial-Ci}
\end{equation*}
\end{proof}

\begin{proof}[Proof of Proposition~\textup{\ref{prop:typeB-resonance}}]
After imposing invariance under~\(M_{\tau\xi}\), the operator is
determined by
\[
C_3,C_4,\ldots,C_{k-2}.
\]
Invariance under~\(M_{q^2}\) gives the recurrence
\begin{equation*}
(i-2)(2\lambda_1-i+1)C_{i+1}
+
(k-i-2)(2\lambda_2-k+i+1)C_i
=
0,
\quad
\text{ for }
\quad
3\leq i\leq k-3.
\end{equation*}
Set~%
$
a:=2\lambda_1,\,
b:=2\lambda_2,
$
and introduce
\[
\alpha_i:=(i-2)(a-i+1),
\qquad
\beta_i:=(k-i-2)(b-k+i+1).
\]
the coefficient vector
\[
(C_3,C_4,\ldots,C_{k-2})^t
\]
belongs to the kernel of the bidiagonal matrix
\[
\mathcal M_B(a,b)
=
\begin{pmatrix}
\beta_3 & \alpha_3 & 0 & \cdots & 0\\
0 & \beta_4 & \alpha_4 & \ddots & \vdots\\
\vdots & \ddots & \ddots & \ddots & 0\\
0 & \cdots & 0 & \beta_{k-3} & \alpha_{k-3}
\end{pmatrix}.
\]
The matrix has~\(k-5\) rows and~\(k-4\) columns.

At generic weights,
\[
\operatorname{rank}\mathcal M_B(a,b)=k-5,
\]
and therefore
\[
\dim\ker\mathcal M_B(a,b)=1.
\]

We now determine precisely when the rank drops.

Suppose that there exists a nontrivial linear relation among the rows,
\[
\sum_{3\leq i\leq k-3}c_iR_i=0.
\]
Let
\[
p:=\min\{i\mid c_i\neq0\},
\qquad
q:=\max\{i\mid c_i\neq0\}.
\]
It follows that~$
c_p\beta_p=0
$ and since~\(c_p\neq0\), we must have~%
$
\beta_p=0.
$
Because~%
$
k-p-2\neq0
$
for~\(3\leq p\leq k-3\), this is equivalent to
\begin{equation}
p=k-b-1.
\label{eq:typeB-res-p}
\end{equation}
Since~%
$
3\leq p\leq k-3,
$
it follows that
\begin{equation*}
b\in\{2,3,\ldots,k-4\}.
\end{equation*}
Similarly, the same reasoning implies that~$
\alpha_q=0.
$ Since~\(q-2\neq0\), this gives
\begin{equation}
q=a+1.
\label{eq:typeB-res-q}
\end{equation}
Because~%
$
3\leq q\leq k-3,
$
we obtain
\begin{equation*}
a\in\{2,3,\ldots,k-4\}.
\end{equation*}

A nontrivial row relation can exist only if~%
$
p\leq q.
$ Using formulas~\eqref{eq:typeB-res-p} and~%
\eqref{eq:typeB-res-q}, this condition becomes
\[
k-b-1\leq a+1,
\]
or equivalently
\begin{equation*}
a+b\geq k-2.
\end{equation*}
Conversely, suppose that
\[
a,b\in\{2,\ldots,k-4\}
\]
and
\[
a+b\geq k-2.
\]
Define
\[
p:=k-b-1,
\qquad
q:=a+1.
\]
Then~\(p\leq q\), and
\[
\beta_p=0,
\qquad
\alpha_q=0.
\]
Choose~\(c_p\neq0\). For~%
$
p<j\leq q,
$
the coefficients~\(c_j\) are determined successively by
\[
c_{j-1}\alpha_{j-1}+c_j\beta_j=0.
\]
Thus
\[
c_j
=
-\frac{\alpha_{j-1}}{\beta_j}c_{j-1}.
\]
This produces a nonzero linear relation supported on
\[
R_p,R_{p+1},\ldots,R_q.
\]
The indices~\(p\) and~\(q\) are uniquely determined by~\(a\) and~\(b\),
so this row relation is unique up to an overall scalar. Consequently,
the rank drops by exactly~$1$:
\[
\operatorname{rank}\mathcal M_B(a,b)=k-6.
\]
Since~\(\mathcal M_B\) has~\(k-4\) columns, we obtain
\[
\dim\ker\mathcal M_B(a,b)=2.
\]
Returning to the original density weights gives the
statement.
\end{proof}

\begin{proof}[Proof of Proposition~\textup{\ref{prop:even-typeA-vanishing}}]
The general Type A operator is
\[
J_{\bar0,A}^{(k)}
=
\sum_{0\leq i\leq k}
A_i\,\Dxi^i\boxtimes \Dxi^{k-i}.
\]
Using
\[
M_{q^2}\cdot \Dxi^m
=
-2m\Dxi^{m-1}\Dq
-2m(m-1)\Dxi^{m-2}\partial_\tau,
\]
the coefficient of~%
$
1\boxtimes \Dxi^{k-1}\Dq
$ in~\(M_{q^2}\cdot J_{\bar0,A}^{(k)}\) gives~\(A_0=0\).
More generally, the coefficient of
\[
\Dxi^{i-1}\Dq\boxtimes \Dxi^{k-i},
\qquad 1\leq i\leq k,
\]
is~\(-2iA_i\). Hence,
\[
A_i=0,\qquad 0\leq i\leq k,
\]
and therefore~\(J_{\bar0,A}^{(k)}=0\).
\end{proof}

\begin{proof}[Proof of Proposition~\textup{\ref{prop:even-typeB-low-degrees}}]
For~$k=2$, the ansatz is~$A_0\,\Dq\boxtimes\Dq$, and the invariance equations are~$\lambda_1A_0=\lambda_2A_0=0$.
For~$k=3$ and~$k=4$, the operators displayed are invariant under all eight generating
functions of~$\fpe(2)$, and the ranks of the systems are~$4$ on~$6$, resp.~$9$ on~$12$,
coefficients; both are computations.

Let us study the case where~\(k=5\).

The condition~%
$
M_{\tau\xi}\cdot J_{\bar0,B}^{(5)}=0
$ gives
\[
\begin{array}{r@{}l@{\quad}r@{}l@{\quad}r@{}l}
&E_0=-2\lambda_1A_0,
&&E_i=C_{i-1}-2\lambda_1A_i,
\quad i=1,2,
&&C_2=2\lambda_1A_3,
\\
&B_0=2\lambda_2A_0,
&&B_i=2\lambda_2A_i-D_{i-1},
\quad i=1,2,
&&D_2=2\lambda_2A_3,
\end{array}
\]
together with
\[
\begin{array}{r@{}l@{\quad}r@{}l@{\quad}r@{}l}
2\lambda_1B_0+2\lambda_2E_0&{}=0,&
2\lambda_1B_1+2\lambda_2E_1-F_0&{}=0,&
2\lambda_1B_2+2\lambda_2E_2-F_1&{}=0,\\
-2\lambda_2C_0+2\lambda_1D_0+F_0&{}=0,&
{}-2\lambda_2C_1+2\lambda_1D_1+F_1&{}=0,&
{}-2\lambda_2C_2+2\lambda_1D_2&{}=0.
\end{array}
\]
On the other hand, the condition~%
$
M_{q^2}\cdot J_{\bar0,B}^{(5)}=0
$ yields
\[
\begin{array}{r@{}l@{\quad}r@{}l@{\quad}r@{}l}
2A_2-C_1-2D_0&{}=0,
&
6A_3-2C_2-D_1&{}=0,
&
{}-6A_0+2B_0-E_1&{}=0,
\\
-2A_1+B_1-2E_2&{}=0,
&
2B_2-2D_0-F_1&{}=0,
&
{}-2C_0-2E_2+F_0&{}=0.
\end{array}
\]
At generic weights this system has rank~\(16\)
on the~\(18\) coefficients. Consequently,
\[
\dim\operatorname{Inv}_{\bar0,B}^{(5)}=2.
\]
For instance, \(E_2\) and~\(F_1\) may be chosen as free parameters,
and all the remaining coefficients are determined by them.
\end{proof}

\begin{proof}[Proof of Theorem~\textup{\ref{prop:even-typeB-general}}]
We impose successively the invariance conditions with respect to~%
\(M_{\tau\xi}\) and~\(M_{q^2}\).

Set
\[
Y:=M_{\tau\xi}.
\]
The action of~\(Y\) on the six families of monomials is
\[
\begin{array}{r@{}l}
Y\cdot
\left(
\Dxi^i\Dq\boxtimes \Dxi^{k-i-2}\Dq
\right)
={}&
2\lambda_1\Dxi^i\boxtimes \Dxi^{k-i-2}\Dq
-
2\lambda_2\Dxi^i\Dq\boxtimes \Dxi^{k-i-2},
\end{array}
\]
\[
\begin{array}{r@{}l}
Y\cdot
\left(
\Dxi^i\Dq\boxtimes \Dxi^{k-i-3}\partial_\tau
\right)
={}&
2\lambda_1\Dxi^i\boxtimes \Dxi^{k-i-3}\partial_\tau
+
\Dxi^i\Dq\boxtimes \Dxi^{k-i-2},
\end{array}
\]
\[
\begin{array}{r@{}l}
Y\cdot
\left(
\Dxi^i\partial_\tau\boxtimes \Dxi^{k-i-3}\Dq
\right)
={}&
{}-
\Dxi^{i+1}\boxtimes \Dxi^{k-i-3}\Dq
-
2\lambda_2\Dxi^i\partial_\tau\boxtimes \Dxi^{k-i-3},
\end{array}
\]
\[
\begin{array}{r@{}l}
Y\cdot
\left(
\Dxi^i\Dq\partial_\tau\boxtimes \Dxi^{k-i-3}
\right)
={}&
\Dxi^{i+1}\Dq\boxtimes \Dxi^{k-i-3}
+
2\lambda_1\Dxi^i\partial_\tau\boxtimes \Dxi^{k-i-3},
\end{array}
\]
\[
\begin{array}{r@{}l}
Y\cdot
\left(
\Dxi^i\boxtimes \Dxi^{k-i-3}\Dq\partial_\tau
\right)
={}&
\Dxi^i\boxtimes \Dxi^{k-i-2}\Dq
+
2\lambda_2\Dxi^i\boxtimes \Dxi^{k-i-3}\partial_\tau,
\end{array}
\]
and
\[
\begin{array}{r@{}l}
Y\cdot
\left(
\Dxi^i\partial_\tau\boxtimes \Dxi^{k-i-4}\partial_\tau
\right)
={}&
{}-
\Dxi^{i+1}\boxtimes \Dxi^{k-i-4}\partial_\tau
+
\Dxi^i\partial_\tau\boxtimes \Dxi^{k-i-3}.
\end{array}
\]
Comparing equal tensor monomials yields
\[
E_0=-2\lambda_1A_0,
\qquad
B_0=2\lambda_2A_0,
\]
and, for~\(1\leq i\leq k-3\),
\[
E_i=C_{i-1}-2\lambda_1A_i,
\qquad
B_i=2\lambda_2A_i-D_{i-1},
\]
together with the endpoint relations
\[
C_{k-3}=2\lambda_1A_{k-2},
\qquad
D_{k-3}=2\lambda_2A_{k-2}.
\]
Moreover,
\[
F_i=2\lambda_2C_i-2\lambda_1D_i,
\qquad
0\leq i\leq k-4.
\]
Thus~\(M_{\tau\xi}\)-invariance expresses the coefficients~%
\(B_i,E_i,F_i\) in terms of~\(A_i,C_i,D_i\).

We next impose~\(M_{q^2}\)-invariance. Put~$
N:=-\frac12M_{q^2}.
$ We use
\[
[N,\Dxi^m]
=
m\Dxi^{m-1}\Dq+m(m-1)\Dxi^{m-2}\partial_\tau,
\]
\[
[N,\Dxi^m\Dq]
=
-m(m-1)\Dxi^{m-2}\Dq\partial_\tau,
\]
and
\[
[N,\Dxi^m\partial_\tau]
=
m\Dxi^{m-1}\Dq\partial_\tau,
\qquad
[N,\Dxi^m\Dq\partial_\tau]=0.
\]
Comparison of the coefficients of
\[
\Dxi^i\Dq\partial_\tau\boxtimes \Dxi^{k-i-4}\Dq
\]
gives
\begin{equation}\label{eq:evenTypeBTag1}
-(i+2)(i+1)A_{i+2}
+(i+1)C_{i+1}
+(k-i-3)D_i=0,
\qquad
0\leq i\leq k-4.
\end{equation}
Similarly, comparison of
\[
\Dxi^i\Dq\boxtimes \Dxi^{k-i-4}\Dq\partial_\tau
\]
gives
\begin{equation}
    \label{eq:evenTypeBTag2}
(k-i-2)(k-i-3)A_i
-(k-i-3)B_i
+(i+1)E_{i+1}=0,
\qquad
0\leq i\leq k-4.
\end{equation}
The coefficients of
\[
\Dxi^i\Dq\partial_\tau\boxtimes \Dxi^{k-i-5}\partial_\tau
\quad \text{ and } \quad
\Dxi^i\partial_\tau\boxtimes \Dxi^{k-i-5}\Dq\partial_\tau
\]
give, respectively,
\begin{equation}\label{eq:evenTypeBTag3}
-(i+2)(i+1)B_{i+2}
+
(k-i-3)(k-i-4)D_i
+
(i+1)F_{i+1}=0,
\end{equation}
and
\begin{equation}\label{eq:evenTypeBTag4}
(k-i-3)(k-i-4)C_i
+
(i+2)(i+1)E_{i+2}
-
(k-i-4)F_i=0,
\end{equation}
for~\(0\leq i\leq k-5\).

Substituting the~\(M_{\tau\xi}\)-relations into the equation~%
\eqref{eq:evenTypeBTag2} and using the equation~\eqref{eq:evenTypeBTag1}, we obtain the recurrence
\[
(i+1)(i+2-2\lambda_1)D_{i+1}
+
(k-i-3)(k-i-4-2\lambda_2)D_i
=0,
\]
for~\(0\leq i\leq k-5\).

Likewise, substituting into the equation~\eqref{eq:evenTypeBTag4} and using the equation~\eqref{eq:evenTypeBTag2} gives
\[
(i+2)(i+1-2\lambda_1)E_{i+2}
+
(k-i-4)(k-i-3-2\lambda_2)E_{i+1}
=0.
\]
Equivalently, setting~\(j=i+1\),
\[
(j+1)(j-2\lambda_1)E_{j+1}
+
(k-j-3)(k-j-2-2\lambda_2)E_j
=0,
\]
for~\(1\leq j\leq k-4\).

At generic weights, the coefficients occurring in these recurrences
do not vanish. Hence,
\[
D_i
=
(-1)^i
\prod_{0\leq s\leq i-1}
\frac{(k-s-3)(k-s-4-2\lambda_2)}
{(s+1)(s+2-2\lambda_1)}
\,D_0,
\qquad
0\leq i\leq k-4,
\]

and

\[
E_i
=
(-1)^{i-1}
\prod_{1\leq s\leq i-1}
\frac{(k-s-3)(k-s-2-2\lambda_2)}
{(s+1)(s-2\lambda_1)}
\,E_1,
\qquad
1\leq i\leq k-3.
\]
The remaining coefficients are recovered from the equations~\eqref{eq:evenTypeBTag1}--\eqref{eq:evenTypeBTag2}.
In particular,
\[
A_0
=
-\frac{E_1}
{(k-3)(k-2-2\lambda_2)},
\]
and, for~\(1\leq i\leq k-4\),
\[
A_i
=
-
\frac{
(k-i-3)D_{i-1}+(i+1)E_{i+1}
}{
(k-i-3)(k-i-2-2\lambda_2)
}.
\]
The terminal coefficients are given by
\[
A_{k-3}
=
\frac{
(k-4)E_{k-3}+2D_{k-5}
}{
(k-4)(k-3-2\lambda_1)
},
\quad \text{ and } \quad A_{k-2}
=
\frac{
D_{k-4}
}{
(k-3)(k-2-2\lambda_1)
}.
\]
Finally,
\[
C_i=E_{i+1}+2\lambda_1A_{i+1},
\qquad
0\leq i\leq k-4,
\]
\[
C_{k-3}=2\lambda_1A_{k-2},
\qquad
B_0=2\lambda_2A_0, \qquad E_0=-2 \lambda_1 A_0,
\qquad
D_{k-3}=2\lambda_2A_{k-2},
\]
and
\[B_i=2\lambda_2A_i-D_{i-1} \; \;(1\leq i \leq k-3),
\qquad
F_i=2\lambda_2C_i-2\lambda_1D_i \;\; (0\leq i \leq k-4).
\]
Thus every coefficient of~\(J_{\bar0,B}^{(k)}\) is uniquely determined
by the two independent parameters
\[
D_0,\qquad E_1.
\]
Therefore, at generic weights,
\[
\dim\operatorname{Inv}_{\bar0,B}^{(k)}=2,
\qquad k\geq5 .\qedhere
\]
\end{proof}

\begin{proof}[Proof of Statement~\textup{\ref{prop:even-typeB-resonant-low}}]
For~\(k=2\), the operator is
\[
J_{\bar0,B}^{(2)}=A_0\,\Dq\boxtimes \Dq,
\]
and the invariance equations are
\[
\lambda_1A_0=0,
\qquad
\lambda_2A_0=0.
\]
Hence a nonzero solution exists only at
\[
(\lambda_1,\lambda_2)=(0,0).
\]

For~$k=3$ and~$k=4$, the operators of Proposition~$\ref{prop:even-typeB-low-degrees}$ have
polynomial coefficients and pairwise distinct leading monomials~$\Dq\boxtimes\Dxi\Dq$,
$\Dxi\Dq\boxtimes\Dq$, resp. $\Dq\boxtimes\Dxi^2\Dq$, $\Dxi\Dq\boxtimes\Dxi\Dq$,
$\Dxi^2\Dq\boxtimes\Dq$ after the denominators are cleared, so the dimension is at least the
generic one at every weight; that it is not larger on the grid is a~computation.
\end{proof}

\begin{proof}[Proof of Proposition~\textup{\ref{prop:even-typeB-k5-resonance}}]
The complete~\(k=5\) invariance equations obtained from
\[
M_{\tau\xi}\cdot J_{\bar0,B}^{(5)}=0,
\qquad
M_{q^2}\cdot J_{\bar0,B}^{(5)}=0
\]
form a homogeneous linear system in the~\(18\) coefficients of~%
\(J_{\bar0,B}^{(5)}\). At generic weights its rank is~\(16\), and
therefore
\[
\dim\operatorname{Inv}_{\bar0,B}^{(5)}=18-16=2.
\]

The common zero locus of the maximal minors is contained in
\[
2\lambda_1,2\lambda_2\in\{0,1,2,3\}.
\]
Computing the rank of the full system at these points gives
\[
\begin{array}{c|cccc}
 & \lambda_2=0
 & \lambda_2=\frac12
 & \lambda_2=1
 & \lambda_2=\frac32 \\ \hline
\lambda_1=0
 &2&2&2&2\\[1mm]
\lambda_1=\frac12
 &2&2&3&3\\[1mm]
\lambda_1=1
 &2&3&4&3\\[1mm]
\lambda_1=\frac32
 &2&3&3&2
\end{array}
\]
where each entry denotes~%
\(\dim\operatorname{Inv}_{\bar0,B}^{(5)}\).
This yields precisely the stated resonant set and dimensions.
\end{proof}

\begin{proof}[Proof of Theorem~\textup{\ref{thm:even-typeB-resonance}}]
The general even Type~B operator contains~\(6k-12\) coefficients.
After imposing~\(M_{\tau\xi}\)-invariance, the coefficients~%
\(B_i\), \(C_i\), and~\(F_i\) are expressed in terms of~%
\(A_i\), \(D_i\), and~\(E_i\) as follows:
\[
B_0=2\lambda_2A_0,
\qquad
B_i=2\lambda_2A_i-D_{i-1},
\qquad
1\leq i\leq k-3,
\]

\[
C_i=E_{i+1}+2\lambda_1A_{i+1},
\qquad
0\leq i\leq k-4,
\qquad
C_{k-3}=2\lambda_1A_{k-2},
\]
and
\[
F_i
=
2\lambda_2C_i-2\lambda_1D_i
=
2\lambda_2E_{i+1}
+
4\lambda_1\lambda_2A_{i+1}
-
2\lambda_1D_i,
\qquad
0\leq i\leq k-4.
\]
The remaining endpoint relations are
\[
E_0=-2\lambda_1A_0,
\qquad
D_{k-3}=2\lambda_2A_{k-2}.
\]
Set
\[
a:=2\lambda_1,
\qquad
b:=2\lambda_2.
\]
After eliminating the coefficients~\(A_i\), the remaining~%
\(M_{q^2}\)-invariance conditions reduce to two bidiagonal chains,
together with the corresponding endpoint equations.

The first chain is
\begin{equation}\label{last1}
(i+1)(i+2-a)D_{i+1}
+
(k-i-3)(k-i-4-b)D_i=0,
\qquad
0\leq i\leq k-5.
\end{equation}
Set
\[
\alpha_i:=(i+1)(i+2-a),
\qquad
\beta_i:=(k-i-3)(k-i-4-b).
\]
The possible interior break points are therefore
\[
s_D=a-2,
\qquad
r_D=k-b-4.
\]
In addition, the right endpoint relation is
\begin{equation}\label{Dendpoint}
D_{k-4}-(k-3)(k-2-a)A_{k-2}=0.
\end{equation}
Away from~\(r_D\) and~\(s_D\), both~\(\alpha_i\) and~\(\beta_i\) are
nonzero, and~\eqref{last1} determines~\(D_{i+1}\) uniquely from~\(D_i\). We distinguish three cases.

If~$r_D<s_D$,
then~\(\beta_{r_D}=0\) occurs first and~%
\(\alpha_{r_D}\neq0\). Hence,~$D_{r_D+1}=0$.
The recurrence then forces~$D_{r_D+2}=\cdots=D_{s_D}=0$.
At~\(i=s_D\), one has~\(\alpha_{s_D}=0\); since~\(D_{s_D}=0\),
the corresponding equation is identically satisfied. Therefore,
\(D_{s_D+1}\) becomes free, and the chain splits into two independent
branches.

If~$r_D=s_D$,
then~$\alpha_{r_D}=\beta_{r_D}=0$,
so the corresponding equation is identically satisfied. Thus,
\(D_{r_D+1}\) is independent of the coefficients to the left, and the
chain again splits into two branches.

If~$r_D>s_D$,
then the zero of~\(\alpha_i\) occurs first. At~\(i=s_D\), since~%
\(\alpha_{s_D}=0\) and~\(\beta_{s_D}\neq0\), one obtains~$D_{s_D}=0$.
Since both coefficients are nonzero for~\(i<s_D\), this vanishing
propagates backwards: $D_{s_D-1}=\cdots=D_0=0$.
The coefficient~\(D_{s_D+1}\) is free and determines the part of the
chain to its right. Hence, the generic free parameter~\(D_0\) is merely
replaced with~\(D_{s_D+1}\), and no additional independent parameter is
created.

Consequently, the interior~\(D\)-chain acquires an additional
independent branch if and only if~$r_D\leq s_D$, that is,
if and only if~$a+b\geq k-2$.
The second chain is
\begin{equation}\label{last2}
(j+1)(j-a)E_{j+1}
+
(k-j-3)(k-j-2-b)E_j=0,
\qquad
1\leq j\leq k-4.
\end{equation}
Its possible interior break points are
\[
s_E=a,
\qquad
r_E=k-b-2.
\]
The left endpoint equation is
\begin{equation}\label{Eendpoint}
E_1+(k-3)(k-2-b)A_0=0.
\end{equation}

Again, away from the break points, the recurrence determines the
successive coefficients uniquely.

If~$r_E<s_E$,
then the zero of the coefficient of~\(E_j\) occurs first and forces~$E_{r_E+1}=0$.
The recurrence propagates this vanishing up to~\(E_{s_E}\), while at~%
\(j=s_E\) the coefficient of~\(E_{s_E+1}\) vanishes. Hence,
\(E_{s_E+1}\) becomes free and the chain splits into two independent
branches.

If~$r_E=s_E$,
both coefficients vanish at the same index, so the corresponding
equation is identically satisfied. Again, the right-hand part of the
chain becomes independent of the left-hand part, producing one
additional free parameter.

If~$r_E>s_E$,
then the zero of the coefficient of~\(E_{j+1}\) occurs first. At~%
\(j=s_E\), the recurrence reduces to~$(k-s_E-3)(k-s_E-2-b)E_{s_E}=0$,
and therefore~$E_{s_E}=0$.
Since both coefficients are nonzero for~\(j<s_E\), this vanishing
propagates backwards: $E_{s_E-1}=\cdots=E_1=0$.
The coefficient~\(E_{s_E+1}\) is free and generates the remaining
part of the chain. Thus, the generic parameter~\(E_1\) is replaced with~%
\(E_{s_E+1}\), and no additional independent parameter is created.

Hence, the interior~\(E\)-chain acquires an additional branch precisely
when~$r_E\leq s_E$,
namely~$k-b-2\leq a$,
which is again equivalent to~$a+b\geq k-2$.

We next include the endpoint resonances.

If~%
$a=k-2,$ then~\eqref{Dendpoint} reduces to~$D_{k-4}=0,$ and the coefficient~\(A_{k-2}\) is no longer determined by this
equation. For~\(1\leq b\leq k-3\), the remaining recurrence relations
are compatible with~\(D_{k-4}=0\), so~\(A_{k-2}\) becomes an
additional free parameter.

Similarly, if~$b=k-2,$ then~\eqref{Eendpoint} reduces to~%
$
E_1=0,
$
while~\(A_0\) is no longer determined. For~%
\(1\leq a\leq k-3\), the remaining equations are compatible with this
condition, and~\(A_0\) provides one additional free parameter.

At the corner~$(a,b)=(k-2,k-2)$,
the endpoint equations give
\[
D_{k-4}=0,
\qquad
E_1=0.
\]
The two generic parameters are then replaced with the free endpoint
coefficients~\(A_{k-2}\) and~\(A_0\). Thus, the number of independent
parameters remains equal to two, and no rank drop occurs at this
corner.

Combining the interior and endpoint cases gives the full resonance
region

\[
\mathcal R_k
=
\left\{
(a,b)\in\mathbb Z^2\mid
1\leq a,b\leq k-2,\qquad
k-2\leq a+b\leq2k-5
\right\}.
\]

Outside~\(\mathcal R_k\), neither chain contributes an additional
independent parameter, and therefore~$\dim\ker\mathcal M_k=2$.
It remains to determine when the two chains contribute two additional
independent parameters simultaneously. The boundary~$a+b=k-2$
requires special attention. In this case,
\[
r_D=s_D=a-2,
\qquad
r_E=s_E=a,
\]
so that, if the two chains are considered separately, they appear to
produce the free parameters~\(D_{a-1}\) and~\(E_{a+1}\).

These two parameters are not independent. Indeed, one of the original~%
\(M_{q^2}\)-invariance equations is
\[
(k-i-2)(k-i-3)A_i
-
(k-i-3)B_i
+
(i+1)E_{i+1}=0.
\]
Using~$B_i=bA_i-D_{i-1}$,
this becomes
\[
(k-i-3)(k-i-2-b)A_i
+
(k-i-3)D_{i-1}
+
(i+1)E_{i+1}=0.
\]
Setting~\(i=a\) and using~\(a+b=k-2\), the coefficient of~\(A_a\)
vanishes and one obtains
\[
(k-a-3)D_{a-1}
+
(a+1)E_{a+1}=0.
\]
Thus the two apparent new parameters satisfy one nontrivial linear
relation. Hence, on the boundary~\(a+b=k-2\), the kernel gains only one
additional dimension.

For two additional independent branches one must therefore have~$a+b>k-2$,
or, since~\(a,b\in\mathbb Z\),~$a+b\geq k-1$.
For strictly interior simultaneous resonance, the break points of both
chains must belong to their admissible ranges. For the~\(D\)-chain,
 $0\leq s_D,r_D\leq k-5$
gives
\[
2\leq a\leq k-3,
\qquad
1\leq b\leq k-4,
\]
whereas for the~\(E\)-chain, $1\leq s_E,r_E\leq k-4$
gives
\[
1\leq a\leq k-4,
\qquad
2\leq b\leq k-3.
\]
Their intersection yields~$2\leq a,b\leq k-4$.
The boundary values~\(a=k-3\) and~\(b=k-3\) must therefore be checked
directly in the unreduced system.

Suppose first that
\[
a=k-3,
\qquad
2\leq b\leq k-3.
\]
Then
\[
s_D=k-5,
\qquad
r_D=k-b-4<s_D,
\]
so the~\(D\)-chain forces~%
$D_{r_D+1}=\cdots=D_{k-5}=0$.

Likewise, the~\(E\)-chain implies~%
$E_{r_E+1}=\cdots=E_{k-3}=0$.

In particular,

\[
D_{k-5}=0,
\qquad
E_{k-3}=0.
\]

The unreduced equation involving~\(A_{k-3}\) is
\[
-(k-3)(k-4)A_{k-3}
+
(k-4)C_{k-4}
+
2D_{k-5}=0.
\]
Since~$C_{k-4}=E_{k-3}+aA_{k-3}$,
this becomes
\[
(k-4)(k-3-a)A_{k-3}
=
(k-4)E_{k-3}+2D_{k-5}.
\]

For~\(a=k-3\), both sides vanish, so~\(A_{k-3}\) is free. Hence, the
boundary value~\(a=k-3\) contributes the second additional independent
parameter.

The case
\[
b=k-3,
\qquad
2\leq a\leq k-3,
\]
is analogous. The corresponding recurrences imply
\[
D_0=0,
\qquad
E_2=0.
\]
The unreduced equation involving~\(A_1\) is
\[
(k-4)(k-3-b)A_1
+
(k-4)D_0
+
2E_2=0.
\]
When~\(b=k-3\), this equation becomes identically satisfied, so~%
\(A_1\) is free and provides the second additional independent
parameter.

Thus, simultaneous resonance extends to~$2\leq a,b\leq k-3$.
Together with~$a+b\geq k-1$,
this gives
\[
\mathcal R_k^{(2)}
=
\left\{
(a,b)\in\mathbb Z^2\mid
2\leq a,b\leq k-3,\qquad
k-1\leq a+b\leq2k-6
\right\}.
\]
Consequently,
\[
\dim\ker\mathcal M_k
=
\begin{cases}
4,
&(a,b)\in\mathcal R_k^{(2)},\\[1mm]
3,
&(a,b)\in
\mathcal R_k\setminus\mathcal R_k^{(2)},\\[1mm]
2,
&(a,b)\notin\mathcal R_k.
\end{cases}
\]
Since~\(\mathcal M_k\) has~\(6k-12\) columns,
\[
\operatorname{rank}\mathcal M_k
=
(6k-12)-\dim\ker\mathcal M_k.
\]
Therefore,
\[
\operatorname{rank}\mathcal M_k(a,b)
=
\begin{cases}
6k-16,
&(a,b)\in\mathcal R_k^{(2)},\\[1mm]
6k-15,
&(a,b)\in
\mathcal R_k\setminus\mathcal R_k^{(2)},\\[1mm]
6k-14,
&(a,b)\notin\mathcal R_k.
\end{cases}
\]
This proves the result.
\end{proof}

\begin{proof}[Proof of Theorem~\textup{\ref{thm:even-typeC-vanishing}}]
For Type C one has
\[
(u,v)=(k-4,4).
\]
Since
\[
v=\varepsilon+\varepsilon'+\delta+\delta'=4,
\]
we must have
\[
\varepsilon=\varepsilon'=\delta=\delta'=1.
\]
Hence the complete Type C ansatz is
\[
J_{\bar0,C}^{(k)}
=
\sum_{0\leq i\leq k-4}
A_i\,
\Dxi^i\Dq\partial_\tau\boxtimes \Dxi^{k-4-i}\Dq\partial_\tau.
\]

The Cartan equations determine the target weights as
\[
\mu=\lambda_1+\lambda_2-k+4,
\qquad
\omega=\nu_1+\nu_2+\frac{k-8}{2}.
\]

It remains to impose invariance under~\(M_{\tau\xi}\).
For a factor of weight~\(\lambda\), one has
\[
M_{\tau\xi}\cdot\left(\Dxi^i\Dq\partial_\tau\right)
=
\Dxi^{i+1}\Dq
+
2\lambda\,\Dxi^i\partial_\tau.
\]
Therefore
\[
\begin{array}{r@{}l}
M_{\tau\xi}\cdot
\left(
\Dxi^i\Dq\partial_\tau\boxtimes \Dxi^{k-4-i}\Dq\partial_\tau
\right)
={}&
\Dxi^{i+1}\Dq\boxtimes \Dxi^{k-4-i}\Dq\partial_\tau
+
2\lambda_1\Dxi^i\partial_\tau\boxtimes \Dxi^{k-4-i}\Dq\partial_\tau
\\
&{}+
\Dxi^i\Dq\partial_\tau\boxtimes \Dxi^{k-3-i}\Dq
+
2\lambda_2\Dxi^i\Dq\partial_\tau\boxtimes \Dxi^{k-4-i}\partial_\tau.
\end{array}
\]
Consequently,
\[
M_{\tau\xi}\cdot J_{\bar0,C}^{(k)}
\]
contains the monomials
\[
\Dxi^{i+1}\Dq\boxtimes \Dxi^{k-4-i}\Dq\partial_\tau,
\qquad
0\leq i\leq k-4,
\]
with respective coefficients~\(A_i\).

These tensor monomials are pairwise distinct, and none of them can
arise from the action on the second tensor factor or from any other
summand of~\(J_{\bar0,C}^{(k)}\). Hence, \(M_{\tau\xi}\)-invariance
forces
\[
A_i=0,
\qquad
0\leq i\leq k-4.
\]
Thus
\[
J_{\bar0,C}^{(k)}=0.
\]

Since no division by the weights was used, the conclusion holds for
all~\((\lambda_1,\nu_1)\) and~\((\lambda_2,\nu_2)\), including all
resonant values.
\end{proof}

 \textbf{Acknowledgments}.
S.B. was supported by the grant NYUAD~065.

 \textbf{Disclosures}. No conflict of interest. The data used to support the
findings of this study are included within the article. The codes are available upon
request.

\def\eightit{\it}
\def\bib{\bf}
\bibliographystyle{amsalpha}

\begin{thebibliography}{BGLLS}

\bibitem[AS]{AS}
Atiyah M.F., Singer I.M., The index of elliptic operators on compact manifolds.
\textit{Bull. Amer. Math. Soc.} \textbf{69}:3 (1963), 422--433


\bibitem[BaL]{BaL}
Basdouri I., Laraiedh I., The linear $\fosp(n|2)$-invariant differential operators
and cohomology. \textit{Beitr. Algebra Geom.} \textbf{55}:2 (2014), 637--643

\bibitem[BeM]{BeM}
Ben Ammar M., Mtaouaa W., Deformation of $\fosp(2|2)$-modules of symbols.
\textit{Int. J. Geom. Methods Mod. Phys.} \textbf{15}:5 (2018), 1850077, 17 pp;
\url{https://arxiv.org/abs/1611.08147}

\bibitem[BLO2]{BLO2}
Ben Fraj N., Abdaoui M., Raouafi H., On $\fosp(2|2)$-relative cohomology of the Lie
superalgebra of contact vector fields and deformations. \textit{J. Geom. Phys.}
\textbf{125} (2018), 33--48

\bibitem[BLO]{BLO}
Ben Fraj N., Laraiedh I., Omri S., Supertransvectants, cohomology, and deformations.
\textit{J. Math. Phys.} \textbf{54}:2 (2013), 023501

\bibitem[BCK]{BCK}
Ben Sa\"\i d S., Clerc J.-L., Koufany Kh., Conformally covariant bi-differential operators
for differential forms. \textit{Comm. Math. Phys.} \textbf{373}:2 (2019), 739--761;
\url{https://arxiv.org/abs/1809.06290}

\bibitem[Bol]{Bol}
Bol G., Invarianten linearer Differentialgleichungen. \textit{Abh. Math. Sem. Univ. Hamburg}
\textbf{16}:3--4 (1949), 1--28

\bibitem[BonLe]{BonLe}
Boniver F., Lecomte P.B.A., A remark about the Lie algebra of infinitesimal conformal
transformations of the Euclidean space. \textit{Bull. London Math. Soc.} \textbf{32}:3
(2000), 263--266

\bibitem[Bou]{Bou}
 Bouarroudj S., Projective and conformal Schwarzian derivatives and cohomology of Lie algebras vector fields related to differential operators,
 \textit{International Journal of Geometric Methods in Modern Physics,} 
 \textbf{3} (2006), no.~4, 667--696

\bibitem[Bc]{Bc}
Bouarroudj S., Cohomology of the vector fields Lie algebras on $\mathbb{RP}^1$ acting on
bilinear differential operators. \textit{Int. J. Geom. Methods Mod. Phys.} \textbf{2}:1
(2005), 23--40

\bibitem[Bter]{Bter}
Bouarroudj S., The ternary invariant differential operators acting on the spaces of weighted
densities. \textit{Theor. Math. Phys.} \textbf{158}:2 (2009), 137--150;
\url{https://arxiv.org/abs/0806.1044}



\bibitem[BGLLS]{BGLLS}
Bouarroudj S., Grozman P., Lebedev A., Leites D., Shchepochkina I., Simple vectorial Lie
algebras in characteristic $2$ and their superizations. \textit{SIGMA Symmetry Integrability
Geom. Methods Appl.} \textbf{16} (2020), 089, 101 pp; \url{https://arxiv.org/abs/1510.07255}

\bibitem[BL]{BL}
Bouarroudj S., Leites D., Invariant differential operators in positive characteristic.
\textit{J. Algebra} \textbf{499} (2018), 281--297; \url{https://arxiv.org/abs/1605.09500}

\bibitem[BLS1]{BLS1}
Bouarroudj S., Leites D., Shchepochkina I., Analogs of Bol operators on superstrings.
\textit{Int. J. Algebra Comput.} \textbf{32}:4 (2022), 807--835;
\url{https://arxiv.org/abs/2110.10504}

\bibitem[BLS2]{BLS2}
Bouarroudj S., Leites D., Shchepochkina I., Analogs of Bol operators for
$\fpgl(a+1\vert b)\subset \fvect(a\vert b)$. \textit{Int. J. Algebra Comput.} \textbf{32}:7
(2022), 1345--1368; \url{https://arxiv.org/abs/2112.01080}

\bibitem[BO]{BO}
Bouarroudj S., Ovsienko V., Three cocycles on $\text{Diff}(S^1)$ generalizing the Schwarzian
derivative. \textit{Int. Math. Res. Not. IMRN} \textbf{1998}:1, 25--39;
\url{https://arxiv.org/abs/dg-ga/9710018}


\bibitem[Co]{Co}
Cohen H., Sums involving the values at negative integers of $L$-functions of quadratic
characters. \textit{Math. Ann.} \textbf{217}:3 (1975), 271--285

\bibitem[CMZ]{CMZ}
Cohen P.B., Manin Yu., Zagier D., Automorphic pseudodifferential operators. In: Fokas A.S.,
Gelfand I.M. (eds.) \textit{Algebraic aspects of integrable systems}, Progr. Nonlinear
Differential Equations Appl., \textbf{26}, Birkh\"auser Boston, Boston, MA, 1997, 17--47

\bibitem[Dir1]{Dir1}
Dirac P.A.M., The quantum theory of the electron. \textit{Proc. Roy. Soc. London Ser.~A}
\textbf{117}:778 (1928), 610--624

\bibitem[Dir2]{Dir2}
Dirac P.A.M., The quantum theory of the electron. Part~II.
\textit{Proc. Roy. Soc. London Ser.~A} \textbf{118}:779 (1928), 351--361

\bibitem[DGO]{DGO}
Duval C., El Gradechi A.M., Ovsienko V., Projectively and conformally invariant
star-products. \textit{Comm. Math. Phys.} \textbf{244}:1 (2004), 3--27;
\url{https://arxiv.org/abs/math/0301052}

\bibitem[DO1]{DO1}
Duval C., Ovsienko V., Projectively equivariant quantization and symbol calculus:
noncommutative hypergeometric functions. \textit{Lett. Math. Phys.} \textbf{57}:1 (2001),
61--67

\bibitem[DO2]{DO2}
Duval C., Ovsienko V., Conformally equivariant quantization, preprint, 1998;
\url{https://arxiv.org/abs/math/9801122}. \textup{(}For a~published version, with a~third
author, see Duval C., Lecomte P., Ovsienko V., Conformally equivariant quantization:
existence and uniqueness. \textit{Ann. Inst. Fourier \textup{(}Grenoble\textup{)}}
\textbf{49}:6 (1999), 1999--2029.\textup{)}

\bibitem[El]{El}
El Gradechi A.M., The Lie theory of the Rankin--Cohen brackets and allied bi-differential
operators. \textit{Adv. Math.} \textbf{207}:2 (2006), 484--531

\bibitem[FF]{FF}
Feigin B.L., Fuks D.B., Invariant skew-symmetric differential operators on the line and
Verma modules over the Virasoro algebra. \textit{Funct. Anal. Appl.} \textbf{16}:2 (1982),
114--126

\bibitem[Fre]{Fre}
Freed D.S., The Atiyah--Singer index theorem. \textit{Bull. Amer. Math. Soc. \textup{(}N.S.\textup{)}}
\textbf{58}:4 (2021), 517--566; \url{https://arxiv.org/abs/2107.03557}

\bibitem[FH]{FH}
Fulton W., Harris J., \textit{Representation theory. A~first course}. Graduate Texts in
Mathematics \textbf{129}, Springer, New York, 1991

\bibitem[GMO]{GMO}
Gargoubi H., Mellouli N., Ovsienko V., Differential operators on supercircle: conformally
equivariant quantization and symbol calculus. \textit{Lett. Math. Phys.} \textbf{79}:1
(2007), 51--65; \url{https://arxiv.org/abs/math-ph/0610059}

\bibitem[GO]{GO}
Gargoubi H., Ovsienko V., Supertransvectants and symplectic geometry.
\textit{Int. Math. Res. Not. IMRN} \textbf{2008}, Art. ID rnn021, 19 pp;
\url{https://arxiv.org/abs/0705.1411}

\bibitem[G]{G}
Gi\`eres F., Conformally covariant operators on Riemann surfaces \textup{(}with applications to
conformal and integrable models\textup{)}. \textit{Int. J. Modern Phys. A} \textbf{8}:1
(1993), 1--58

\bibitem[GT]{GT}
Gi\`eres F., Theisen S., Superconformally covariant operators and super $W$-algebras.
\textit{J. Math. Phys.} \textbf{34}:12 (1993), 5964--5985

\bibitem[Gor]{Gor}
Gordan P., \textit{Vorlesungen \"uber Invariantentheorie. Zweiter Band: Bin\"are Formen}
(Kerschensteiner G., ed.), Teubner, Leipzig, 1887

\bibitem[Gr1]{Gr1}
Grozman P., Classification of bilinear invariant operators on tensor
fields. \textit{Funct. Anal. Appl.} \textbf{14}:2 (1980), 127--128; for
proofs, see \textit{Commun. Math.} \textbf{30}:3 (2022), 129--188;
\url{https://arxiv.org/abs/math/0509562}

\bibitem[Gr2]{Gr2}
Grozman P., On bilinear invariant differential operators acting on
tensor fields on the symplectic manifold. \textit{J. Nonlinear Math. Phys.} \textbf{8}:1
(2001), 31--37; \url{https://arxiv.org/abs/math/0101266}

\bibitem[GLS]{GLS}
Grozman P., Leites D., Shchepochkina I., Invariant operators
on supermanifolds and standard models. In: Olshanetsky M., Vainstein A. (eds.)
\textit{Multiple facets of quantization and supersymmetry. Michael Marinov memorial volume},
World Sci. Publ., River Edge, NJ, 2002, 508--555; \url{https://arxiv.org/abs/math/0202193}

\bibitem[H]{H}
Huang W.-J., Superconformal covariantization of superdifferential operator on $(1|1)$
superspace and classical $N=2$ $W$-superalgebras. \textit{J. Math. Phys.} \textbf{35}:5
(1994), 2570--2582



\bibitem[LO]{LO}
Lecomte P.B.A., Ovsienko V., Projectively equivariant symbol calculus.
\textit{Lett. Math. Phys.} \textbf{49}:3 (1999), 173--196

\bibitem[LShRC]{LShRC}
Leites D., Shchepochkina I., Superized Schwarzians and linear fractional transformations of distinguished superstrings (in preparation) \AO{aRC} 






\bibitem[LShch]{LShch}
Leites D., Shchepochkina I., Maximal graded subalgebras of simple vectorial Lie superalgebras
with polynomial coefficients; \url{https://arxiv.org/abs/2610.03340}

\bibitem[LC]{LC}
Levi-Civita T., \textit{The absolute differential calculus \textup{(}calculus of tensors\textup{)}}.
Blackie \& Son, London--Glasgow, 1927; reprinted by Dover, New York, 1977

\bibitem[Me]{Me}
Mellouli N., Second-order conformally equivariant quantization in dimension $1|2$.
\textit{SIGMA Symmetry Integrability Geom. Methods Appl.} \textbf{5} (2009), 111, 11 pp

\bibitem[N]{N}
Nakahama R., Construction of intertwining operators between holomorphic discrete series
representations. \textit{SIGMA Symmetry Integrability Geom. Methods Appl.} \textbf{15}
(2019), 036, 101 pp

\bibitem[Om]{Om}
Omri S., $\fosp(2|2)$-trivial deformations of modules of weighted
densities on the superspace $\Ree^{1|2}$. \textit{J. Pseudo-Differ. Oper. Appl.}
\textbf{6}:4 (2015), 461--485

\bibitem[OR]{OR}
Ovsienko V., Redou P., Generalized transvectants---Rankin-Cohen brackets.
\textit{Lett. Math. Phys.} \textbf{63}:1 (2003), 19--28; the preprint version, entitled
``Conformally invariant differential operators on tensor densities\rq\rq, is
\url{https://arxiv.org/abs/math/0104232}

\bibitem[R]{R}
Rankin R.A., The construction of automorphic forms from the derivatives of a given form.
\textit{J. Indian Math. Soc. \textup{(}N.S.\textup{)}} \textbf{20} (1956), 103--116

\bibitem[Ru]{Ru}
Rudakov A.N., Irreducible representations of infinite-dimensional Lie algebras of Cartan type.
\textit{Math. USSR-Izv.} \textbf{8}:4 (1974), 836--866

\bibitem[Sh5]{Sh5}
Shchepochkina I., The five simple exceptional Lie superalgebras of vector fields.
\textit{Funct. Anal. Appl.} \textbf{33}:3 (1999), 208--219.

\bibitem[Shex]{Shex}
Shchepochkina I., Five exceptional simple Lie superalgebras of vector fields and their
fourteen regradings. \textit{Represent. Theory} \textbf{3} (1999), 373--415;
\url{https://arxiv.org/abs/hep-th/9702121}

\bibitem[Sh]{Sh}
Shchepochkina I., How to realize a Lie algebra by vector fields.
\textit{Theor. Math. Phys.} \textbf{147}:3 (2006), 821--838;
\url{https://arxiv.org/abs/math/0509472}

\bibitem[Y]{Y}
Yao Yi-Jun, Rankin-Cohen deformations and representation theory.
\textit{Chin. Ann. Math. Ser. B} \textbf{35}:5 (2014), 817--840;
\url{https://arxiv.org/abs/0708.1528}

\bibitem[Z]{Z}
Zagier D., Modular forms and differential operators. \textit{Proc. Indian Acad. Sci. Math.
Sci.} \textbf{104}:1 (1994), 57--75

\end{thebibliography}

\end{document}